\documentclass[11pt, a4paper, twoside,leqno]{amsart}
\usepackage[centering, totalwidth = 420pt, totalheight = 610pt]{geometry}
\usepackage{amssymb, amsmath, amsthm, enumerate, microtype, stmaryrd, mathpartir, url,mathrsfs}
\usepackage[latin1]{inputenc}
\usepackage[dvips, arrow, matrix, tips, curve]{xy}
\usepackage{color}
\definecolor{darkgreen}{rgb}{0,0.45,0}
\usepackage[pagebackref,colorlinks,citecolor=darkgreen,linkcolor=darkgreen]{hyperref}
\SelectTips{cm}{10}

\newcommand{\twocong}[2][0.5]{\ar@{}[#2] \save ?(#1)*{\cong}\restore}
\newcommand{\twoeq}[2][0.5]{\ar@{}[#2] \save ?(#1)*{=}\restore}
\newcommand{\rtwocell}[3][0.5]{\ar@{}[#2] \ar@{=>}?(#1)+/l 0.2cm/;?(#1)+/r 0.2cm/^{#3}}
\newcommand{\ltwocell}[3][0.5]{\ar@{}[#2] \ar@{=>}?(#1)+/r 0.2cm/;?(#1)+/l 0.2cm/^{#3}}
\newcommand{\ltwocello}[3][0.5]{\ar@{}[#2] \ar@{=>}?(#1)+/r 0.2cm/;?(#1)+/l 0.2cm/_{#3}}
\newcommand{\dltwocell}[3][0.5]{\ar@{}[#2] \ar@{=>}?(#1)+/ur  0.2cm/;?(#1)+/dl 0.2cm/^{#3}}
\newcommand{\dtwocell}[3][0.5]{\ar@{}[#2] \ar@{=>}?(#1)+/u  0.2cm/;?(#1)+/d 0.2cm/^{#3}}
\newcommand{\dthreecell}[3][0.5]{\ar@{}[#2] \ar@3{->}?(#1)+/u  0.2cm/;?(#1)+/d 0.2cm/^{#3}}
\newcommand{\utwocell}[3][0.5]{\ar@{}[#2] \ar@{=>}?(#1)+/d 0.2cm/;?(#1)+/u 0.2cm/_{#3}}
\newcommand{\dtwocelltarg}[3][0.5]{\ar@{}#2 \ar@{=>}?(#1)+/u  0.2cm/;?(#1)+/d 0.2cm/^{#3}}
\newcommand{\utwocelltarg}[3][0.5]{\ar@{}#2 \ar@{=>}?(#1)+/d  0.2cm/;?(#1)+/u 0.2cm/_{#3}}

\DeclareMathOperator{\ob}{ob}

\newcommand{\cat}[1]{\mathbf{#1}}
\newcommand{\op}{\mathrm{op}}
\newcommand{\id}{\mathrm{id}}
\newcommand{\thg}{{\mathord{\text{--}}}}

\newcommand{\spn}[1]{{\left<{#1}\right>}}

\newcommand{\defn}[1]{\textit{#1}}
\newcommand{\cd}[2][]{\vcenter{\hbox{\xymatrix#1{#2}}}}
\newcommand{\cdl}[2][]{\xymatrix@1#1{#2}}
\newcommand{\defeq}{\mathrel{\mathop:}=}

\newcommand{\A}{{\mathcal A}}
\newcommand{\B}{{\mathcal B}}
\newcommand{\C}{{\mathcal C}}
\newcommand{\D}{{\mathcal D}}
\newcommand{\E}{{\mathcal E}}

\newcommand{\J}{{\mathcal J}}

\renewcommand{\L}{{\mathcal L}}

\newcommand{\R}{{\mathcal R}}
\newcommand{\T}{{\mathcal T}}

\renewcommand{\c}{,\,\,}
\renewcommand{\t}{\,\vdash\,}

\newtheorem{Thm}{Theorem}[subsection]
\newtheorem{Prop}[Thm]{Proposition}
\newtheorem{Cor}[Thm]{Corollary}
\newtheorem{Lem}[Thm]{Lemma}

\theoremstyle{definition}
\newtheorem{Defn}[Thm]{Definition}

\theoremstyle{remark}

\newtheorem{Rk}[Thm]{Remark}
\newtheorem{Not}[Thm]{Notation}

\newcommand{\coop}{\mathrm{coop}}
\newcommand{\ty}{\mathsf{type}}
\newcommand{\Id}{\mathsf{Id}}
\newcommand{\tcat}[1]{\mathbf{#1}}
\newcommand{\subst}{\mathsf{subst}}
\newcommand{\ext}{\mathsf{ext}}

\renewcommand{\r}{\mathrm{r}}
\renewcommand{\J}{\mathrm{J}}
\newcommand{\Ell}{\mathrm{E}}

\newcommand{\ctxt}{\mathsf{ctxt}}
\renewcommand{\vec}[1]{#1}
\newcommand{\Ss}{\mathrm S}
\newcommand{\con}{\mathcal{C}}
\newcommand{\typ}{\mathcal{T}}
\newcommand{\Con}{\mathfrak{C}}
\newcommand{\Typ}{\mathfrak{T}}
\newcommand{\Tt}{\mathrm T}

\newcommand{\mli}{\ensuremath{\mathrm{ML}_I}}
\newcommand{\mle}{\ensuremath{\mathrm{ML}_E}}
\newcommand{\mlt}{\ensuremath{\mathrm{ML}_2}}
\renewcommand{\sb}[2]{{#1}^\ast{#2}}
\newcommand{\ro}[1]{\iota}
\DeclareMathOperator{\cod}{cod}

\begin{document}

\title{Two-dimensional models of type theory}
\author{Richard Garner}
\thanks{Supported by a Research Fellowship of St John's College, Cambridge and a Marie Curie Intra-European Fellowship, Project No.\ 040802.}
\address{Department of Pure Mathematics and Mathematical Statistics, 
Wilberforce Road, Cambridge CB3 0WB, UK}
\email{rhgg2@cam.ac.uk}
\date{11 September 2008; Revised 21 January 2009} \maketitle
\begin{abstract}
We describe a non-extensional variant of Martin-L\"of type theory which we call
\emph{two-dimensional type theory}, and equip it with a sound and complete
semantics valued in $2$-categories.
\end{abstract}
\bibliographystyle{acm}

\section{Introduction}

This is the second in a series of papers detailing the author's investigations
into the intensional type theory of Martin-L\"of, as described
in~\cite{Nordstrom1990Programming}. The first of these
papers,~\cite{Garner2008strength}, investigated syntactic issues relating to
its dependent product types. The present paper is a contribution to its
categorical semantics.

In~\cite{Seely1984Locally} it was proposed that the correct categorical models
for extensional Martin-L\"of type theory should be locally cartesian closed
categories: these being categories $\C$ with finite limits in which each of the
functors \mbox{$f^\ast \colon \C / X \to \C / Y$} induced by pulling back along
a morphism $f\colon Y \to X$ has a right adjoint. The idea is to think of each
object $X$ of a locally cartesian closed category $\C$ as a closed type, each
morphism as a term, and each object of the slice category $\C / X$ as a type
dependent upon $X$. Now substitution of terms in types may be interpreted by
pullback between the slices of $\C$; dependent sum and product types by left
and right adjoints to pullback; and the equality type on $X$ by the diagonal
morphism $\Delta \colon X \to X \times X$ in $\C / X \times X$. It was later
pointed out in~\cite{Hofmann1995interpretation} that this picture, whilst very
appealing, is not wholly accurate, since in the syntax, the operation which to
each morphism of types $f \colon Y \to X$ assigns the corresponding
substitution operation $\cat{Type}(X) \to \cat{Type}(Y)$ is strictly functorial
in $f$; whilst in the semantics, the corresponding assignation $(f \colon Y \to
X) \mapsto (f^\ast \colon \C/X \to \C/Y)$ is rarely so. Thus this notion of
model is not \emph{sound} for the syntax, and we are forced to refine it
slightly: essentially by equipping our locally cartesian closed category with a
split fibration $\T \to \C$ equivalent to its codomain fibration $\C^\mathbf 2
\to \C$. Types over $X$ are now interpreted as objects of the fibre category
$\T(X)$; and since $\T \to \C$ is a split fibration, the interpretation is
sound for substitution.

The question of how the above should generalise from extensional to intensional
Martin-L\"of type theory is a delicate one. It is possible to paraphrase the
syntax of intensional type theory in categorical language and so arrive at a
notion of model---as done in~\cite{Dybjer1996Internal}
or~\cite{Hofmann1995Extensional}, for example---but then we lose sight of a
key aspect of the extensional semantics, namely that dependent sums and
products may be characterised universally, as left and right adjoints to
substitution. To obtain a similar result for the intensional theory requires a
more refined sort of semantics. More specifically, we are thinking of a
semantics valued in higher-dimensional categories, motivated by work such as
\cite{Awodey2008Homotopy,Gambino2008identity,Hofmann1998groupoid} which
identifies in intensional type theory certain higher-dimensional features. The
idea is that, in such a semantics, we should be able to characterise dependent
sums and products universally in terms of weak, higher-dimensional adjoints to
substitution.

Eventually we expect to be able to construct a sound and complete semantics for
intensional type theory valued in weak $\omega$-categories. At the moment the
theory of weak $\omega$-categories is insufficiently well-developed for us to
describe this semantics; yet we can at least take steps towards it, by
describing semantics valued in simpler kinds of higher-dimensional category. In
this paper, we describe such a semantics valued in $2$-categories; which, as
well as objects representing types, and morphisms $f \colon X \to Y$
representing terms, has $2$-cells $\alpha \colon f \Rightarrow g$ representing
witnesses for the propositional equality of terms $f$ and $g$. Intuitively, the
$2$-categorical models we consider provide a notion of two-dimensional locally
cartesian closed category; though bearing in mind the above concerns regarding
the functoriality of substitution, it is in fact a ``split'' notion of
two-dimensional model which we will describe here. Relating this to a notion of
two-dimensional local cartesian closed category will require a $2$-categorical
coherence result along the lines of~\cite{Hofmann1995interpretation}, and we
defer this to a subsequent paper.

Our $2$-categorical semantics is sound and complete neither for intensional nor
extensional type theory, but rather for an intermediate theory which we call
\emph{two-dimensional type theory}. Recall that extensional type theory
distinguishes itself from intensional type theory by its admission of an
\emph{equality reflection rule}, which states that any two terms of type $A$
which are propositionally equal, are also definitionally equal. The
two-dimensional type theory that we will consider admits instances of the
equality reflection rule at just those types which are themselves identity
types. The effect this has is to collapse the higher-dimensional aspects of the
intensional theory, but only above dimension two: and it is this which allows a
complete semantics in \mbox{$2$-categories}. The leading example of model for
our semantics is the \emph{groupoid model} of~\cite{Hofmann1998groupoid};
indeed, it plays the same fundamental role for two-dimensional type theory as
the $\cat{Set}$-based model does for extensional type theory. However, we
expect there to be many more examples: on the categorical side, \emph{prestack}
and \emph{stack} models, which will provide two-dimensional analogues of the
\emph{presheaf} and \emph{sheaf} models of extensional type theory; and on the
type-theoretic side, an \emph{$E$-groupoid} model~\cite{Aczel1994Notes}, which
extends to two dimensions the \emph{setoid} model of extensional type
theory~\cite{Hofmann1994Elimination}. Once again, the task of describing these
models will be deferred to a subsequent paper.

Our hope is that the semantics we describe in this paper will provide a useful
guide in setting up more elaborate semantics for intensional type theory: both
of the $\omega$-categorical kind outlined above, and of the homotopy-theoretic
kind espoused in~\cite{Awodey2008Homotopy}. Indeed, most of the problematic
features of these higher-dimensional semantics are fully alive in the
two-dimensional case---in particular, the rather subtle issues regarding
stability of structure under substitution---and the analysis we give of them
here should prove useful in understanding these more general situations.

The paper is set out as follows. In Section~\ref{sec1}, we review the syntax of
intensional and extensional Martin-L\"of type theory and describe our
intermediate two-dimensional theory, $\mlt$. In Sections~\ref{sec2}
and~\ref{sec3} we describe a \mbox{$2$-categorical} structure built from the
syntax of \mlt. Section~\ref{sec2} makes use of the non-logical rules of \mlt\
together with the rules for identity types in order to construct a $2$-category
$\Con$ of contexts; a two-dimensional fibration $\Typ \to \Con$ of types over
contexts; and a comprehension $2$-functor $E \colon \Typ \to \Con^\mathbf 2$,
sending each type-in-context $\Gamma \vdash A \ \ty$ to the corresponding
\emph{dependent projection} map $(\Gamma \c x : A) \to \Gamma$. So far we have
given nothing more than a simple-minded extension of the one-dimensional
semantics; the twist is that each dependent projection in our $2$-categorical
model carries the structure of a \emph{normal isofibration}. This can be seen
as the semantic correlate of the \emph{Leibniz rule} in dependent type theory.
Section~\ref{sec3} considers the extra structure imposed on this basic
framework by the logical rules of \mlt. The identity types are characterised as
\emph{arrow objects} in the slices of the $2$-category of contexts; whilst the
unit type, dependent sums and dependent products admit description in terms of
a notion of weak $2$-categorical adjointness which we call \emph{retract
biadjunction}. Where a plain adjunction concerns itself with isomorphisms of
hom-sets $\mathbb C(FX, Y) \cong \mathbb D(X, GY)$, a retract biadjunction
instead requires retract equivalences of hom-categories $\tcat C(FX, Y) \simeq
\tcat D(X, GY)$. In particular, dependent sums and products are characterised
as left and right retract biadjoints to weakening. These syntactic
investigations lead us to define a notion of \emph{model} for two-dimensional
type theory, this being an arbitrary $2$-fibration $\Typ \to \Con$ equipped
with the structure outlined above; and the results of Sections~\ref{sec2}
and~\ref{sec3} can be summarised as saying that to each type theory $\Ss$
extending $\mlt$ we may assign a \emph{classifying} two-dimensional model
$\mathbb C(\Ss)$. In Section~\ref{sec4}, we provide a converse to this result
by showing that to each two-dimensional model $\mathbb C$ we can assign a
two-dimensional type theory $\Ss(\mathbb C)$ which represents the model
faithfully. We call this type theory the \emph{internal language} of $\mathbb
C$. Finally, we show that these two constructions---classifying model and
internal language---give rise to a \emph{functorial semantics} in the sense
of~\cite{Lawvere1968Some}: which is to say that they induce an equivalence
between suitably defined categories of two-dimensional type theories, and of
two-dimensional models.

\textbf{Acknowledgements} The author thanks the anonymous referees for their
helpful suggestions.

\newcommand{\rl}{\textsc}
\section{Intensional, extensional and two-dimensional type theory}\label{sec1}
\subsection{Intensional type theory} \label{intensional}
By \emph{intensional Martin-L\"of type theory}, we mean the logical calculus
set out in Part~II of~\cite{Nordstrom1990Programming}. In this paper, we
consider only the core calculus \mli, with type-formers for dependent sums,
dependent products, identity types and the unit type. We now summarise this
calculus, partly to fix notation and partly because there are few peculiarities
which are worth commenting on. The calculus has four basic forms of judgement:
$A \ \ty$ (``$A$ is a type''); $a : A$ (``$a$ is an element of the type $A$'');
$A = B \ \ty$ (``$A$ and $B$ are definitionally equal types''); and $a = b : A$
(``$a$ and $b$ are definitionally equal elements of the type $A$''). These
judgements may be made either absolutely, or relative to a context $\Gamma$ of
assumptions, in which case we write them as
\begin{equation*}
  \Gamma \vdash A\ \ty\text, \qquad
  \Gamma \vdash a : A\text, \qquad
  \Gamma \vdash A = B\ \ty \qquad \text{and} \qquad
  \Gamma \vdash a = b: A
\end{equation*}
respectively. Here, a \emph{context} is a list $\Gamma = (x_1 : A_1\c x_2 :
A_2\c \dots\c x_n : A_n)$, wherein each $A_i$ is a type relative to the context
$(x_1 : A_1\c \dots\c x_{i-1} : A_{i-1})$. There are now some rather natural
requirements for well-formed judgements: in order to assert that $a : A$ we
must first know that $A \ \ty$; to assert that $A = B \ \ty$ we must first know
that $A \ \ty$ and $B \ \ty$; and so on. We specify intensional Martin-L\"of
type theory as a collection of inference rules over these forms of judgement.
Firstly we have the \emph{equality rules}, which assert that the two judgement
forms $A = B \ \ty$ and $a = b : A$ are congruences with respect to all the
other operations of the theory; then we have the \emph{structural rules}, which
deal with weakening, contraction, exchange and substitution\footnote{Note in
particular that we take substitution to be a \emph{primitive}, rather than a
\emph{derived} operation: as done in~\cite{Jacobs1999Categorical}, for
instance.}; and finally, the \emph{logical rules}, which we list in Table
\ref{fig1}. Note that we commit the usual abuse of notation in leaving implicit
an ambient context $\Gamma$ common to the premisses and conclusions of each
rule. We also omit the rules expressing stability under substitution in this
ambient context.
\begin{table}
\emph{Dependent sum types}
\begin{equation*}
  \inferrule*[right=$\Sigma$-form;]{A\ \ty \\ x : A \t B(x) \ \ty}{\Sigma x:A .\, B(x) \ \ty}\qquad
  \inferrule*[right=$\Sigma$-intro;]{a : A\\ b : B(a)}{\spn{a, b} : \Sigma x:A .\, B(x)}
\end{equation*}
\
\begin{equation*}
    \inferrule*[right=$\Sigma$-elim;]{z : \Sigma x:A .\, B(x) \t C(z) \ \ty \\ x : A\c y : B(x) \t d(x, y) : C(\spn{x, y})}
    {z : \Sigma x:A .\, B(x) \t \Ell_d(z) : C(z)}
\end{equation*}\
\begin{equation*}
    \inferrule*[right=$\Sigma$-comp.]{z : \Sigma x:A .\, B(x) \t C(z) \ \ty \\ x : A\c y : B(x) \t d(x, y) : C(\spn{x, y})}
    {x : A\c y : B(x) \t \Ell_d(\spn{x, y}) = d(x, y) : C(\spn{x, y})}
\end{equation*}

\vskip\baselineskip\emph{Unit type}
\begin{equation*}
    \inferrule*[right=$\mathbf 1$-form;]{\ }{\mathbf 1 \ \ty} \qquad
    \inferrule*[right=$\mathbf 1$-intro;]{\ }{\mathord \star : \mathbf 1}
\end{equation*}
\begin{equation*}
    \inferrule*[right=$\mathbf 1$-elim;]{z : \mathbf 1 \t C(z) \ \ty \\ d : C(\mathord \star)}
    {z : \mathbf 1 \t \mathrm U_d(z) : C(z)}
\end{equation*}
\
\begin{equation*}
    \inferrule*[right=$\mathbf 1$-comp.]{z : \mathbf 1 \t C(z) \ \ty \\ d : C(\mathord \star)}
    {\mathrm U_d(\mathord \star) = d : C(\mathord \star)}
\end{equation*}

\vskip\baselineskip \emph{Identity types}
\begin{equation*}
    \inferrule*[right=$\Id$-form;]{A\ \ty \\ a, b : A}{\Id_A(a, b) \ \ty} \qquad
    \inferrule*[right=$\Id$-intro;]{A\ \ty \\ a : A}{\r(a) : \Id_A(a, a)}
\end{equation*}
\ \begin{equation*}
\inferrule*[right=$\Id$-elim;]{x, y : A \c p : \Id_A(x, y) \t C(x, y, p) \ \ty\\
  x : A \t d(x) : C(x, x, \r(x))}
  {x, y : A \c p : \Id_A(x, y) \t \J_{d}(x, y, p) : C(x, y, p)}
\end{equation*}
\
\begin{equation*}
\inferrule*[right=$\Id$-comp.]{x, y : A \c p : \Id_A(x, y) \t C(x, y, p) \ \ty\\
  x : A \t d(x) : C(x, x, \r(x))}
  {x : A \t \J_{d}(x, x, \r(x)) = d(x) : C(x, x, \r(x))}
\end{equation*}

\vskip\baselineskip\emph{Dependent product types}
\begin{equation*}
    \inferrule*[right=$\Pi$-form;]{A\ \ty \\ x : A \t B(x) \ \ty}{\Pi x:A .\, B(x) \ \ty} \qquad
    \inferrule*[right=$\Pi$-abs;]{x : A \t f(x) : B(x)}{\lambda x.\, f(x) : \Pi x:A .\, B(x)}
\end{equation*}
\
\begin{equation*}
    \inferrule*[right=$\Pi$-app;]{m : \Pi x:A .\, B(x)}
    {y : A \t m \cdot y : B(y)}\qquad
    \inferrule*[right=$\Pi$-$\beta$.]{x : A \t f(x) : B(x)}
    {y : A \t \big(\lambda x.\, f(x)\big) \cdot y = f(y) : B(y)}
\end{equation*}
\

\caption{Logical rules of intensional Martin-L\"of type theory
(\mli)}\label{fig1}
\end{table}
We will find it convenient to use the following extended forms of the identity
elimination and computation rules:
\begin{equation}\label{extended}
\begin{gathered}
\inferrule{x, y : A \c p : \Id_A(x, y) \c \Delta \t C(x, y, p) \ \ty\\
  x : A \c \Delta[x, x, \r(x) / x, y, p] \t d(x) : C(x, x, \r(x))}
  {x, y : A \c p : \Id_A(x, y) \c \Delta \t \J_{d}(x, y, p) : C(x, y, p)}
\\
\inferrule{x, y : A \c p : \Id_A(x, y) \c \Delta \t C(x, y, p) \ \ty\\
  x : A \c \Delta[x, x, \r(x) / x, y, p] \t d(x) : C(x, x, \r(x))}
  {x : A \c \Delta[x, x, \r(x) / x, y, p] \t \J_{d}(x, x, \r(x)) = d(x) : C(x, x, \r(x))}
\end{gathered}
\end{equation}
These rules may be derived from the elimination and computation rules in
Table~\ref{fig1} by using the $\Pi$-types to shift the additional contextual
parameter $\Delta$ onto the right-hand side of the turnstile.

\begin{Not}\label{notation}
We may omit from the premisses of a rule or deduction any hypothesis which may
be inferred from later hypotheses of that rule. Where it improves clarity we
may omit brackets in function applications, writing $hgfx$ in place of
$h(g(f(x)))$, for example. We may drop the subscript $A$ in an identity type
$\Id_A(a, b)$ where no confusion seems likely to occur. We may write a sum type
$\Sigma x : A.\, B(x)$ as $\Sigma(A, B)$, a product type $\Pi x : A.\, B(x)$ as
$\Pi(A, B)$, and a $\lambda$-abstraction $\lambda x.\, f(x)$ as $\lambda(f)$
(or using our applicative convention, simply $\lambda f$). It will occasionally
be useful to perform lambda-abstraction at the meta-theoretic level, for
instance writing $[x]\, f(x)$ to denote a term $f$ of the form $x : A \t f(x) :
B(x)$. We may write $\Gamma \t a \approx b : A$ to indicate that the type
$\Gamma \t \Id_A(a, b)$ is inhabited, and say that $a$ and $b$ are
\emph{propositionally equal}. We will also make use of \emph{vector notation}
in the style of~\cite{DeBruijn1991Telescopic}. Given a context $\Gamma = (x_1
: A_1, \dots, x_n : A_n)$, we may abbreviate a series of judgements:
\begin{equation*}
\t a_1 : A_1\text, \qquad
\t a_2 : A_2(a_1)\text, \qquad \dots \qquad
\t a_n : A_n(a_1, \dots, a_{n-1})\text,
\end{equation*}
as $\vdash \vec a : \Gamma$, where $\vec a
\defeq (a_1, \dots, a_n)$, and say that $\vec a$ is a \emph{global element} of $\Gamma$. We may also use this notation to abbreviate sequences of
hypothetical elements on the left-hand side of the turnstile; so, for example,
we may specify a dependent type in context $\Gamma$ as $\vec x : \Gamma \t
A(\vec x) \ \ty$. We will also make use of~\cite{DeBruijn1991Telescopic}'s
notion of \emph{telescope}. Given $\Gamma$ a context as before, this allows us
to abbreviate the series of judgements
\begin{equation*}
\begin{aligned}
\vec x : \Gamma & \t B_1(\vec x) \ \ty\text,\\
\vec x : \Gamma \c y_1 : B_1 & \t B_2(\vec x, y_1) \ \ty\text,\\
& \dots \\
\vec x : \Gamma \c y_1 : B_1 \c \dots \c y_{m-1} : B_{m-1} & \t B_m(\vec x, y_1, \dots y_{m-1}) \
\ty\text.
\end{aligned}
\end{equation*}
as $
  \vec x : \Gamma \t \Delta(\vec x) \ \ctxt
$, where $\Delta(\vec x) \defeq (y_1 : B_1(\vec x)\c y_2 : B_2(\vec x, y_1)\c
\dots )$. We say that $\Delta$ is a \emph{context in context $\Gamma$}, or a
\emph{context dependent upon $\Gamma$}, and refer to contexts like $\Delta$ as
\emph{dependent contexts}, and to those like $\Gamma$ as \emph{closed
contexts}. Given a dependent context $\vec x : \Gamma \t \Delta(\vec x) \
\ctxt$, we may abbreviate the series of judgements
\begin{equation*}
\begin{gathered}
\vec x : \Gamma \t f_1(\vec x) : B_1(\vec x)\\
\dots \\
\vec x : \Gamma \t f_m(\vec x) : B_m(\vec x, f_1(\vec x), \dots, f_{m-1}(\vec x))\text,
\end{gathered}
\end{equation*}
as $    \vec x : \Gamma \t \vec f(\vec x) : \Delta(\vec x)$, and say that $\vec
f$ is a \defn{dependent element} of $\Delta$. We can similarly assign a meaning
to the judgements
  $\vec x : \Gamma \t \Delta(\vec x) = \Theta(\vec x) \ \ctxt$ and $\vec x : \Gamma \t \vec f(\vec x) = \vec g(\vec x) : \Delta(\vec
  x)$,
expressing the definitional equality of two dependent contexts, and the
definitional equality of two dependent elements of a dependent context.
\end{Not}

\subsection{Extensional type theory} \label{extensional} We obtain extensional Martin-L\"of type
theory \mle\ by augmenting the intensional theory with the two \emph{equality
reflection rules}:
\begin{mathpar}
\inferrule{
  a, b : A \\ \alpha : \Id(a, b)}
  {a = b : A}

\inferrule{
  a, b : A \\ \alpha : \Id(a, b)}
  {\alpha = \r(a) : \Id(a, b)}
\end{mathpar}
together with the rule of \emph{function extensionality}:
\begin{equation*}
\inferrule{
  m, n : \Pi(A, B) \\ x : A \t m \cdot x = n \cdot x}
  {m = n : \Pi(A, B)}\ \text.
\end{equation*}
The addition of these three rules yields a type theory which is intuitively
simpler, and more natural from the perspective of categorical models, but
proof-theoretically unpleasant: we lose the decidability of definitional
equality and the decidability of type-checking. Note that if one develops
Martin-L\"of type theory in a framework admitting higher-order inference rules
(such as the Logical Framework of~\cite{Nordstrom1990Programming}) then the
above three rules are equipotent with the definitional $\eta$-rule.

\subsection{Two-dimensional type theory}\label{funextsec1} The type theory we investigate in this paper lies between
the intensional theory of~\S \ref{intensional} and the extensional theory of~\S
\ref{extensional}. We denote it by \mlt, and call it \emph{two-dimensional type
theory}, because as we will see, it has a natural semantics in two-dimensional
categories. It is obtained by augmenting intensional type theory with the rules
of Tables~\ref{fig2}~and~\ref{fig3}. These provide restricted versions of the
equality reflection rules (Table~\ref{fig2}) and the function extensionality
rules (Table~\ref{fig3}). To motivate the rules in Table~\ref{fig2}, we
introduce the notion of a discrete type. We say that $\Gamma \t A\ \ty$ is
\emph{discrete} if the judgements
\begin{mathpar}
  \inferrule*[right=$\Id$-refl$_1$-$A$;]{
    \Gamma \t a, b : A \\
    \Gamma \t p:\Id(a, b)
  }{
    \Gamma \t a = b : A
  }

  \inferrule*[right=$\Id$-refl$_2$-$A$]{
    \Gamma \t a, b : A \\
    \Gamma \t p : \Id(a, b)
  }{
    \Gamma \t p = \r(a) : \Id(a, b)
  }
\end{mathpar}
are derivable. Thus the intensional theory says that no types need be discrete;
the extensional theory says that all types are discrete; and the
two-dimensional theory says that all identity types are discrete. Note that
although two-dimensional type theory suffers from the same proof-theoretic
deficiencies of the extensional theory, it does so in a less severe manner:
indeed, only those types of $\mlt$ in whose construction the identity types
have been used will have undecidable definitional equality. As we ascend to
higher-dimensional variants of type theory, this undecidability will be pushed
further and further up the hierarchy of iterated identity types; but it is only
in the limit---which is intensional type theory---that we regain complete
decidability.

The necessity of the rules in Table~\ref{fig3} will become clear when we
reach~\S\ref{depprodssec}. We require them in order to obtain a satisfactory
notion of two-dimensional categorical model, in which dependent product
formation is right adjoint to substitution (in a suitably weak $2$-categorical
sense). The first of the rules in Table~\ref{fig3} is a propositional version
of the function extensionality principle of \S\ref{extensional}; whilst the
second and the third express coherence properties of the first. To understand
the third rule we must first explain the symbol $\ast$ appearing in it. It is a
definable constant which expresses that two propositionally equal elements of a
$\Pi$-type are pointwise propositionally equal. Explicitly, it satisfies the
following introduction and computation rules:
\begin{mathpar}
\inferrule*[right=$\ast$-intro;]{
  m, n : \Pi(A, B) \\ p : \Id(m, n) \\ a : A}
  {p \ast a : \Id(m \cdot a, n \cdot a)}

\inferrule*[right=$\ast$-comp;]{
  m : \Pi(A, B) \\ a : A}
  {\r(m) \ast a = \r(m \cdot a) : \Id(m \cdot a, m \cdot a)}
\end{mathpar}
and we may define it by $\Id$-elimination, taking $p \ast a \defeq \J_{[x]\r(x
\cdot a)}(m, n, p)$.
\begin{table}
\begin{mathpar}\inferrule*[right=$\Id$-disc$_1$;]{
  a, b : A \\ p, q : \Id(a, b) \\ \alpha : \Id(p, q)}
  {p = q : \Id(a, b)}\\
\inferrule*[right=$\Id$-disc$_2$.]{
  a, b : A \\ p, q : \Id(a, b) \\ \alpha : \Id(p, q)}
  {\alpha = \r(p) : \Id(p, q)}
\end{mathpar}
\caption{Rules for discrete identity types}\label{fig2}
 \end{table}
\begin{table}
\begin{mathpar}
\inferrule*[right=$\Pi$-ext{\text;}]{
  m, n : \Pi(A, B) \\ x : A \t p(x) : \Id(m \cdot x\c n \cdot x)}
  {\ext(m, n, p) : \Id(m, n)}

\inferrule*[right=$\Pi$-ext-comp\text;]{
  m : \Pi(A,B)}
  {\ext(m, m, [x]\,\r(m \cdot x)) = \r(m): \Id(m, m)}

\inferrule*[right=$\Pi$-ext-app{\text.}]{
  m, n : \Pi(A, B) \\ x : A \t p(x) : \Id(m \cdot x\c n \cdot x)}
  {x : A \t \ext(m, n, p) \ast x = p(x) : \Id(m \cdot x\c n \cdot x)}
\end{mathpar}
\caption{Rules for function extensionality}\label{fig3}
 \end{table}

\newcommand{\lift}[1]{{#1}^{\text{\tiny\textbullet}}}
\section{Categorical models for \mlt: structural aspects}\label{sec2}
\looseness=-1 The remainder of this paper will describe a notion of categorical
semantics for \mlt. In this section and the following one, we define a
syntactic category and enumerate its structure; whilst in Section~\ref{sec4},
we consider an arbitrary category endowed with this same structure, and derive
from it a type theory incorporating the rules of \mlt. This yields a semantics
which is both complete and sound. In this section, we define the basic
syntactic category and look at the structure induced on it by the non-logical
rules of \mlt. In the next section, we consider the logical rules. As mentioned
in the Introduction, the syntactic category we define will in fact be a
\emph{$2$-category}, whose objects will be (vectors of) types; whose morphisms
will be (vectors) of terms between those types; and whose $2$-cells will be
(vectors of) identity proofs between these terms. The various forms of $2$-cell
composition will be obtained using the identity elimination rules; whilst the
rules for discrete identity types given in Table~\ref{fig2} ensure that these
compositions satisfy the $2$-category axioms. For basic terminology and
notation relating to $2$-categories we refer to~\cite{Kelly1974Review}.

\subsection{One-dimensional semantics of type dependency}\label{onedimsemantics}
We begin by recalling the construction of a one-dimensional categorical
structure from the syntax of a dependent type theory. The presentation we have
chosen follows~\cite{Jacobs1993Comprehension} in its use of (full)
\emph{comprehension categories}. There are various other, essentially
equivalent, presentations that we could have used:
see~\cite{Cartmell1986Generalised,Dybjer1996Internal,Ehrhard1988Une,Hyland1989theory,Taylor1999Practical}
for example. We use comprehension categories because they afford a
straightforward passage to a two-dimensional structure.

So suppose given an arbitrary dependently-typed calculus $\Ss$ admitting the
same four basic judgement types and the same structural rules as the calculus
$\mli$. We define its \emph{category of contexts} $\con_\Ss$ to have as
objects, contexts $\Gamma$,~$\Delta$,~\dots, in $\Ss$, considered modulo
$\alpha$-conversion and definitional equality (so we identify $\Gamma$ and
$\Delta$ whenever $\t \Gamma = \Delta \ \ctxt$ is derivable); and as morphisms
$\Gamma \to \Delta$, judgements \mbox{$\vec x : \Gamma \t \vec f (\vec x) :
\Delta$}, considered modulo $\alpha$-conversion and definitional equality (so
we identity $f, g \colon \Gamma \to \Delta$ whenever $\vec x : \Gamma \t \vec f
(\vec x) = \vec g (\vec x)$ is derivable). To avoid further repetition, we
introduce the convention that any further categorical structures we define
should also be interpreted modulo $\alpha$-equivalence and definitional
equality. The identity map on $\Gamma$ is given by $\vec x : \Gamma \t \vec x :
\Gamma$; whilst composition is given by substitution of terms. Note that
$\con_\Ss$ has a terminal object, given by the empty context $(\ )$.

For each context $\Gamma$ we now define the category $\typ_\Ss(\Gamma)$ of
\emph{types-in-context-$\Gamma$}, whose objects $A$ are judgements $\vec x :
\Gamma \t A(\vec x)\ \ty$ and whose morphisms $A \to B$ are judgements $\vec x
: \Gamma\c y : A(\vec x) \t f(\vec x, y) : B(\vec x)$. Each morphism $f \colon
\Gamma \to \Delta$ of $\con_\Ss$ induces a functor $\typ_\Ss(f) \colon
\typ_\Ss(\Delta) \to \typ_\Ss(\Gamma)$ which sends a type $A$ in context
$\Delta$ to the type $\sb f A$ in context $\Gamma$ given by $\vec x : \Gamma \t
A(\vec f(\vec x)) \ \ty$. The assignation $f \mapsto \typ_\Ss(f)$ is itself
functorial in $f$, and so we obtain an indexed category $\typ_\Ss(\thg) \colon
\con_\Ss^\op \to \cat{Cat}$; which via the Grothendieck construction, we may
equally well view as a split fibration $p \colon \typ_\Ss \to \con_\Ss$. We
refer to this as the \emph{fibration of types over contexts}.

Explicitly, objects of $\typ_\Ss$ are pairs $(\Gamma, A)$ of a context and a
type in that context; whilst morphisms $(\Gamma, A) \to (\Delta, B)$ are pairs
$(f, g)$ of a context morphism $f \colon \Gamma \to \Delta$ together with a
judgement \mbox{$\vec x : \Gamma\c y : A(\vec x) \t g(\vec x, y) : B(\vec
f(\vec x))$}. The chosen cartesian lifting of a morphism $f \colon \Gamma \to
\Delta$ at an object $(\Delta, B)$ is given by $(f, \ro f) \colon (\Gamma, \sb
f B) \to (\Delta, B)$, where $\ro f$ denotes the judgement $\vec x : \Gamma\c y
: B (\vec f \vec x )\t y : B(\vec f \vec x)$. Now, for each object $(\Gamma,
A)$ of $\typ_\Ss$ we have the extended context \mbox{$\big(\vec x : \Gamma \c y
: A(\vec x)\big)$}, which we denote by $\Gamma . A$; and we also have the
judgement $\vec x : \Gamma\c y : A(\vec x) \t \vec x : \Gamma$, corresponding
to a context morphism $\pi_A \colon \Gamma.A \to \Gamma$ which we call the
\emph{dependent projection} associated to $A$. In fact, the assignation
$(\Gamma, A) \mapsto \pi_A$ provides the action on objects of a functor $E
\colon \typ_\Ss \to \con_\Ss^\mathbf 2$ (where $\mathbf 2$ denotes the arrow
category $0 \to 1$), whose action on maps sends the morphism $(f, g) \colon
(\Gamma, A) \to (\Delta, B)$ of $\typ_\Ss$ to the morphism
\begin{equation}\label{notpbsquare}
  \cd[@C+1em]{
    \Gamma.A \ar[d]_{\pi_A} \ar[r]^{f.g} &
    \Delta.B \ar[d]^{\pi_B} \\
    \Gamma \ar[r]_f &
    \Delta
  }
\end{equation}
of $\con_\Ss^\mathbf 2$, where $f.g$ denotes the judgement $\vec x : \Gamma \c
y : A \t (\vec f(\vec x), g(\vec x, y)) : \Delta.B$.

We can make two observations about this functor $E$. Firstly, it is fully
faithful, which says that every morphism $h \colon \Gamma.A \to \Delta.B$
fitting into a square like \eqref{notpbsquare} is of the form $f.g$ for a
unique $(f, g) \colon (\Gamma, A) \to (\Delta, B)$. Secondly, for a cartesian
morphism $(f, \ro f) \colon (\Gamma, \sb f B) \to (\Delta, B)$, the
corresponding square~\eqref{notpbsquare} is a pullback square. Indeed, given an
arbitrary commutative square
\begin{equation*}
  \cd[@C+1em]{
    \Lambda \ar[d]_{h} \ar[r]^{k} &
    \Delta.B \ar[d]^{\pi_\Delta} \\
    \Gamma \ar[r]_{f} &
    \Delta\text,
  }
\end{equation*}
commutativity forces $k$ to be of the form $\vec z : \Lambda \t (\vec f \vec h
\vec z, k' \vec z) : \Delta.B$ for some $\vec z : \Lambda \t k'(\vec z) :
B(\vec f \vec h \vec z)$; and so the required factorisation $\Lambda \to
\Gamma.\sb f B$ is given by the judgement $\vec z : \Lambda \t (\vec h \vec z,
k' \vec z) : \Gamma.\sb f B$. We may abstract away from the above situation as
follows. We define a \emph{full split comprehension category}
(cf.~\cite{Jacobs1993Comprehension}) to be  given by a category $\C$ with a
specified terminal object, together with a split fibration $p \colon \T \to \C$
 and a full and faithful functor $E
\colon \T \to \C^\mathbf 2$ rendering commutative the triangle
\begin{equation*}
  \cd[@!C]{
    \T \ar[dr]_p \ar[rr]^{E} & &
    \C^\mathbf 2 \ar[dl]^{\cod} \\ &
    \C\text,
  }
\end{equation*}
and sending cartesian morphisms in $\T$ to pullback squares in $\C^\mathbf 2$.
The preceding discussion shows that to any suitable dependent type theory $\Ss$
we may associate a full split comprehension category $\mathbb C(\Ss)$, which we
will refer to as the \emph{classifying comprehension category} of $\Ss$.

\begin{Not}\label{notconv}
We will make use of the notation developed above in arbitrary comprehension
categories $(p \colon \T \to \C, E \colon \T \to \C^\mathbf 2)$. Thus we write
chosen cartesian liftings as $(f, \ro f) \colon (\Gamma, \sb f B) \to (\Delta,
B)$, and write the image of $(\Gamma, A) \in \typ$ under $E$ as $\pi_A \colon
\Gamma.A \to \Gamma$. We will find it convenient to develop a little more
notation. Given $\Gamma \in \C$ and $A \in \T(\Gamma)$, we call a map $a \colon
\Gamma \to \Gamma.A$ satisfying $\pi_A a = \id_\Gamma$ a \emph{global section}
of $A$, and denote it by $a \in_\Gamma A$. Given further a morphism $f \colon
\Delta \to \Gamma$ of $\con$, we write $\sb f a \in_\Delta \sb f A$ for the
section of $\pi_{\sb f A}$ induced by the universal property of pullback in the
following diagram:
\begin{equation}\label{univproppullback}
    \cd{
      \Delta \ar@/^1em/[drr]^{a f} \ar@/_1em/[ddr]_{\id} \ar@{.>}[dr] \\ &
      \Delta.\sb f A \ar[r]^{f.\ro f} \ar[d]_[0.4]{\pi_{\sb f A}} &
      \Gamma.A \ar[d]^{\pi_A} \\ &
      \Delta \ar[r]_f &
      \Gamma\text.
    }
\end{equation}
\end{Not}

\subsection{A $2$-category of types}\label{sec2cattypes}
We will now extend the classifying comprehension category $\mathbb C(\Ss)$
defined above to a classifying comprehension $2$-category. We will not need the
full strength of two-dimensional type theory, \mlt, for this. Rather, for the
rest of this section we fix an arbitrary dependently typed theory $\Ss$ which
admits the structural rules required in the previous subsection together with
the identity type rules from Table~\ref{fig1} and the discrete identity rules
of Table~\ref{fig2}. Our first task will be to construct a $2$-category of
closed types in $\Ss$. We will do this by enriching the category $\typ_\Ss(\ )$
of closed types with $2$-cells derived from the $2$-category of \emph{strict
internal groupoids} in $\Ss$. A strict internal groupoid in $\Ss$ is given by a
closed type $A_0$; a family $A_1(x, y)$ of types over $x, y : A_0$; and
operations of unit, composition and inverse:
\begin{align*}
  x : A_0 & \t \id_x : A_1(x, x)\\
  x, y, z : A_0 \c p : A_1(x, y) \c q : A_1(y, z) & \t q \circ p : A_1(x, z)\text,\\
  x, y : A_0 \c p : A_1(x, y) & \t p^{-1} : A_1(y, x)\text,
\end{align*}
which obey the usual five groupoid axioms up to definitional equality. For
instance, the left unit axiom requires that
\begin{equation*}
  x, y : A_0 \c p : A_1(x, y) \t \id_y \circ p = p : A_1(x, y)
\end{equation*}
should hold. We will generally write that $(A_0, A_1)$ is an internal groupoid
in $\Ss$, leaving the remaining structure understood. Now an \emph{internal
functor} $F \colon (A_0, A_1) \to (B_0, B_1)$ between internal groupoids is
given by judgements
\begin{align*}
  x : A_0 & \t F_0(x) : B_0\\
  x, y : A_0 \c p : A_1(x, y) & \t F_1(p) : B_1(F_0 x, F_0 y)\text,
\end{align*}
subject to the usual two functoriality axioms (up to definitional equality
again); whilst an \emph{internal natural transformation} $\alpha \colon F
\Rightarrow G$ is given by a family of components $x : A_0 \t \alpha(x) :
B_1(F_0 x, G_0 x)$ subject to the (definitional) naturality axiom.

\begin{Prop}\label{strict2cat}
The strict groupoids, functors and natural transformations internal to $\Ss$
form a $2$-category $\tcat{Gpd}(\Ss)$ which is \emph{locally groupoidal}, in
the sense that every $2$-cell is invertible.
\end{Prop}

\begin{proof}
Recall that for any category $\E$, we can define a $2$-category
$\tcat{Gpd}(\E)$ of groupoids internal to that category\footnote{One commonly
requires the category $\E$ to have pullbacks, but this is inessential.}. In
particular, we have the $2$-category $\tcat{Gpd}(\con_\Ss)$ of groupoids
internal to the category of contexts of $\Ss$. Now, each strict internal
groupoid $\A$ in $\Ss$ gives rise to such an internal groupoid $\A'$ in
$\con_\Ss$ whose object of objects is the context $(x : A_0)$ and whose object
of morphisms is the context $(x : A_0\c y : A_0\c p : A_1(x, y))$. We can check
that internal functors $\A \to \B$ in $\Ss$ correspond bijectively with
internal functors $\A' \to \B'$ in $\con_\Ss$; and that this correspondence
extends to the natural transformations between them. Thus we may take
$\tcat{Gpd}(\Ss)$ to be the $2$-category whose objects are strict internal
groupoids in $\Ss$, whose hom-categories are given by $\tcat{Gpd}(\Ss)(\A, \B)
:= \tcat{Gpd}(\con_\Ss)(\A', \B')$, and whose remaining structure is inherited
from $\tcat{Gpd}(\con_\Ss)$. Note that every $2$-cell of $\tcat{Gpd}(\con_\Ss)$
is invertible, so that the same obtains for $\tcat{Gpd}(\Ss)$
\end{proof}

Our method for obtaining the $2$-category of closed types will be to construct
a functor $\typ_\Ss(\ ) \to \tcat{Gpd}(\Ss)$, and to lift the $2$-cell
structure of $\tcat{Gpd}(\Ss)$ along it.

\begin{Prop} \label{functor}
To each closed type $A$ in $\Ss$ we may assign a strict internal groupoid $(A,
\Id_A)$; and the assignation $A \mapsto (A, \Id_A)$ underlies a functor
$\typ_\Ss(\ ) \to \tcat{Gpd}(\Ss)$.
\end{Prop}

\begin{proof}
The proof of this result is essentially due to~\cite{Hofmann1998groupoid}. We
repeat it because we will need the details. We first show that $(A, \Id_A)$ has
the structure of a strict internal groupoid. For identities, we take $x : A \t
\id_x := \r(x) : \Id(x, x)$. For composition, we require a judgement
\begin{equation*}
  x, y, z : A \c p : \Id(x, y) \c q : \Id(y, z) \t q \circ p : \Id(x, z)\text;
\end{equation*}
and by $\Id$-elimination on $p$---in the extended form given in
equation~\eqref{extended}---it suffices to define this when $y = z$ and $q =
r(y)$, for which we take $r(y) \circ p := p$. Similarly, to give the judgement
\begin{equation*}
  x, y : A \c p : \Id(x, y) \t p^{-1} : \Id(y, x)
\end{equation*}
providing inverses, it suffices to consider the case $x = y$ and $p = \r(x)$;
for which we take $\r(x)^{-1} = \r(x)$. We must now check the five groupoid
axioms. The first unitality axiom $\id_y \circ p = p$ follows from the
$\Id$-computation rule. For the other unitality axiom, it suffices, by discrete
identity types, to show that
\begin{equation*}
  x, y : A \c p : \Id(x, y) \t p \circ \r(x) \approx p : \Id(x, y)
\end{equation*}
holds; and by $\Id$-elimination, it suffices to do this in the case $x = y$ and
$p = \r(x)$, for which we have that $\r(x) \circ \r(x) = \r(x)$ as required.
Likewise, for the associativity axiom, it suffices to show that
\begin{multline*}
  w, x, y, z : A \c p : \Id(w, x) \c q : \Id(x, y) \c s : \Id(y, z) \\
  \t s \circ (q \circ p) \approx (s \circ q) \circ p : \Id(w, z)\text;
\end{multline*}
and again by $\Id$-elimination, it suffices to do this when $y = z$ and $s =
\r(y)$, when have that $\r(y) \circ (q \circ p) = q \circ p = (\r(y) \circ q)
\circ p$ as required. Note that again, we require the extended form of
$\Id$-elimination of equation~\eqref{extended}, and in future we will use this
without further comment. The invertibility axioms are similar. Suppose now that
in addition to $A$ we are given another type $B$ together with a judgement $x :
A \t f(x) : B$ between them. We will extend this to an internal functor $(f,
\lift f) \colon (A, \Id_A) \to (B, \Id_B)$. We define the action on hom-types
\begin{equation*}
  x, y : A \c p : \Id(x, y) \t \lift f(p) : \Id(fx, fy)
\end{equation*}
by $\Id$-elimination on $p$: for when $x = y$ and $p = \r(x)$, we may take
$\lift f(\r(x)) := \r(f(x))$. We must now check the functoriality axioms. That
$(f, \lift f)$ preserves identities follows from the $\Id$-computation rule;
whilst to to show that it preserves binary composition, it suffices by discrete
identity types to show that
\begin{equation*}
  x, y, z : A \c p : \Id(x, y) \c q : \Id(y, z)
  \t \lift f(q \circ p) \approx \lift f(q) \circ \lift f(p) : \Id(fx, fz)
\end{equation*}
holds; and this follows by $\Id$-elimination on $q$, since when $y = z$ and $q
= r(y)$, we have that $\lift f(\r(y) \circ p) = \lift f(p) = \r(f(y)) \circ
\lift f(p) = \lift f(\r(y)) \circ \lift f(p)$ as required. We must now check
that the assignation $f \mapsto (f, \lift f)$ is itself functorial. To show
that it preserves identities, we must show that for any closed type $A$,
\begin{equation*}
  x, y : A \c p : \Id(x, y) \t \lift{(\id_A)}(p) = p : \Id(x, y)
\end{equation*}
holds. By discrete identity types, it suffices to show this up to mere
propositional equality; and by $\Id$-elimination, we need only do so in the
case when $x = y$ and $p = \r(x)$, when we have that $\lift{(\id_A)}(\r(x)) =
\r(\id_A(x)) = \r(x)$ as required. To show that $f \mapsto (f, \lift f)$
respects composition, we must show that for maps of closed types $f \colon A
\to B$ and $g \colon B \to C$, the judgement
\begin{equation*}
  x, y : A \c p : \Id(x, y) \t \lift{(gf)}(p) = \lift g (\lift f (p)) : \Id(g f x, g f y)
\end{equation*}
holds. Again, it suffices to do this only up to propositional equality, and
this only in the case where $x = y$ and $p = \r(x)$; whereupon we have that
$\lift{(gf)}(\r(x)) = \r(g(f(x))) = \lift g (\r(f(x)) = \lift g(\lift f(\r(x))$
as required.
\end{proof}

\begin{Cor}\label{2cattypes}
The category $\typ_\Ss(\ )$ of closed types in $\Ss$ may be extended to a
locally groupoidal $2$-category $\Typ_\Ss(\ )$ whose $2$-cells $\alpha \colon f
\Rightarrow g \colon A \to B$ are judgements $x : A \t \alpha(x) : \Id_B(fx,
gx)$\text.
\end{Cor}

\begin{proof} If we view $\typ_\Ss(\ )$ as a $2$-category with only
identity $2$-cells, then the functor of the previous proposition may be seen as
a $2$-functor $\typ_\Ss(\ ) \to \tcat{Gpd}(\Ss)$. We can factorise this
2-functor as a composite
\begin{equation*}
  \typ_\Ss(\ ) \to \Typ_\Ss(\ ) \to \tcat{Gpd}(\Ss)\text,
\end{equation*}
whose first part is bijective on objects and $1$-cells and whose second part is
fully faithful on $2$-cells; and we now define $\Typ_\Ss(\ )$ to be the
intermediate \mbox{$2$-category} in this factorisation. We must check that this
definition agrees with the description of $\Typ_\Ss(\ )$ given above. Clearly
this is so for the objects and morphisms; whilst for the $2$-cells, we must
show that for any $f, g \colon A \to B$, each judgement
  $x : A \t \alpha(x) : \Id_B(fx, gx)$
satisfies the axiom for an internal natural transformation $\alpha \colon (f,
\lift f) \Rightarrow (g, \lift g)$. By discrete identity types, this amounts to
validating the judgement
\begin{equation*}
  x, y : A \c p : \Id_A(x, y) \t \lift g(p) \circ \alpha(x) \approx \alpha(y) \circ \lift f(p) : \Id_B(f x, g y)\text;
\end{equation*}
and by $\Id$-elimination on $p$, it suffices to do this in the case where $x =
y$ and $p = \r(x)$: for which we have that $\lift g(\r(x)) \circ \alpha(x) =
\r(g(x)) \circ \alpha(x) = \alpha(x) = \alpha(x) \circ \r(f(x)) = \alpha(x)
\circ \lift f(\r(x))$, as required.
\end{proof}

\begin{Cor}\label{2cattypesdep}
For any context $\Gamma$ in $\Ss$, the category $\typ_\Ss(\Gamma)$ of
types-in-context-$\Gamma$ may be extended to a locally groupoidal $2$-category
$\Typ_\Ss(\Gamma)$ wherein $2$-cells $\alpha \colon f \Rightarrow g$ are
judgements $ \vec x : \Gamma \c y : A \t \alpha(\vec x, y) : \Id_B(f(\vec x,
y)\c g(\vec x, y))$.
\end{Cor}

\begin{proof}
We consider the \emph{slice theory} $\Ss / \Gamma$, whose closed types are the
types of $\Ss$ in context $\Gamma$. It is easy to see that $\Ss / \Gamma$
admits the same inference rules as $\Ss$---and in particular has discrete
identity types---so that the result follows upon identifying $\Typ_\Ss(\Gamma)$
with $\Typ_{\Ss/\Gamma}(\ )$.
\end{proof}

\subsection{A $2$-category of contexts}\label{sec2catctxts}
In this section, we generalise the construction of the $2$-category of closed
types in order to construct a $2$-category of contexts. The method will be a
direct transcription of the one used in the previous section, but in order for
it to make sense, we need to extend the identity type constructor to a
``meta-constructor'' which operates on entire contexts rather than single
types.

\begin{Prop}\label{idcontexts}
The following inference rules are definable in $\Ss$.
\begin{equation*}
  \inferrule*[right=$\Id$-form';]{
    \Phi\ \ctxt \\
    \vec a, \vec b : \Phi
  }{
    \Id_\Phi(\vec a, \vec b) \ \ctxt
  } \qquad
  \inferrule*[right=$\Id$-intro';]{
    \Phi\ \ctxt \\
    \vec a : \Phi
  }{
    \vec\r(\vec a) : \Id_\Phi(\vec a, \vec a)
  }
\end{equation*}\
\begin{equation*}
  \inferrule*[right=$\Id$-elim';]{
    \vec x, \vec y : \Phi \c \vec p : \Id_\Phi(\vec x, \vec y) \c \Delta
      \t \Theta(\vec x, \vec y, \vec p) \ \ctxt\\
    \vec x : \Phi \c \Delta[\vec x, \vec x, \r(\vec x) / \vec x, \vec y, \vec p]
      \t \vec d(\vec x) : \Theta(\vec x, \vec x, \r(\vec x))
  }{
    \vec x, \vec y : \Phi \c \vec p : \Id_\Phi(\vec x, \vec y) \c \Delta
      \t \J_{\vec d}(\vec x, \vec y, \vec p) : \Theta(\vec x, \vec y, \vec p)
  }
\end{equation*}\
\begin{equation*}
  \inferrule*[right=$\Id$-comp'.]{
    \vec x, \vec y : \Phi \c \vec p : \Id_\Phi(\vec x, \vec y) \c \Delta
      \t \Theta(\vec x, \vec y, \vec p) \ \ctxt\\
    \vec x : \Phi \c \Delta[\vec x, \vec x, \r(\vec x) / \vec x, \vec y, \vec p]
      \t \vec d(\vec x) : \Theta(\vec x, \vec x, \r(\vec x))
  }{
    \vec x : \Phi \c \Delta[\vec x, \vec x, \r(\vec x) / \vec x, \vec y, \vec p]
      \t \J_{\vec d}(\vec x, \vec x, \r(\vec x)) = \vec d(\vec x) : \Theta(\vec x, \vec x, \r(\vec x))
  }
\end{equation*}
\end{Prop}
In order to prove this result, we will make use of the following well-known
consequence of the identity type rules:
\begin{Prop}[The Leibniz rule]\label{leibnitz}
Given $A\ \ty$ and $x : A \t B(x) \ \ty$ in $\Ss$, the following rules are
derivable:
\begin{mathpar}
  \inferrule*[right=$\Id$-subst;]{
    a_1, a_2 : A \\
    p : \Id(a_1, a_2) \\
    b_2 : B(a_2)
  }{
    p^\ast( b_2) : B(a_1)
  }

  \inferrule*[right=$\Id$-subst-comp.]{
    a : A \\
    b : B(a)
  }{
    \r(a)^\ast (b) = b : B(a)
  }
\end{mathpar}
\end{Prop}

\begin{proof}
By $\Id$-elimination on $p$, it suffices to derive the first rule in the case
where $a_1 = a_2$ and $p = \r(a_1)$: in which case we can take $\r(a_1)^\ast
(b) := b$. The second rule now follows from the $\Id$-computation rule.
\end{proof}

The key idea behind the proof of Proposition~\ref{idcontexts} can be
illustrated by considering a context $\Phi = (x : A\c y : B(x))$ of length~$2$.
The corresponding identity context $\Id_\Phi$ will be given by
\begin{equation*}
\Id_\Phi\big((x, y), (x', y')\big) := \big(p : \Id_A(x, x')\c q : \Id_{B(x)}(y,
p^\ast y')\big)\text.
\end{equation*}
We use substitution along the first component $p$ to make the second component
$q$ type-check. This can be seen as a type-theoretic analogue of the
Grothendieck construction for fibrations. Indeed, it is possible to show that
there is a propositional isomorphism between this identity context $\Id_\Phi$
and the identity type $\Id_{\Sigma(A, B)}$. Thus in principle it is unnecessary
to introduce identity contexts; however, we prefer to do so in order to obtain
a cleaner separation between the identity rules and the $\Sigma$-type rules.

\begin{proof}[Proof of Proposition~\ref{idcontexts}]
The proof has two stages. First, we define the generalised $\Id$-inference
rules in the special case where the context $\Phi$ has length~$1$; and then we
use these to define them in the general case. We will reduce syntactic clutter
by proving our results only in the case where the postcontext $\Delta$ is
empty: the reader may readily supply the annotations for the general case. For
the first part of the proof, we suppose ourselves given a context $\Phi = (x :
A)$ of length 1. The inference rules \textsc{$\Id$-form'} and
\textsc{$\Id$-intro'} for $\Phi$ are just the usual $\Id$-formation and
$\Id$-introduction rules for $A$. However, \textsc{$\Id$-elim'} corresponds to
the following generalised elimination rule:
\begin{equation}\label{genelim}
\inferrule{x, y : A \c p : \Id(x, y) \t \Theta(x, y, p) \ \ctxt\\
  x : A \t \vec d(x) : \Theta(x, x, \r(x))}
  {x, y : A \c p : \Id(x, y) \t \J_{\vec d}(x, y, p) : \Theta(x, y, p)}
\end{equation}
with \textsc{$\Id$-comp'} stating that $\J_{\vec d}(a, a, \r(a)) = \vec d(a)$.
We will define the elimination rule by induction on the length $n$ of the
context $\Theta$. When $n = 0$, this is trivial, and when $n = 1$, we use the
usual identity elimination rule. So suppose now that we have defined the rule
for all contexts $\Theta$ of length $n$, and consider a context $
  x, y : A \c p : \Id(x, y)
    \t \Theta(x, y, p) \ \ctxt
$ of length $n+1$. Thus $\Theta$ is of the form
\begin{equation*}
  \Theta(x, y, p) = (\vec u : \Lambda(x, y, p) \c v : D(x, y, p, \vec u))
\end{equation*}
for some context $\Lambda$ of length $n$ and type $D$. It follows that to make
a judgement \mbox{$x : A \t \vec d(x) : \Theta(x, x, \r(x))$} is equally well
to make a pair of judgements
\begin{equation}\label{2.14}
    \begin{aligned}
        x : A & \t {\vec d}_1(x) : \Lambda(x, x, \r(x)) \\
        x : A & \t d_2(x) : D\big(x, x, \r(x), {\vec d}_1(x)\big)\text.
    \end{aligned}
\end{equation}
By the inductive hypothesis, we may apply the elimination rule~\eqref{genelim}
for the context $\Lambda$ with eliminating family ${\vec d}_1$ to deduce the
existence of a term
\begin{equation}\label{2.15}
  x, y : A \c p : \Id(x, y)
    \t \J_{{\vec d}_1}(x, y, p) : \Lambda(x, y, p)\text,
\end{equation}
satisfying $\J_{{\vec d}_1}(x, x, \r(x)) = \vec{d}_1(x)$. Now we consider the
dependent type
\begin{equation}\label{2.16}
  x, y : A \c p : \Id(x, y) \t C(x, y, p) := D(x, y, p, \J_{{\vec d}_1}(x, y, p)\big) \ \ty\text.
\end{equation}
We have that $C(x, x, \r(x)) = D\big(x, x, \r(x), \J_{{\vec d}_1}(x, x,
\r(x))\big) = D\big(x, x, \r(x), {\vec d}_1(x)\big)$ and so from \eqref{2.14}
we can derive the judgement
\begin{equation}\label{2.17}
   x : A \t d_2(x) : C\big(x, x, \r(x)\big)\text.
\end{equation}
Now applying the standard $\Id$-elimination rule to \eqref{2.16} and
\eqref{2.17} yields a judgement
\begin{equation}\label{2.18}
   x, y : A \c p : \Id(x, y) \t \J_{d_2}(x, y, p) : D\big(x, y, p, \J_{{\vec d}_1}(x, y, p)\big)
\end{equation}
satisfying $\J_{d_2}(x, x, \r(x)) = d_2(x)$. But to give \eqref{2.15} and
\eqref{2.18} is equally well to give a dependent element $x, y : A \c p :
\Id(x, y) \t \J_{\vec d}(x, y, p) : \Theta(x, y, p)$; and the respective
computation rules for $\J_{{\vec d}_1}$ and $\J_{d_2}$ now imply the
computation rule for $\J_{\vec d}$. This completes the first part of the proof.

We now construct the generalised inference rules for an arbitrary
context~$\Phi$. Once again the proof will be by induction, this time on the
length of~$\Phi$. For the base case, the only context of length $0$ is $(\ )$,
the empty context. For this, we take the identity context $\Id_{(\ )}$ also to
be the empty context. The introduction rule is vacuous, whilst the elimination
rule requires us to provide, for each closed context $\Theta$ and global
element $\vec d : \Theta$, a global element $\J_{\vec d} : \Theta$, satisfying
the computation rule $\J_{\vec d} = \vec d : \Theta$. Thus we simply take
$\J_{\vec d} := \vec d$ and are done. Suppose now that we have defined identity
contexts for all contexts of length $n$, and consider a context $\Phi =
\big(\vec x_1 : \Lambda, x_2 : D(\vec x_1)\big)$ of length $n+1$. In order to
define $\Id_\Phi$, we first apply the inductive hypothesis to $\Lambda$ in
order to define its Leibniz rule. Thus given $\vec x : \Lambda \t \Upsilon(\vec
x) \ \ctxt$, we may define a judgement
\begin{equation*}
  \vec x, \vec y : \Lambda \c \vec p : \Id_\Lambda(\vec x, \vec y) \c \vec z : \Upsilon(\vec y) \t \vec p^\ast(\vec z)  : \Upsilon(\vec
  x)\text,
\end{equation*}
satisfying $\r(\vec x)^\ast(\vec z) = \vec z : \Upsilon(\vec x)$. The proof is
as in Proposition~\ref{leibnitz}. Now, to give the formation rule for
$\Id_\Phi$ is equally well to give a judgement
\begin{equation*}
  \vec x_1 : \Lambda \c y_1 : D(\vec x_1) \c \vec x_2 : \Lambda \c y_2 : D(\vec x_2) \t
  \Id_\Phi(\vec x_1, y_1, \vec x_2, y_2)\ \ctxt\text,
\end{equation*}
which we do by setting
\begin{equation*}
   \Id_\Phi(\vec x_1, y_1, \vec x_2, y_2) := \big({\vec p} : \Id_\Lambda({\vec x}_1, {\vec x}_2) \c q : \Id_{D({\vec x}_1)}(y_1, \vec p^\ast y_2)\big)\text.
\end{equation*}
Next, to define the introduction rule for $\Id_\Phi$ is equally well to give
judgements
\begin{equation*}
    \begin{aligned}
        \vec x : \Lambda \c y : D(\vec x) & \t r_1(\vec x, y) : \Id_\Lambda({\vec x}, {\vec x})\\
        \vec x : \Lambda \c y : D(\vec x) & \t r_2(\vec x, y) : \Id_{D({\vec x})}\big(y\c \vec r_1(\vec x, y)^\ast (y)\big)
    \end{aligned}
\end{equation*}
which we do by setting $
    \vec r_1(\vec x, y) := \r(\vec x)$ and $r_2(\vec x, y) :=
    \r(y)$, where for the second of these we make
use of the fact that $\Id_{D({\vec x})}\big(y, \vec \r(\vec x)^\ast (y)\big) =
\Id_{D({\vec x})}(y, y)$. In order to define the elimination rule for
$\Id_\Phi$, we first define a context
%$\Delta$ by
%\begin{equation*}
%   \b I_1 := \big(\ \vec a_1 : \Lambda \c \vec b_1 : \Lambda \c \vec p_1 : \Id_\Lambda(\vec a_1, \vec b_1)\ \big)\text,
%\end{equation*}
dependent on $
   \vec x_1, \vec x_2 : \Lambda$ and $\vec p : \Id_\Lambda(\vec x_1, \vec x_2)$
%     \t \b \Delta(\vec a_1, \vec a_2, \vec p_1) \ \ctxt
by
\begin{equation*}
    \Delta(\vec x_1, \vec x_2, \vec p) := \big(y_1 : D(\vec x_1) \c y_2 : D(\vec x_2) \c q : \Id_{D({\vec x}_1)}(y_1, \vec p^\ast y_2)\big)\text.
\end{equation*}
We may then write the premisses of the elimination rule for $\Id_\Phi$ as:
\begin{equation}\label{2.24}
  \vec x_1, \vec x_2 : \Lambda \c \vec p : \Id_\Lambda(\vec x_1, \vec x_2) \c \vec z : \Delta(\vec x_1, \vec x_2, \vec p)
    \t \Theta(\vec x_1, \vec x_2, \vec p, \vec z)\ \ctxt
\end{equation}
and
\begin{equation}\label{2.25}
  \vec x : \Lambda \c y : D(\vec x)
      \t \vec d(\vec x, y) : \Theta(\vec x, \vec x, \r(\vec x), y, y, \r(y))\text.
\end{equation}
We would like to apply the elimination rule for $\Id_\Lambda$ (with postcontext
$\Delta$) to equation \eqref{2.24}. In order to do so, we need to exhibit a
generating family
\begin{equation}\label{2.26}
  \vec x : \Lambda \c \vec z : \Delta(\vec x, \vec x, \r(\vec x))
    \t \vec d'(\vec x, \vec z) : \Theta(\vec x, \vec x, \r(\vec x), \vec z)\text;
\end{equation}
which is equivalently a family
\begin{equation*}
  \vec x : \Lambda \c y_1, y_2 : D(\vec x) \c q : \Id_{D(\vec x)}(y_1,
  y_2) \t \vec d'(\vec x, y_1, y_2, q) : \Theta(\vec x, \vec x, \r(\vec x), y_1, y_2, q)
\end{equation*}
since we have that $\r(\vec x)^\ast (y_2) = y_2$. But we may obtain such a
family by applying the generalised elimination rule \eqref{genelim} for
$\Id_{D(\vec x)}$ to the dependent context
\begin{equation*}
  \vec x : \Lambda \c y_1, y_2 : D(\vec x) \c q : \Id_{D(\vec x)}(y_1, y_2)
    \t \Theta(\vec x, \vec x, \r(\vec x), y_1, y_2, q) \ \ctxt
\end{equation*}
 with eliminating family \eqref{2.25}. This yields a judgement \eqref{2.26} as
required, whilst the computation rule says that $\vec d'(\vec x, y, y, \r(y)) =
\vec d(\vec x, y)$. Now applying the elimination rule for $\Id_\Lambda$ to
\eqref{2.24} and \eqref{2.26} yields a judgement
\begin{multline*}
  \vec x_1, \vec x_2 : \Lambda \c \vec p : \Id_\Lambda(\vec x_1, \vec x_2) \c \vec z :
  \Delta(\vec x_1, \vec x_2, \vec p)
    \t \J_{\vec d'}(\vec x_1, \vec x_2, \vec p, \vec z) : \Theta(\vec x_1, \vec x_2, \vec p, \vec
    z)\text,
\end{multline*}
of the correct form to provide the conclusion of the elimination rule for
$\Id_\Phi$. From the computation rule for $\Id_\Lambda$, this will satisfy
$\J_{\vec d'}(\vec x, \vec x, \r(\vec x), \vec z) = \vec d'(\vec x, \vec z)$,
and so in particular, we obtain that
\begin{equation*}
  \J_{\vec d'}(\vec x, \vec x, \r(\vec x), y, y, \r(y)) = \vec d'(\vec x, y, y,
  \r(y)) = \vec d(\vec x, y)
\end{equation*}
which gives us the computation rule for $\Id_\Phi$.
\end{proof}

Using Proposition~\ref{idcontexts} we can now construct the $2$-category of
contexts in $\Ss$ by mimicking the developments of~\S\ref{sec2cattypes}. We
first define a
\defn{strict groupoid context} in $\Ss$ to be given by a context
$\Gamma_0$ together with a dependent family $
    \vec x, \vec y : \Gamma_0 \t \Gamma_1(\vec x, \vec y) \ \ctxt
$ of hom-contexts, and operations of unit, composition and inverse satisfying
the groupoid axioms as before. It is still the case that any groupoid context
induces an internal groupoid object in the category of contexts $\con_\Ss$; and
so with the obvious definition of functor and natural transformation, we obtain
a $2$-category $\tcat{GpdCtxt}(\Ss)$ of groupoid contexts in $\Ss$. Following
Proposition~\ref{functor}, we now define a functor $\con_\Ss \to
\cat{GpdCtxt}(\Ss)$ sending $\Gamma$ to $(\Gamma, \Id_\Gamma)$. A small
subtlety we must check in order for this to go through is that $\Ss$ has not
only discrete identity types, but also discrete identity \emph{contexts}; and
this follows by a straightforward induction on the length of a context.
Thereafter, the argument of Proposition~\ref{2cattypes} carries over to give:

\begin{Cor}\label{2catctxt}
The category $\con_\Ss$ of contexts in $\Ss$ may be extended to a locally
groupoidal $2$-category $\Con_\Ss$ whose $2$-cells $\alpha \colon f \Rightarrow
g \colon \Gamma \to \Delta$ are judgements $\vec x : \Gamma \t \alpha(\vec x) :
\Id_\Delta(\vec f \vec x, \vec g \vec x)$.
\end{Cor}

We end this section with a simple observation:

\begin{Prop}
The $2$-category $\Con_\Ss$ has a $2$-terminal object given by the empty
context $(\ )$.
\end{Prop}
\begin{proof}
It is clear that every context $\Gamma$ admits a unique morphism $! \colon
\Gamma \to (\ )$, which makes $(\ )$ a terminal object. For it to be
$2$-terminal, we must also show that for any $2$-cell $\alpha \colon !
\Rightarrow ! \colon \Gamma \to (\ )$ we have $\alpha = \id_!$. But this
follows because we defined $\Id_{(\ )} := (\ )$ in the proof of
Proposition~\ref{idcontexts}.
\end{proof}

\subsection{A $2$-fibration of types over contexts}\label{2fibsec}
The next stage in our development will be to extend the fibration of types over
contexts to a $2$-fibration of types over contexts. In~\S\ref{onedimsemantics},
we built the one-dimensional fibration by first defining an indexed category of
types over contexts, and then applying the Grothendieck construction. In the
two-dimensional case it turns out that the indexed $2$-category of types over
contexts has a structure so elaborate (it is given by a trihomomorphism
$\Typ_\Ss(\thg) \colon \Con_\Ss^\coop \to \cat{Gray}$, where $\cat{Gray}$ is
the tricategory of $2$-categories, $2$-functors, pseudo-natural
transformations, and modifications; see~\cite{Gordon1995Coherence}) that it is
that it is significantly less work to construct the associated $2$-fibration
directly. We begin by recalling from~\cite{Hermida1999Some} the definition of
$2$-fibration. Of the several equivalent formulations given there, the most
convenient for our purposes is the following:
\begin{Defn}\label{2fibdef}(cf.~\cite[Theorem 2.8]{Hermida1999Some})
Let $\tcat E$ and $\tcat B$ be $2$-categories. We say that a $2$-functor $p
\colon \tcat{E} \to \tcat{B}$ is a \emph{cloven $2$-fibration} if the following
four conditions are satisfied:
\begin{enumerate}[(i)]
\item The underlying ordinary functor of $p$ is a cloven fibration of
    categories;
\item Each cartesian 1-cell $f \colon y \to z$ of $\tcat E$ has the
    following two-dimensional universal property: that whenever we are
    given a $2$-cell $\alpha \colon g \Rightarrow h \colon x \to z$ of
    $\tcat{E}$
together with a factorisation
\begin{equation*}
  p(\alpha) =
  \cd{
    p(x) \ar@/^1em/[r]^{k} \ar@/_1em/[r]_{l} \dtwocell{r}{\gamma} & p(y) \ar[r]^{p(f)} & p(z)\text,
  }
\end{equation*}
we may lift this to a unique factorisation
\begin{equation*}
  \alpha =
  \cd{
    x \ar@/^1em/[r]^{k'} \ar@/_1em/[r]_{l'} \dtwocell{r}{\gamma'} & y \ar[r]^{f} & z
  }
\end{equation*}
satisfying $p(\gamma') = \gamma$.
\item For each $x, y \in \tcat E$, the induced functor $p_{x, y} \colon
    \tcat E(x, y) \to \tcat B(px, py)$ is a cloven fibration of categories;
\item For each $x, y, z \in \tcat E$ and $f \colon x \to y$, the functor
    $(\thg) \cdot f \colon \tcat E(y, z) \to \tcat E(x, z)$ %\emph{strictly}
    preserves cartesian liftings of $2$-cells.
\end{enumerate}
We say further that a cloven $2$-fibration is \emph{globally split} if its
underlying fibration of categories in (i) is a split fibration.
\end{Defn}

We will now show that the split fibration $p \colon \typ_\Ss \to \con_\Ss$ of
types over contexts extends to a globally split $2$-fibration $p \colon
\Typ_\Ss \to \Con_\Ss$. The first step will be to construct the total
$2$-category $\Typ_\Ss$. Before doing this we prove a useful lemma.
\begin{Lem}\label{pilift}
For a dependent projection $\pi_A \colon \Gamma.A \to \Gamma$ of $\Con_\Ss$,
its lifting to an internal functor $(\pi_A, \lift {\pi_A})$, as defined in
Proposition~\ref{functor}, satisfies
\begin{equation*}
    (\vec x, y)\c (\vec x', y') : \Gamma.A \c (\vec p, q) : \Id_{\Gamma.A}\big((\vec x, y)\c (\vec x', y')\big) \t \lift {\pi_A} (\vec p, q) = \vec p : \Id_\Gamma(\vec x, \vec x')\text.
\end{equation*}
\end{Lem}
\begin{proof}
By discrete identity types, it suffices to show that $\lift {\pi_A}(\vec p, q)
\approx \vec p$; and by $\Id$-elimination on $\Gamma.A$, we need only consider
the case where \mbox{$\vec x = \vec x'$}, \mbox{$y = y'$}, \mbox{$\vec p =
\r(\vec x)$} and \mbox{$q = r(y)$}. But here, by definition of $\lift {\pi_A}$,
we have $\lift {\pi_A}(\r(\vec x), \r(y)) = \r(\pi_A(\vec x, y)) = \r(\vec x)$
as required.
\end{proof}

\begin{Prop}\label{total2cat}
The category $\typ_\Ss$ defined in~\S\ref{onedimsemantics} extends to a locally
groupoidal $2$-category $\Typ_\Ss$ whose $2$-cells $(\alpha, \beta) \colon (f,
g) \Rightarrow (f', g') \colon (\Gamma, A) \rightarrow (\Delta, B)$ are given
by pairs of judgements
\begin{equation}\label{2cellsoftyp}
\begin{aligned}
  \vec x : \Gamma & \t \vec \alpha(\vec x) : \Id_\Delta(\vec f \vec x, \vec f' \vec x) \\
  \vec x : \Gamma \c y : A(\vec x) & \t \beta(\vec x, y) : \Id_{B(\vec f \vec x)}\big(g(\vec x, y)\c \vec \alpha(\vec x)^\ast(g'(\vec x, y))\big)\text.
\end{aligned}
\end{equation}
\end{Prop}
\begin{proof}If we view $\typ_\Ss$ as a $2$-category with only identity
$2$-cells, then the functor $E \colon \typ_\Ss \to \con_\Ss^\mathbf 2$ defined
in~\S\ref{onedimsemantics} may be viewed as a $2$-functor $\typ_\Ss \to
\Con_\Ss^\mathbf 2$. We can factorise this 2-functor as a composite
\begin{equation}\label{compositefunctor}
  \typ_\Ss \to \Typ_\Ss \to \Con_\Ss^\mathbf 2\text,
\end{equation}
whose first part is bijective on objects and $1$-cells and whose second part is
bijective on $2$-cells. We claim that the intermediate $2$-category is the
$\Typ_\Ss$ of the Proposition. Clearly it has the correct objects and
$1$-cells, whilst for the $2$-cells, we must show that given maps $(f, g), (f',
g') \colon (\Lambda, A) \to (\Delta, B)$ of $\typ_\Ss$, pairs of judgements as
in~\eqref{2cellsoftyp} are in bijection with diagrams
\begin{equation}\label{2celldigram}
  \cd[@+1em]{
    \Gamma. A \ar[d]_{\pi_A} \ar@/^1em/[r]^{f.g} \ar@/_1em/[r]_{f'.g'} \dtwocell{r}{\gamma} &
    \Delta.B \ar[d]^{\pi_B} \\
    \Gamma \ar@/^1em/[r]^{f} \ar@/_1em/[r]_{f'} \dtwocell{r}{\alpha}  &
    \Delta
  }
\end{equation}
in $\Con_\Ss$ satisfying $\pi_B \gamma = \alpha \pi_A$. For a diagram
like~\eqref{2celldigram}, the $2$-cell $\gamma \colon f.g \Rightarrow f'.g'$
corresponds---by the definition of $\Id_{\Delta.B}$ given in
Proposition~\ref{idcontexts}---to a pair of judgements
\begin{equation}\label{compartison2}
\begin{aligned}
  \vec x : \Gamma \c y : A(\vec x) & \t \vec \gamma_1(\vec x, y) : \Id_\Delta(\vec f \vec x, \vec f' \vec x) \\
  \vec x : \Gamma \c y : A(\vec x) & \t \gamma_2(\vec x, y) : \Id_{B(\vec f \vec x)}\big(g(\vec x, y)\c \vec \gamma_1(\vec x, y)^\ast(g'(\vec x, y))\big)\text,
\end{aligned}
\end{equation}
whilst the equality $\pi_B \gamma = \alpha \pi_A$ corresponds to the validity
of the judgement
\begin{equation*}
  \vec x : \Gamma \c y : A(\vec x) \t \vec \alpha(\vec x) = \lift {\pi_B}(\gamma(\vec x, y)) : \Id_\Delta(\vec f \vec x, \vec f' \vec x)\text.
\end{equation*}
But by Lemma~\ref{pilift}, we have $\lift {\pi_B}(\gamma(\vec x, y)) =
\gamma_1(\vec x, y)$, so that $\alpha(\vec x) = \gamma_1(\vec x, y)$, and we
may identify~\eqref{compartison2} with~\eqref{2cellsoftyp} upon taking $\beta
:= \gamma_2$.
\end{proof}
\begin{Cor}
The fully faithful functor $E \colon \typ_\Ss \rightarrow \con_\Ss^\mathbf 2$
of~\S\ref{onedimsemantics} extends to a $2$-fully faithful (i.e., bijective on
$1$- and $2$-cells) $2$-functor $E \colon \Typ_\Ss \to \Con_\Ss^\mathbf 2$.
\end{Cor}
\begin{proof}
We take $E$ to be the second half of the factorisation
in~\eqref{compositefunctor}.
\end{proof}
We now define $p \colon \Typ_\Ss \to \Con_\Ss$ to be the composite of the
$2$-functor $E$ of the previous Proposition with the codomain $2$-functor
$\Con_\Ss^\mathbf 2 \to \Con_\Ss$; explicitly, $p$ is the $2$-functor sending
$(\Gamma, A)$ to $\Gamma$, $(f, g)$ to $f$ and $(\alpha, \beta)$ to $\alpha$.
We intend to show that $p$ is a (globally split) $2$-fibration; and will do so
by making using of two further properties of the $2$-functor $E \colon \Typ_\Ss
\to \Con_\Ss^\mathbf 2$. The first of these generalises directly the
one-dimensional situation described in~\S\ref{onedimsemantics}. Its proof is
less straightforward than one might think.

\begin{Prop}\label{2pullbacklemma}
For each $(\Delta, B) \in \Typ_\Ss$ and $f \colon \Gamma \to \Delta$ in
$\Con_\Ss$, the following pullback square in $\Con_\Ss$ is also a $2$-pullback:
\begin{equation}\label{2pullbacksquare}
\cd[@C+1em]{  \Gamma. \sb f B \ar[d]_{\pi_{\sb f B}} \ar[r]^{f.\ro f} & \Delta. B \ar[d]^{\pi_B} \\
  \Gamma \ar[r]_{f}  & \Delta\text.}
\end{equation}
\end{Prop}
\begin{proof}
We begin by introducing a piece of local notation: for the duration of this
proof, we will write applications of the Leibniz rule as
\begin{equation*}
  \inferrule*[right=$\Id$-subst.]{
    a_1, a_2 : A \\
    p : \Id(a_1, a_2) \\
    b_2 : B(a_2)
  }{
    \mathsf{subst}_B(p, b_2) : B(a_1)
  }
\end{equation*}
We do this in order to make explicit the family $B$ in which substitution is
occurring. Now, to say that \eqref{2pullbacksquare} is not just a pullback but
also a $2$-pullback is to say that, whenever we are given maps $h, k \colon
\Lambda \to \Gamma.\sb f B$ and $2$-cells
\begin{equation}\label{arb2}
\cd[@+3em]{
  \Lambda \ar@/_1em/[d]_{\pi_{\sb f B} h} \rtwocell{d}{\alpha} \ar@/^1em/[d]^{\pi_{\sb f B} k} \ar@/^1em/[r]^{(f.\ro f) h} \dtwocell{r}{\beta} \ar@/_1em/[r]_{(f.\ro f) k} &
  \Delta.B \ar[d]^{\pi_B} \\
  \Gamma \ar[r]_{f}  & \Delta}
\end{equation}
in $\Con_\Ss$ satisfying $f\alpha = \pi_B \beta$, we can find a unique $2$-cell
$\gamma \colon h \Rightarrow k \colon \Lambda \to \Gamma.\sb f B$ satisfying
$\pi_{\sb f B} \circ \gamma = \alpha$ and $f.\ro f \circ \gamma = \beta$. In
order to show this, we will first need to understand how $f.\ro f$ lifts to an
internal functor
\begin{equation*}
    (f.\ro f\c \lift{(f.\ro f)}) \colon
(\Gamma.\sb f B\c \Id_{\Gamma.\sb f B}) \to (\Delta.B\c
\Id_{\Delta.B})\text.
\end{equation*}
So suppose given elements $(\vec x_1, y_1)$ and $(\vec x_2, y_2) : \Gamma.\sb f
B$; now a typical element $(\vec p, q) : \Id_{\Gamma.\sb f B}\big((\vec x_1,
x_2)\c (\vec y_1, y_2)\big)$ is given by a pair of judgements
\begin{equation}\label{from}
  \vec p : \Id_\Gamma(\vec x_1, \vec y_1) \qquad \text{and} \qquad
  q : \Id_{B (\vec f\vec x_1)}\big(x_2\c \mathsf{subst}_{\sb f B}(\vec p, y_2)\big)\text.
\end{equation}
This is sent by $\lift{(f.\ro f)}$ to some element $(\vec u, v) :
\Id_{\Delta.B}\big((\vec f \vec x_1, y_1)\c (\vec f \vec y_1, y_2)\big)$, which
is equally well a pair of judgements
\begin{equation}\label{to}
  \vec u : \Id_\Gamma(\vec f \vec x_1, \vec f \vec y_1) \qquad \text{and} \qquad
  v : \Id_{B(\vec f \vec x_1)}\big(x_2\c \mathsf{subst}_{B}(\vec u, y_2)\big)\text.
\end{equation}
Since we have $\pi_B \circ f.\ro f = f \circ \pi_{\sb f B}$, we have by
Lemma~\ref{pilift} that
\begin{equation*}
    \vec u = \lift{\pi_B}(\vec u, v) = \lift{(\pi_\B \circ f.\ro f)}(\vec p, q) =
    \lift{(f \circ \pi_{\sb f B})}(\vec p, q) = \lift{\vec f}(\vec p)\text;
\end{equation*}
and so it remains only to describe $v$. We will do this by reduction to a
special case. Suppose that we have $x_2 = \mathsf{subst}_{\sb f B}(\vec p,
y_2)$ and $q = \r(\mathsf{subst}_{\sb f B}(\vec p, y_2))$. We denote the
corresponding $v$ by
\begin{equation}\label{thetas}
    \theta(\vec p, y_2) : \Id_{B(\vec f \vec x_1)}\big(\mathsf{subst}_{\sb f B}(\vec p, y_2)\c \mathsf{subst}_{B}(\lift {\vec f}(\vec p), y_2)\big)\text.
\end{equation}
Note that in the case where $x_1 = y_1$ and $p = \r(x_1)$, we have by
$\Id$-computation that $\theta(\r(x_1), y_2) = \r(y_2)$. We now
use~\eqref{thetas} to describe the general case. We claim that given $\vec p$
and $q$ as in~\eqref{from}, the corresponding $v$ as in~\eqref{to} satisfies
\begin{equation*}
    v = \theta(\vec p, y_2) \circ q : \Id_{B(\vec f \vec x_1)}\big(x_2\c \mathsf{subst}_{B}(\lift{\vec f}(\vec p), y_2)\big)\text.
\end{equation*}
Now, by discrete $\Id$-types, it suffices to show this up to propositional
equality; and by $\Id$-elimination on $\Gamma.\sb f B$, this only in the case
where $\vec x_1 = \vec y_1$, $\vec p = \r(\vec x_1)$, $x_2 = y_2$ and $q =
\r(x_2)$. Here, by definition of $\lift{(f.\ro f)}$ and $\Id$-computation, we
have on the one hand that $v = \r(x_2)$; but on the other that $\theta(\r(\vec
x_1), x_2) \circ \r(x_2) = \r(x_2) \circ \r(x_2) = \r(x_2)$ as required. This
completes the proof of the claim.

We are now ready to show that~\eqref{2pullbacksquare} is a $2$-pullback. So
suppose given maps $h, k \colon \Lambda \to \Gamma.\sb f B$ and $2$-cells
$\alpha, \beta$ as in~\eqref{arb2}. To give $h$ is to give judgements $x :
\Lambda \t \vec h_1(x) : \Gamma$ and $x : \Lambda \t h_2(x) : B(\vec f \vec h_1
x)$---and correspondingly for $k$---whilst to give $\alpha$ and $\beta$ as in
\eqref{arb2} satisfying $f \alpha = \pi_B \beta$ is to give judgements
\begin{equation*}
\begin{aligned}
  x : \Lambda & \t \vec \alpha (x) : \Id_\Gamma(\vec h_1 x , \vec k_1 x) \\
  x : \Lambda & \t \vec \beta_1 (x) : \Id_\Delta(\vec f \vec h_1 x, \vec f \vec k_1 x) \\
  x : \Lambda & \t \beta_2 (x) : \Id_{B(\vec f \vec h_1 x)}(h_2 x\c \mathsf{subst}_B(\vec \beta_1 x\c k_2 x))
\end{aligned}
\end{equation*}
satisfying
\begin{equation*}
    x : \Lambda \t \lift{\vec f}(\alpha x) = \lift{\pi_B} (\vec \beta_1 x, \beta_2 x) : \Id_\Delta(\vec f \vec h_1 x, \vec f \vec k_1 x)\text.
\end{equation*}
By Lemma~\ref{pilift}, we have that $\lift{\pi_B} (\vec \beta_1 x, \beta_2 x) =
\beta_1(x)$; and so to give~\eqref{arb2} satisfying $f \alpha = \pi_B \beta$ is
equally well to give a pair of judgements
\begin{equation*}
\begin{aligned}
  x : \Lambda & \t \vec \alpha (x): \Id_\Gamma(\vec h_1 x, \vec k_1 x) \\ \text{and} \qquad
  x : \Lambda & \t \beta_2 (x): \Id_{B(\vec f \vec h_1 x)}(h_2 x\c \mathsf{subst}_B(\lift{\vec f} \vec \alpha x, k_2 x)\big)\text.
\end{aligned}
\end{equation*}
From this we are required to find a unique $2$-cell $\gamma \colon h
\Rightarrow k \colon \Lambda \to \Gamma.\sb f B$ satisfying $\pi_{\sb f B}
\circ \gamma = \alpha$ and $(f.\ro f) \circ \gamma = \beta$; which is equally
well a pair of judgements
\begin{equation*}
\begin{aligned}
  x : \Lambda & \t \vec \gamma_1 (x) : \Id_\Gamma(\vec h_1 x, \vec k_1 x) \\ \text{and} \qquad
  x : \Lambda & \t \gamma_2(x) : \Id_{B(\vec f \vec h_1 x)}(h_2 x\c \mathsf{subst}_{\sb f B}(\vec \gamma_1 x, k_2 x)\big)
\end{aligned}
\end{equation*}
satisfying $\lift{(\pi_{\sb f B})}(\gamma_1 x, \gamma_2 x) = \alpha(x)$ and
$\lift{(f.\ro f)}(\gamma_1 x, \gamma_2 x) = (\lift{\vec f} \alpha x, \beta_2
x)$. Now by Lemma~\ref{pilift}, we have $\lift{(\pi_{\sb f B})}(\gamma_1 x,
\gamma_2 x) = \gamma_1(x)$, whence we must take $\gamma_1 \defeq \alpha$;
whilst from our investigations above, we have
\begin{equation*}
\lift{(f.\ro f)}(\gamma_1 x, \gamma_2 x) = \lift{(f.\ro f)}(\alpha x,
\gamma_2 x) = \big(\lift{\vec f}(\alpha x)\c \theta(\alpha x, k_2 x) \circ \gamma_2
(x)\big)
\end{equation*}
which tells us that we must have $\gamma_2(x) \defeq \theta(\alpha x, k_2
x)^{-1} \circ \beta_2 (x)$. Uniqueness of $\gamma$ follows easily.
\end{proof}

The second property of $E$ we consider has no one-dimensional analogue, as it
involves the inherently $2$-categorical notion of \emph{isofibration}:

\begin{Defn}\label{isofibdefn}
Let $\tcat K$ be a $2$-category. A morphism $p \colon X \to Y$ in $\tcat K$ is
said to be a \emph{cloven isofibration} if for every invertible $2$-cell
\begin{equation}\label{isofibdiag}
    \cd[@!C]{
      W \ar[rr]^{g} \ar[dr]_{f} & \rtwocell{d}{\alpha} &
      X \ar[dl]^p \\ &
      Y\text,
    }
\end{equation}
we are given a choice of $1$-cell $s_\alpha \colon W \to X$ and $2$-cell
$\sigma_\alpha \colon s_\alpha \Rightarrow g$ satisfying \mbox{$p \circ
s_\alpha = f$} and $p \circ \sigma_\alpha = \alpha$; and these choices are
natural in $W$, in the sense that given further $k \colon W' \to W$, we have
$s_{\alpha k} = s_\alpha \circ k$ and $\sigma_{\alpha k} = \sigma_\alpha \circ
k$. A cloven isofibration is said to be \emph{normal} if for any $g \colon W
\to X$, we have $s_{\id_{pg}} = g$ and $\sigma_{\id_{pg}} = \id_g$.
\end{Defn}

\begin{Prop}\label{normalisofibration}
Every dependent projection $\pi_B \colon \Delta.B \to \Delta$ in $\Con_\Ss$ may
be equipped with the structure of a normal isofibration.
\end{Prop}
\begin{proof}
Suppose given an invertible $2$-cell
\begin{equation}\label{invert2cell}
    \cd[@!C]{
      \Gamma \ar[rr]^{g} \ar[dr]_{f} & \rtwocell{d}{\alpha} &
      \Delta.B \ar[dl]^{\pi_B} \\ &
      \Delta
    }
\end{equation}
of $\Con_\Ss$. We must find a $1$-cell $s_\alpha \colon \Gamma \to \Delta.B$
and $2$-cell $\sigma_\alpha \colon s_\alpha \Rightarrow g$ satisfying
\mbox{$\pi_B \circ s_\alpha = f$} and $\pi_B \circ \sigma_\alpha = \alpha$.
Now, to give a $2$-cell as in \eqref{invert2cell} is equally well to give
judgements
\begin{equation*}
\begin{aligned}
  \quad \vec x : \Gamma & \t \vec f(\vec x) : \Delta\text, &
  \quad \vec x : \Gamma & \t \vec g_1(\vec x) : \Delta\text,\\
  \quad \vec x : \Gamma & \t g_2(\vec x) : B(\vec g_1 \vec x)\text, &
  \quad \vec x : \Gamma & \t \vec \alpha(\vec x) : \Id(\vec f \vec x, \vec g_1 \vec x)\text.
\end{aligned}
\end{equation*}
So we may take $s_\alpha \colon \Gamma \to \Delta.B$ to be given by the pair of
judgements
\begin{equation}\label{j1}
  \vec x : \Gamma  \t \vec f(\vec x) : \Delta \quad \text{and} \quad
  \vec x : \Gamma  \t \vec (\alpha \vec x)^\ast(g_2 \vec x) : B(\vec f \vec x)\text,
\end{equation}
and take $\sigma_\alpha \colon s_\alpha \Rightarrow g$ to be given by the pair
of judgements
\begin{equation}\label{j2}
\begin{aligned}
  \vec x : \Gamma &
    \t \vec \alpha(\vec x) : \Id(\vec f \vec x, \vec g_1 \vec x)\\
  \vec x : \Gamma\c y : A(\vec x) &
    \t \r\big((\alpha x)^\ast(g_2 \vec x)\big) :
    \Id\big(\vec (\alpha \vec x)^\ast(g_2 \vec x)\c (\alpha \vec x)^\ast(g_2 \vec x)\big)\text.
\end{aligned}
\end{equation}
Given further $k \colon \Lambda \to \Gamma$, the equalities $s_{\alpha k} =
s_\alpha \circ k$ and $\sigma_{\alpha k} = \sigma_\alpha \circ k$ correspond
precisely to the stability of \eqref{j1} and \eqref{j2} under substitution in
$\vec x$. Thus $\pi_B$ is a cloven isofibration; and it remains to check
normality. But when $\alpha$ is an identity $2$-cell we have $f(\vec x) = \vec
g_1(\vec x)$ and $\alpha(\vec x) = \r(\vec g_1(\vec x))$ and so by the Leibniz
computation rule, \eqref{j1} reduces to $g$ and \eqref{j2} to $\id_g$ as
required.
\end{proof}
We will refer to the isofibration structure described in
Proposition~\ref{normalisofibration} as the \emph{canonical isofibration
structure} on a dependent projection.

\begin{Rk}\label{linkage}
Proposition~\ref{normalisofibration} provides a link between the
$2$-categorical semantics of this paper and the homotopy-theoretic semantics
espoused by Awodey and Warren in~\cite{Awodey2008Homotopy}. The key idea of
that paper is that a suitable environment for modelling intensional type theory
should be a category equipped with a \emph{weak factorisation system} $(\L,
\R)$ (in the sense of~\cite{Bousfield1977Constructions}) whose right-hand
class of maps $\R$ is used to model dependent projections. Now, any finitely
complete $2$-category carries a weak factorisation system $(\L, \R)$ wherein
$\R$ is the class of normal isofibrations; it forms one half of
what~\cite[Section~4]{Gambino2008Homotopy} calls the ``dual of the natural
model structure on a $2$-category''. Thus our two-dimensional semantics fits
naturally into the framework outlined in~\cite{Awodey2008Homotopy}.

This result can also be seen as a special case of a result obtained
in~\cite{Gambino2008identity}. The main result of that paper is that the
classifying category of any intensional type theory may be equipped with a weak
factorisation system whose right class of maps is generated by the dependent
projections; and it is shown (Lemma 13) that the maps in this right class are
``type-theoretic normal isofibrations''. Our
Proposition~\ref{normalisofibration} can be seen as a two-dimensional collapse
of this result.
\end{Rk}

Using Propositions \ref{2pullbacklemma} and \ref{normalisofibration}, we may
now show that:
\begin{Prop}\label{basic2fib}
The $2$-functor $p \colon \Typ_\Ss \to \Con_\Ss$ is a globally split
$2$-fibration.
\end{Prop}
\begin{proof}
We check the four clauses in Definition~\ref{2fibdef}. Clause (i) is immediate,
since the underlying ordinary functor of $p \colon \Typ_\Ss \to \Con_\Ss$ is
the split fibration $p \colon \typ_\Ss \to \con_\Ss$. For clause (ii), it
suffices to consider a chosen cartesian lifting $(f, \ro f) \colon (\Gamma, \sb
f B) \rightarrow (\Delta, B)$ of $\Typ_\Ss$. Taking advantage of the $2$-fully
faithfulness of $E \colon \Typ_\Ss \to \Con_\Ss^\mathbf 2$, we may express the
property we are to verify as follows: that for each diagram
\begin{equation*}
  \cd[@+1em]{
    \Lambda.A \ar[d]_{\pi_A} \ar@/^1em/[rr]^{h_1} \ar@/_1em/[rr]_{h_2} \dtwocell{rr}{\beta} & &
    \Delta.B \ar[d]^{\pi_B} \\
    \Lambda \ar@/^1em/[r]^{g_1} \ar@/_1em/[r]_{g_2} \dtwocell{r}{\alpha} &
    \Gamma \ar[r]_f & \Delta
  }
\end{equation*}
in $\Con_\Ss$ with $\pi_B \beta = f \alpha \pi_A$, there is a unique
factorisation
\begin{equation*}
  \beta =
  \cd[@C+1em]{
    \Lambda.A \ar@/^1em/[r]^{h_1'} \ar@/_1em/[r]_{h_2'} \dtwocell{r}{\beta'} & \Gamma.\sb f B \ar[r]^{f.\ro f} & \Delta.B
  }
\end{equation*}
with $\pi_{\sb f B} \beta' = \beta \pi_A$. But this follows without difficulty
from the fact that diagram \eqref{2pullbacksquare} is a $2$-pullback. For
clause (iii) in the definition of $2$-fibration, we suppose given $(\Gamma, A)$
and $(\Delta, B)$ in $\Typ_\Ss$ and are required to show that the functor
$\Typ_\Ss\big((\Gamma, A), (\Delta, B)\big) \to \Con_\Ss(\Gamma, \Delta)$ is a
fibration. Using once more the $2$-fully faithfulness of $E$, it suffices to
show that for each commutative square
\begin{equation*}
  \cd{
    \Gamma.A \ar[d]_{\pi_A} \ar[r]^{g.h} &  %\ar@/^1em/[r]^{k} \dtwocell{r}{\gamma} &
    \Delta.B \ar[d]^{\pi_B} \\
    \Gamma \ar[r]_{g} &
    \Delta
  }
\end{equation*}
in $\Con_\Ss$ and $2$-cell $\alpha \colon f \Rightarrow g$, we can find a
$1$-cell $k \colon \Gamma.A \to \Delta.B$ and a \mbox{$2$-cell} $\beta \colon k
\Rightarrow g.h$ satisfying $\pi_B k = f \pi_A$ and $\pi_B \beta = \alpha
\pi_A$. This follows using the canonical isofibration structure of $\pi_B$.
Finally, for clause (iv), we must show that each \mbox{$(\thg) \cdot f \colon
\Typ_\Ss(y, z) \to \Typ_\Ss(x, z)$} preserves cartesian liftings of $2$-cells.
As every $2$-cell of $\Typ_\Ss$ is invertible, and hence cartesian, this is
automatic.
\end{proof}

We end this section by considering the pullback stability of the canonical
isofibration structures of Proposition~\ref{normalisofibration}. To this end,
consider a square like~\eqref{2pullbacksquare}. Both vertical arrows $\pi_B$
and $\pi_{\sb f B}$ have their canonical isofibration structures; but we also
have a second isofibration structure on $\pi_{\sb f B}$, obtained by pulling
back the canonical structure of $\pi_B$ along $f$. A careful examination of the
proof of Proposition~\ref{normalisofibration} reveals that \emph{these two
structures on $\pi_{\sb f B}$ need not coincide}. In other words, the canonical
isofibration structures of Proposition~\ref{normalisofibration} are not
necessarily stable by pullbacks. At first glance, this may appear surprising,
since stability by pullbacks tends to be an automatic consequence of stability
under substitution. However, a more careful analysis shows that in this case,
stability under substitution corresponds to a more restricted form of pullback
stability, which we now describe.

Suppose we are given $\Delta \in \Con_\Ss$, $A \in \Typ_\Ss(\Delta)$ and $B \in
\Typ_\Ss(\Delta.A)$. We can view the dependent projection $\pi_B \colon
\Delta.A.B \to \Delta.A$ not only as a map of $\Con_\Ss$, but also as a map
\begin{equation}\label{diagtopullback}
  \cd[@!C@C-2em]{
    \Delta.A.B \ar[dr]_{\pi_A \pi_B} \ar[rr]^{\pi_B} & &
    \Delta.A \ar[dl]^{\pi_A} &
    \\ &
    \Delta
  }
\end{equation}
of $\Con_\Ss / \Delta$. It is easy to see that the forgetful $2$-functor
$\Con_\Ss / \Delta \to \Con_\Ss$ creates normal isofibrations, so
that~\eqref{diagtopullback} is canonically a normal isofibration in $\Con_\Ss /
\Delta$. Suppose we are now given a morphism $f \colon \Gamma \to \Delta$ of
$\Con_\Ss$. By pulling back~\eqref{diagtopullback} along $f$, we obtain the map
\begin{equation}\label{diagtopullback2}
  \cd[@!C@C-2.5em]{
    \Gamma.\sb f A.\sb f B \ar[dr]_{\pi_{\sb f A} \pi_{\sb f B}} \ar[rr]^{\pi_{\sb f B}} & &
    \Gamma.\sb f A \ar[dl]^{\pi_{\sb f A}} &
    \\ &
    \Gamma
  }
\end{equation}
of $\Con_\Ss / \Gamma$ (note that we are abusing notation slightly: we should
really write the left-hand vertex as $\Gamma.\sb f A.\sb {(f.\ro f)} B$), and
this now has two isofibration structures on it: the one induced by the
canonical isofibration structure on $\pi_{\sb f B}$, and the one obtained by
pulling back the isofibration structure of~\eqref{diagtopullback}. The
following Proposition now tells us that these two isofibration structures
on~\eqref{diagtopullback2} \emph{do} coincide.

\begin{Prop}\label{stabilityiso}
Suppose given $\Delta \in \Con_\Ss$, $A \in \Typ_\Ss(\Delta)$ and $B \in
\Typ_\Ss(\Delta.A)$ and $f \colon \Gamma \to \Delta$ as above. With reference
to the $2$-pullback square
\begin{equation}\label{special2pullback}
\cd[@C+1em]{  \Gamma. \sb f A. \sb f B \ar[d]_{\pi_{\sb f B}} \ar[r]^{f.\ro f.\ro f} & \Delta. A . B \ar[d]^{\pi_B} \\
  \Gamma . \sb f A \ar[r]_{f.\ro f}  & \Delta. A\text,}
\end{equation}
the canonical isofibration structure on $\pi_{\sb f B}$ qua map of $\Con_\Ss /
\Gamma$ agrees with the pullback along $f$ of the canonical isofibration
structure on $\pi_B$ qua map of $\Con_\Ss / \Delta$.
\end{Prop}
\begin{proof}
As in the proof of Proposition~\ref{2pullbacklemma}, we will use $\subst$
notation in applications of the Leibniz rule, in order to make clear the
dependent family in which substitution is taking place. Now, to prove the
Proposition, it suffices to show the following. Suppose give an invertible
$2$-cell
\begin{equation}\label{intert}
    \cd[@!C@C-2em]{
      \Lambda \ar[rr]^{k} \ar[dr]_{h} & \rtwocell{d}{\alpha} &
      \Gamma.\sb f A.\sb f B \ar[dl]^{\pi_{\sb f B}} \\ &
      \Gamma.\sb f A
    }
\end{equation}
of $\Con_\Ss / \Gamma$ (i.e., one satisfying $\pi_{\sb f A} \alpha =
\id_{\pi_{\sb f A} h}$). Let us write $\alpha'
\defeq f.\ro f \circ \alpha$ and $k' \defeq f.\ro f.\ro f \circ k$. Then we
must show that
\begin{equation}\label{equalities}
%\begin{aligned}
    s_{\alpha'}  = f.\ro f.\ro f \circ s_{\alpha} \colon \Lambda \to \Delta.A.B
   \quad \text{and} \quad \sigma_{\alpha'} = f.\ro f.\ro f \circ \sigma_{\alpha} \colon s_{\alpha'} \Rightarrow k'\text,
%\end{aligned}
\end{equation}
where we obtain $(s_\alpha, \sigma_\alpha)$ from the canonical isofibration
structure on $\pi_{\sb f B}$, and $(s_{\alpha'}, \sigma_{\alpha'})$ from that
on $\pi_B$. So suppose given a $2$-cell as in~\eqref{intert}, with $h$, $k$ and
$\alpha$ given as follows:
\begin{align*}
x : \Lambda &\t h(x) \defeq (h_1x, h_2 x) : \Gamma.\sb f A \\
x : \Lambda &\t k(x) \defeq (h_1x, k_2 x, k_3x) : \Gamma.\sb f A.\sb f B \\
\text{and} \qquad
x : \Lambda &\t \alpha(x) \defeq (\r h_1 x, \alpha_2 x) : \Id_{\Gamma.\sb f A}\big((h_1x, h_2x), (h_1x, k_2x)\big)\text.
\end{align*}
We first compute the pair $(s_\alpha, \sigma_\alpha)$. The map $s_\alpha \colon
\Lambda \to \Gamma.\sb f A.\sb f B$ is given by
\begin{equation*}
    x : \Lambda \t \big(h_1 x\c h_2 x\c \subst_{[u, v] B(fu, v)}((\r h_1 x, \alpha_2 x), k_3 x)\big) : \Gamma.\sb f A.\sb f B\text;
\end{equation*}
which, by unfolding the inductive description of the $\Id$-elimination rule
given in the proof of Proposition~\ref{idcontexts}, is equal to
\begin{equation}\label{sa}
    x : \Lambda \t \big(h_1 x\c h_2 x\c \subst_{[v] B(fh_1 x, v)}(\alpha_2 x, k_3 x)\big) : \Gamma.\sb f A.\sb f B\text.
\end{equation}
The corresponding $2$-cell $\sigma_\alpha \colon s_\alpha \Rightarrow k$ is now
given by
\begin{equation}\label{siga}
    x : \Lambda \t \big(\r h_1 x\c \alpha_2 x\c \r( \subst_{[v] B(fh_1 x, v)}(\alpha_2 x, k_3 x) )\big) : \Id(s_\alpha x, kx)\text.
\end{equation}
Next we compute the pair $(s_{\alpha'}, \sigma_{\alpha'})$. By the proof of
Proposition~\ref{2pullbacklemma} we have
\begin{align*}
    \alpha'(x) &\defeq \lift{(f.\ro f)}(\r h_1 x, \alpha_2 x) \\ & = (\lift f \r h_1 x, \theta(\r h_1 x, k_2 x) \circ \alpha_2 x)
    \\ &= (\r f h_1 x, \r(k_2 x) \circ \alpha_2 x) = (\r f h_1 x, \alpha_2 x)\text.
\end{align*}
Thus the morphism $s_{\alpha'} \colon \Gamma \to \Delta.A.B$ is given by
\begin{equation*}
    x : \Lambda \t \big(fh_1 x\c h_2 x\c \subst_{[u, v] B(u, v)}((\r f h_1 x, \alpha_2 x), k_3 x)\big) : \Delta.A.B\text;
\end{equation*}
which, by unfolding the description of $\Id$-elimination, is definitionally
equal to
\begin{equation}\label{saa}
    x : \Lambda \t \big(fh_1 x\c h_2 x\c \subst_{[v] B(f h_1 x, v)}(\alpha_2 x, k_3 x)\big) : \Delta.A.B\text.
\end{equation}
The corresponding $2$-cell $\sigma_{\alpha'} \colon s_{\alpha'} \Rightarrow k'$
is now given by
\begin{equation}\label{sigaa}
    x : \Lambda \t \big(\r f h_1 x\c \alpha_2 x\c \r( \subst_{[v] B(fh_1 x, v)}(\alpha_2 x, k_3 x) )\big) : \Id(s_{\alpha'} x, k'x)\text.
\end{equation}
It remains to verify the equalities in~\eqref{equalities}. The first equality
follows immediately from inspection of~\eqref{sa} and~\eqref{saa}. For the
second, we will need a calculation. Suppose given $(x, y, s)$ and $(x, z, t) :
\Lambda.\sb f A.\sb f B$ together with identity proofs
    $p : \Id(y, z)$ and $q : \Id(s, \subst_{[v] B(fx, v)}(p, t))$.
We claim that:
\begin{equation}\label{claim}
    \lift{(f.\ro f.\ro f)}(\r(x), p, q) = (\r(fx), p, q) : \Id_{\Delta.A.B}\big((fx, y, s), (fx, z, t)\big)\text.
\end{equation}
By discrete identity types, it suffices to prove this up to propositional
equality; and by applying $\Id$-elimination twice, first on $p$ and then on
$q$, it suffices for this to show that, given $(x, y, s) : \Gamma.\sb f A.\sb f
B$, we have
\begin{equation*}
    \lift{(f.\ro f.\ro f)}(\r x, \r y, \r k) \approx (\r fx, \r y , \r k) : \Id_{\Delta.A.B}\big((fx, y, s), (fx, y, s)\big)\text.
\end{equation*}
But this follows by the $\Id$-computation rule and the definition of
$\lift{(f.\ro f.\ro f)}$. Thus we have~\eqref{claim} as claimed. We now use
this to affirm the second equality in~\eqref{equalities}. Given $x : \Lambda$,
we have that:
\begin{align*}
    (f.\ro f.\ro f \circ \sigma_\alpha)(x) &= \lift{(f.\ro f.\ro f)}\big(\r h_1 x\c \alpha_2 x\c \r(\subst_{[v] B(fh_1 x, v)}(\alpha_2 x,
k_3 x))\big)\\
&= \big(\r f h_1 x\c \alpha_2 x\c \r(\subst_{[v] B(fh_1 x, v)}(\alpha_2 x,
k_3 x))\big) \\
&= \sigma_{f.\ro f \circ \alpha}(x) = \sigma_{\alpha'}(x)\text. %\qedhere
\end{align*}
\end{proof}

\begin{Rk}
Although it may seem somewhat technical, the result we have just proven is
absolutely crucial for obtaining a sound notion of two-dimensional model.
Without it, our models would not necessarily be sound for the rules expressing
stability of the elimination rules under change of ambient context. It is not
just the identity type rules would be afflicted either: we will see in \S
\ref{intidtypes}--\S \ref{intprodtypes} below that
Proposition~\ref{stabilityiso} is used in verifying the stability under
substitution of \emph{all} the type-theoretic elimination rules. One of the key
issues in giving higher-dimensional and homotopy-theoretic semantics for
intensional type theory will be finding an appropriate counterpart of this
Proposition.
\end{Rk}

\subsection{Comprehension $2$-categories}\label{compr2cat}
We may abstract away from the syntactic investigations of the preceding
sections as follows. We define a \emph{full split comprehension $2$-category}
$\mathbb C$ to be given by the following data: a locally groupoidal
$2$-category $\Con$ with a specified $2$-terminal object; a globally split
$2$-fibration $p \colon \Typ \to \Con$, with $\Typ$ also locally groupoidal;
and a $2$-fully faithful $2$-functor $E \colon \Typ \to \Con^\mathbf 2$
rendering commutative the triangle
\begin{equation*}
  \cd[@!C]{
    \Typ \ar[dr]_p \ar[rr]^{E} & &
    \Con^\mathbf 2 \ar[dl]^{\cod} \\ &
    \Con\text.
  }
\end{equation*}
Moreover, the $2$-functor $E$ should send cartesian morphisms in $\Typ$ to
$2$-pullback squares in $\Con$; should send each object of $\Typ$ to a normal
isofibration in $\Con$; and should satisfy the stability conditions of
Proposition~\ref{stabilityiso}. The preceding developments show that we may
associate a full split comprehension $2$-category to each dependent type theory
$\Ss$ satisfying the rules for identity types in Table~\ref{fig1} and the
discreteness rules of Table~\ref{fig2}. We denote this comprehension
$2$-category by $\mathbb C(\Ss)$, and call it the \emph{classifying
comprehension $2$-category} of $\Ss$.

\section{Categorical models for \mlt: logical aspects}\label{sec3}
\subsection{Identity types}\label{idtypessec}
In this section, we will examine the structure induced on the syntactic
comprehension $2$-category of the previous section by the logical rules of
two-dimensional type theory. Once again we consider a fixed dependently typed
calculus $\Ss$ which we now suppose to admit all of the rules in
Tables~\ref{fig1}, \ref{fig2} and \ref{fig3}. We begin by investigating the
identity types. Given how deeply intertwined these have been with the
construction of the syntactic comprehension $2$-category, it is perhaps
unsurprising that their characterisation is rather intrinsic. It will be given
in terms of the $2$-categorical notion of arrow object. Given a $2$-category
$\tcat K$, an \emph{arrow object} for $X \in \tcat K$ is given by an object $Y
\in \tcat K$ such that $1$-cells into $Y$ correspond naturally to $2$-cells
into~$X$. That is, we have an isomorphism of categories
\begin{equation}\label{isoarrow}
    \tcat K(A, Y) \cong \tcat K(A, X)^\mathbf 2\text,
\end{equation}
$2$-natural in $A$. In particular, under the bijection~\eqref{isoarrow}, the
identity map $\id_Y \colon Y \to Y$ corresponds to a $2$-cell
\begin{equation*}
  \cd{
    Y \ar@/^1em/[rr]^{s} \ar@/_1em/[rr]_{t} \dtwocell{rr}{\kappa} & & X\text;}
\end{equation*}
and $2$-naturality of~\eqref{isoarrow} says that any other such $2$-cell into
$X$ factors uniquely through $\kappa$. In the language of enriched category
theory~\cite{Kelly1982Basic}, an arrow object is a certain kind of (weighted)
limit, namely a \emph{power} (or sometimes \emph{cotensor}) with the category
$\mathbf 2$. For a full treatment of $2$-categorical limits, we refer the
reader to~\cite{Kelly1989Elementary}.

We now introduce a small abuse of notation. Given $\Gamma \in \Con_\Ss$ and $A
\in \Typ_\Ss(\Gamma)$,  we write $\Gamma.A.A$ for the context $\big(\vec x :
\Gamma, y : A(\vec x), z : A(\vec x)\big)$---this rather than the more correct
$\Gamma.A.\pi_A^\ast A$---and write $\pi_{1}$ and $\pi_{2}$ for the context
morphisms $\Gamma.A.A \to \Gamma.A$ projecting onto the first or second copy of
$A$.

\begin{Prop}\label{arrowobj}
For every context $\Gamma$ and type $\Gamma \vdash A\ \ty$ in $\Ss$, the
context $\Gamma.A.A.\Id_A$, together with the projections $\pi_1 \pi_{\Id_A},
\pi_2 \pi_{\Id_A} \colon \Gamma.A.A.\Id_A \to \Gamma.A$, can be made into an
arrow object for $\Gamma.A$ in the slice $2$-category $\Con_\Ss / \Gamma$.
\end{Prop}

\begin{proof}
Let us write $s \defeq \pi_1 \pi_{\Id_A}$ and $t \defeq \pi_2 \pi_{\Id_A}$. We
are to find a $2$-cell
\begin{equation}\label{universal2cell}
  \cd[@!C@C-3em@R+1em]{
    \Gamma.A.A.\Id_A \ar[dr]_{\pi} \ar@/^1em/[rr]^{s} \ar@/_1em/[rr]_{t} \dtwocell{rr}{\kappa} & &
    \Gamma.A \ar[dl]^{\pi_A} \\ &
    \Gamma
  }
\end{equation}
in $\Con_\Ss$ which is over $\Gamma$ in the sense that $\pi_A s = \pi_A t =
\pi$ and $\pi_A \kappa = \id_{\pi}$, and such that any other $2$-cell
\begin{equation}\label{generaluniv}
  \cd[@!C@C-1em]{
    \Lambda \ar[dr]_{h} \ar@/^1em/[rr]^{f} \ar@/_1em/[rr]_{g} \dtwocell{rr}{\alpha} & &
    \Gamma.A \ar[dl]^{\pi_A} \\ &
    \Gamma
  }
\end{equation}
over $\Gamma$ factors through $\kappa$ via a unique morphism $\bar \alpha
\colon \Lambda \to \Gamma.A.A.\Id_A$. The universal property of $\kappa$ also
has a two-dimensional aspect. Suppose we are given a commutative diagram
\begin{equation}\label{2dimup}
    \cd{
        f \ar@2[r]^\beta \ar@2[d]_{\alpha} &
        f' \ar@2[d]^{\alpha'} \\
        g \ar@2[r]_\gamma &
        g'
    }
\end{equation}
of $1$- and $2$-cells $\Lambda \to \Gamma.A$ over $\Gamma$. Then we should be
able to find a unique $2$-cell $\delta \colon \bar \alpha \Rightarrow \bar
{\alpha}' \colon \Lambda \to \Gamma.A.A.\Id_A$ with $\beta = s \delta$ and
$\gamma = t \delta$. We begin by defining $\kappa$ as
in~\eqref{universal2cell}. For this we are required to give a judgement
\begin{equation*}
  \vec x : \Gamma\c y, z : A(\vec x) \c p : \Id(y, z) \t \kappa(\vec x, y, z, p) : \Id(y, z)\text;
\end{equation*}
which we do by taking $\kappa(\vec x, y, z, p) \defeq p$. We now verify the
universal property of $\kappa$. Suppose given an $\alpha$ as
in~\eqref{generaluniv}: then the commutativity conditions $\pi_A f = \pi_A g =
h$ mean that $f$ and $g$ correspond to judgements
\begin{equation*}
  \vec x : \Lambda \t f(\vec x) : A(\vec h \vec x) \qquad \text{and} \qquad
  \vec x : \Lambda \t g(\vec x) : A(\vec h \vec x)\text,
\end{equation*}
whereupon---by Lemma~\ref{pilift}---the condition $\pi_A \alpha = \id_h$ allows
us to view $\alpha$ as a judgement $\vec x : \Lambda \t \alpha(\vec x) : \Id(f
\vec x, g \vec x)$. We now define a morphism $\bar \alpha \colon \Lambda \to
\Gamma.A.A.\Id_A$ by $  \vec x : \Lambda \t (\vec h \vec x, f \vec x, g \vec x,
\alpha \vec x) : \Gamma.A.A.\Id_A$. It is immediate from the definition of
$\kappa$ that $\kappa \bar \alpha = \alpha$, and moreover that if $\kappa m =
\alpha$ for some $m \colon \Lambda \to \Gamma.A.A.\Id_A$ then we have $\bar
\alpha = m$. It still remains to verify the two-dimensional universal property
of $\kappa$. So suppose given $1$- and $2$-cells as in~\eqref{2dimup}. We are
required to define a $2$-cell $\delta \colon \bar \alpha \Rightarrow \bar
{\alpha}' \colon \Lambda \to \Gamma.A.A.\Id_A$ satisfying $s \delta = \beta$
and $t \delta = \gamma$. In order to satisfy these last two requirements,
$\delta$, if it exists, must be given by a judgement
\begin{equation*}
  \vec x : \Lambda \t (\r \vec h \vec x, \beta \vec x, \gamma \vec x, \delta_4 \vec x)
  : \Id_{\Gamma.A.A.\Id_A}\big((\vec h \vec x, f \vec x, g \vec x, \alpha \vec x),
    (\vec h \vec x, f' \vec x, g' \vec x, \alpha' \vec x)\big)
\end{equation*}
for some $\vec x : \Lambda \t \delta_4(\vec x) : \Id_{\Id(f \vec x, g \vec
x)}\big(\alpha \vec x\c (\vec \beta \vec x, \vec \gamma \vec x)^\ast(\alpha'
\vec x)\big)$. By discrete identity types, this is only possible if in fact
$\alpha(\vec x) = (\vec \beta \vec x, \vec \gamma \vec x)^\ast(\alpha' \vec
x)$, whereupon we can take $\delta_4(\vec x) = \r(\vec \alpha \vec x)$. We
claim that in fact $(\vec \beta \vec x, \vec \gamma \vec x)^\ast(\alpha' \vec
x) = (\gamma \vec x)^{-1} \circ (\alpha'\vec x \circ \beta \vec x)$, so that we
will be done if we can show that $\alpha(\vec x) = (\gamma \vec x)^{-1} \circ
(\alpha'\vec x \circ \beta \vec x)$: and this follows from the equation $\gamma
\alpha = \alpha' \beta$ using the groupoid laws for $\Id_A$. It remains only to
prove the claim, which follows from the more general result that
\begin{multline*}
  \vec x : \Gamma\c y, z, y', z' : A(\vec x) \c p : \Id(y, y') \c q : \Id(z, z') \c s : \Id(y', z') \\
  \t (p, q)^\ast(s) = q^{-1} \circ (s \circ p) : \Id(y,z)\text.
\end{multline*}
By discrete identity types, it suffices to prove this up to propositional
equality; and by $\Id$-elimination on $p$ and $q$, it suffices to consider the
case where $y = y'$, $z = z'$, $p = \r(y)$ and $q = \r(z)$, where we have that
$(\r(y), \r(z))^\ast(s) = s = \r(z)^{-1} \circ (s \circ \r(y))$ as required.
\end{proof}

\begin{Prop}[Stability for identity types]\label{bcids} Let $\Gamma$, $\Delta$ be contexts in $\Ss$,
let $f \colon \Gamma \to \Delta$ be a context morphism, and let $\vec x :
\Delta \t B(\vec x) \ \ty$. Then the comparison morphism
\begin{equation*}
\Gamma.\sb f B.\sb f B.\sb {(f.\ro f.\ro f)} {(\Id_B)} \to
\Gamma.\sb f B.\sb f B.\Id_{\sb f B}
\end{equation*}
induced by the universal property of $\Id_{\sb f B}$ is an identity.
\end{Prop}

\begin{proof}
Immediate from the stability of identity types under substitution.
\end{proof}

\subsection{Digression on $2$-categorical adjoints}\label{digression1}
Our characterisation of the remaining type constructors of \mlt\ will be given
in terms of weak $2$-categorical adjoints. We therefore break off at this point
in order to give a brief summary of the \mbox{$2$-categorical} notions
necessary for this characterisation. Let $\tcat K$ be a $2$-category. By a
\emph{retract equivalence} in $\tcat K$, we mean a pair of objects $x, y \in
\tcat K$, a pair of morphisms $i \colon x \to y$ and $p \colon y \to x$
satisfying $pi = \id_x$, and an invertible $2$-cell $\theta \colon \id_y
\Rightarrow ip$ satisfying $\theta i = \id_i$ and $p \theta = \id_p$. In these
circumstances, we may call $i$ an \emph{injective equivalence}---with the
understanding that the extra data $(p, \theta)$ is provided as part of this
assertion---or call $p$ a \emph{surjective equivalence} (with the same
understanding). Given now a $2$-functor $U \colon \tcat K \to \tcat L$ and an
object $x \in \tcat L$, we define a \emph{retract bireflection} of $x$ along
$U$ to be an object $Fx \in \tcat K$ and morphism $\eta_x \colon x \to UFx$
such that for each $y \in \tcat K$, the functor
\begin{equation*}
  \tcat K(Fx, y) \xrightarrow{U_{Fx, y}} \tcat L(UFx, Uy) \xrightarrow{(\thg) \circ \eta_x} \tcat L(x, Uy)
\end{equation*}
is a surjective equivalence of categories. By a \emph{left retract biadjoint}
$F$ for $U$, we mean a choice for every $x \in \tcat L$ of a retract
bireflection $Fx$ of $x$ along $U$. Note that if $F$ is a left retract
biadjoint for $U$, then the assignation $x \mapsto Fx$ will not in general
extend to a $2$-functor $F \colon \tcat L \to \tcat K$; rather, it gives a
\emph{pseudo-functor}, which preserves identities and composition only up to
invertible $2$-cells. Likewise, the maps $\eta_x \colon x \to UFx$ do not
provide components of a $2$-natural transformation $\eta \colon \id_\tcat L
\Rightarrow UF$ but merely of a \emph{pseudo-natural} transformation, whose
naturality squares commute only up to invertible $2$-cells. We could give a
definition of left retract biadjoint in terms of a pseudo-functor $\tcat K \to
\tcat L$ and unit and counit transformations $\eta$ and $\epsilon$ satisfying
weakened versions of the triangle laws
(see~\cite[Section~1]{Street1974Fibrations} for the details); but the above
description is both simpler and, as we will see, closer to the type theory. In
fact, the above definitions admit a further simplification, using the
observation that the surjective equivalences of categories are precisely those
functors $F \colon \C \to \D$ which are fully faithful and whose object
function $\ob F \colon \ob \C \to \ob \D$ is a split epimorphism:

\begin{Prop}\label{characterisation}
To give a retract bireflection of $x \in \tcat L$ along $U \colon \tcat K \to
\tcat L$ is to give an object $Fx \in \tcat K$ and map $\eta_x \colon x \to
UFx$, together with, for each $f \colon x \to Uy$ in $\tcat L$, a choice of map
$\bar f \colon Fx \to y$ in $\tcat K$ satisfying $U\bar f \circ \eta_x = f$;
all subject to the requirement that, for every $h, k \colon Fx \to y$ in $\tcat
K$ and every $\alpha \colon Uh \circ \eta_x \Rightarrow Uk \circ \eta_x$ in
$\tcat L$, there is a unique $\bar \alpha \colon h \Rightarrow k$ with $U\bar
\alpha \circ \eta_x = \alpha$.
\end{Prop}

Given a $2$-functor $U \colon \tcat K \to \tcat L$ and $x \in \tcat L$ as
before, we have the dual notion of \emph{retract bicoreflection} of $x$ along
$U$: this being given by an object $Gx \in \tcat K$, together with a morphism
$\epsilon_x \colon UGx \to x$ such that for each $y \in \tcat K$, the functor
\begin{equation*}
  \tcat K(y, Gx) \xrightarrow{U_{y, Gx}} \tcat L(Uy, UGx) \xrightarrow{\epsilon_x \circ (\thg)} \tcat L(Uy, x)
\end{equation*}
is a surjective (not injective!) equivalence of categories. Now a \emph{right
retract biadjoint} for $U$ is of course given by a choice for every $x \in
\tcat L$ of a retract bicoreflection along $U$. As before, we have an
elementary characterisation of retract bicoreflections:

\begin{Prop}
To give a retract bicoreflection of $x \in \tcat L$ along $U \colon \tcat K \to
\tcat L$ is to give an object $Gx \in \tcat K$ and map $\epsilon_x \colon UGx
\to x$, together with, for each $f \colon Uy \to x$ in $\tcat L$, a choice of
map $\bar f \colon y \to Gx$ in $\tcat K$ satisfying $\epsilon_x \circ U\bar f
 = f$; all subject to the requirement that, for every $h, k \colon
y \to Gx$ in $\tcat K$ and every $\alpha \colon \epsilon_x \circ Uh \Rightarrow
\epsilon_x \circ Uk$ in $\tcat L$, there is a unique $\bar \alpha \colon h
\Rightarrow k$ with $\epsilon_x \circ U\bar \alpha  = \alpha$.
\end{Prop}

\subsection{Unit types}
Our first application of the $2$-categorical adjoint notions developed above
will be to the unit types of~$\Ss$---which we recall is an arbitrary dependent
type theory admitting all the rules listed in
Tables~\ref{fig1},~\ref{fig2}~and~\ref{fig3}. In the following result, we
denote by $E(\Gamma) \colon \Typ_\Ss(\Gamma) \to \Con_\Ss / \Gamma$ the
$2$-functor obtained by restricting $E \colon \Typ_\Ss \to \Con_\Ss$ to the
fibre over $\Gamma \in \Con_\Ss$.

\begin{Prop}\label{unitytpes}
For each context $\Gamma$ of $\Ss$, the object $\mathbf 1_\Gamma \in
\Typ_\Ss(\Gamma)$ given by $\Gamma \t \mathbf 1 \ \ty$ provides a retract
bireflection of $\id_\Gamma \colon \Gamma \to \Gamma$ along the $2$-functor
$E(\Gamma) \colon \Typ_\Ss(\Gamma) \to \Con_\Ss / \Gamma$.
\end{Prop}

\begin{proof}
The unit of the bireflection $\eta_\Gamma \colon \Gamma \to \Gamma.\mathbf
1_\Gamma$ (over $\Gamma$) is given by the judgement $\vec x : \Gamma \t \star :
\mathbf 1$. Given now a morphism $f \colon \Gamma \to \Gamma.A$ over
$\Gamma$---which is equally well a judgement $\vec x : \Gamma \t f(\vec x) :
A(\vec x)$---we obtain a factorisation $\bar f \colon \Gamma.\mathbf 1_\Gamma
\to \Gamma.A$ over $\Gamma$ by $\mathbf 1$-elimination, taking $\bar f$ to be
the term \mbox{$\vec x : \Gamma\c z : \mathbf 1 \t \mathrm U_{f(\vec x)}(z) :
A(\vec x)$}. That this satisies $\bar f \eta_\Gamma = f$ is now precisely the
computation rule $\vec x : \Gamma \t \mathrm U_{f(\vec x)}(\star) = f(\vec x)$.
It remains to check that for maps $h, k \colon \Gamma.\mathbf 1 \to \Gamma.A$
over $\Gamma$, every 2-cell $\alpha \colon h \eta_\Gamma \Rightarrow k
\eta_\Gamma$ over $\Gamma$ is of the form $\bar \alpha \eta_\Gamma$ for a
unique $\bar \alpha \colon h \Rightarrow k$. Now, to give $h$, $k$ and $\alpha$
is to give judgements
\begin{align*}
  \vec x : \Gamma\c z : \mathbf 1 &\t h(\vec x, z) : A(\vec x)\\
  \vec x : \Gamma\c z : \mathbf 1 &\t k(\vec x, z) : A(\vec x)\\
  \vec x : \Gamma &\t \alpha(\vec x) : \Id\big(h(\vec x, \star), k(\vec x, \star)\big)\text;
\end{align*}
from which we must determine $\vec x : \Gamma\c z : \mathbf 1 \t \bar
\alpha(\vec x, z) : \Id\big(h(\vec x, z), k(\vec x, z)\big)$. We do this by
$\mathbf 1$-elimination, taking $\bar \alpha(\vec x, z) := \mathrm
U_{\alpha(\vec x)}(z)$. The equality $\bar \alpha \eta_\Gamma = \alpha$ now
follows from the $\mathbf 1$-computation rule. It remains to check uniqueness
of $\bar \alpha$. So suppose we are given $\vec x : \Gamma\c z : \mathbf 1 \t
\beta(\vec x, z) : \Id\big(h(\vec x, z), k(\vec x, z)\big)$ satisfying
$\beta(\vec x, \star) = \alpha(\vec x)$. We must show that $\beta(\vec x, z) =
\bar \alpha(\vec x, z)$. By discrete identity types, it suffices to show this
up to propositional equality; and by $\mathbf 1$-elimination, this only in the
case where $z = \star$, for which we have that $\beta(\vec x, \star) =
\alpha(\vec x) = \bar \alpha(\vec x, \star)$ as required.
\end{proof}

\begin{Prop}[Stability for unit types]\label{bcunits} For each $k \colon \Gamma \to \Delta$ in $\Con$,
we have $k^\ast(\mathbf 1_\Delta) = \mathbf 1_\Gamma$; we have $\eta_\Gamma =
k^\ast(\eta_\Delta) \colon \Gamma \to \Gamma.\mathbf 1_\Gamma$; and for each $f
\colon \Delta \to \Delta.B$ over $\Delta$, have $k^\ast(\bar f) =
\overline{k^\ast(f)} \colon \Gamma.\mathbf 1_\Gamma \to \Gamma.k^\ast B$.
\end{Prop}

\begin{proof}
By the stability of unit types under substitution.
\end{proof}

\begin{Rk}
Note carefully what the previous result does \emph{not} say: it does not say
that for a context morphism $k \colon \Gamma \to \Delta$, the comparison map
$\mathbf 1_\Gamma \to k^\ast \mathbf 1_\Delta$ of $\Typ_\Ss(\Gamma)$ is an
identity; indeed, this map will in general only be \emph{isomorphic} to the
identity, since it corresponds to the judgement $\vec x : \Gamma\c z : \mathbf
1 \t \mathrm U_{\star}(z) : \mathbf 1$.
\end{Rk}

\subsection{Dependent sum types}
We next consider the dependent sum types.

\begin{Prop}\label{depsums1}
For each context $\Gamma$ and type $\Gamma \t A \ \ty$ of $\Ss$, the
\mbox{$2$-functor} $\Delta_A := \Typ_\Ss(\pi_A) \colon \Typ_\Ss(\Gamma) \to
\Typ_\Ss(\Gamma.A)$ has a left retract biadjoint $\Sigma_A$.
\end{Prop}

\begin{proof}
We must provide, for each $B \in \Typ_\Ss(\Gamma.A)$ a retract bireflection
$\Sigma_A(B)$ of $B$ along $\Delta_A$. So we take $\Sigma_A(B) \in
\Typ_\Ss(\Gamma)$ to be given by the judgement \mbox{$\Gamma \t \Sigma(A, B)\
\ty$} (where for readability we suppress explicit mention of dependencies on
the variables in $\Gamma$); and the unit map $\eta \colon B \to \Delta_A
\Sigma_A(B)$ of $\Typ_\Ss(\Gamma.A)$ to be given by the judgement $\Gamma \c y
: A \c z : B(y) \t \spn{y, z} : \Sigma(A, B)$. Now given a type $C \in
\Typ_\Ss(\Gamma)$ and a map $f \colon B \to \Delta_A C$ of
$\Typ_\Ss(\Gamma.A)$, we must provide a morphism $\bar f \colon \Sigma_A(B) \to
C$ of $\Typ_\Ss(\C)$ satisfying $\Delta_A(\bar f) \circ \eta = f$. But to give
$f$ is to give a judgement $\Gamma\c y : A\c z : B(y) \t f(y, z) : C$, whilst
to give $\bar f$ is to give a judgement $\Gamma\c s : \Sigma(A, B) \t \bar f(s)
: C$. Thus by using $\Sigma$-elimination we may define $\bar f(s)
\defeq \mathrm E_f(s)$. The equality $\Delta_A(\bar f) \circ \eta = f$ follows by
the \mbox{$\Sigma$-computation} rule. It remains to show, given two morphisms
$h, k \colon \Sigma_A(B) \to D$ in $\Typ_\Ss(\Gamma)$, that each $2$-cell
$\alpha \colon \Delta_A(h) \circ \eta \Rightarrow \Delta_A(k) \circ \eta$ is of
the form $\Delta_A(\bar \alpha) \circ \eta$ for a unique $\bar \alpha \colon h
\Rightarrow k$. This follows by an argument analogous to that given in the
proof of Proposition~\ref{unitytpes}.
\end{proof}

Whilst Proposition~\ref{depsums1} is very natural from a categorical
perspective, it fails to capture the full strength of the elimination rule for
dependent sums (even though it requires the full strength of that elimination
rule in its proof). In order to do this, we need the following result:

\begin{Prop}\label{depsums2}
Suppose given a context $\Gamma$ in $\Ss$ and types $\Gamma \t A \ \ty$ and
$\Gamma\c x : A \t B(x) \ \ty$ in $\Ss$, and consider the morphism
\begin{equation}\label{canonn}
  \cd{
    \Gamma.A.B \ar[r]^-{i} \ar[d]_{\pi_B} &
    \Gamma.\Sigma_A(B) \ar[d]^{\pi_{\Sigma_A(B)}} \\
    \Gamma.A \ar[r]_{\pi_A} &
    \Gamma
  }
\end{equation}
in $\Con_\Ss^\mathbf 2$ corresponding to the unit morphism $\eta \colon B \to
\Delta_A \Sigma_A(B)$ in $\Typ_\Ss(\Gamma.A)$. The map $i$ appearing in this
diagram is an injective equivalence in $\Con_\Ss / \Gamma$.
\end{Prop}

\begin{proof}
We construct a pseudoinverse retraction for $i$ over $\Gamma$ as follows. The
map $p \colon \Gamma.\Sigma_A(B) \to \Gamma.A.B$ over $\Gamma$ is given by the
projections out of the sum:
\begin{align*}
  \Gamma\c s : \Sigma(A, B) &\t s.1 : A\\
  \Gamma\c s : \Sigma(A, B) &\t s.2 : B(s.1)
\end{align*}
(where again, we suppress explicit mention of the dependency on $\Gamma$). We
define these by $\Sigma$-elimin\-ation on $s$, the first being given by $s.1
\defeq \mathrm E_{[y, z] y}(s)$ and the second by $s.2 \defeq \mathrm E_{[y, z]
z}(s)$. The equality $pi = \id_{\Gamma.A.B}$ follows from the
$\Sigma$-computation rule. We must now give a $2$-cell $\theta \colon
\id_{\Gamma.\Sigma_A(B)} \Rightarrow ip$; which is equally well a judgement
$\Gamma\c s : \Sigma(A, B) \t \theta(s) : \Id(s, \spn{s.1, s.2})$. By
\mbox{$\Sigma$-elimination} on $s$, it suffices to define $\theta$ when $s =
\spn{y, z}$; whereupon we have $\spn{s.1, s.2} = \spn{\spn{y, z}.1, \spn{y,
z}.2} = \spn{y, z}$ so that we can take $\theta(\spn{y, z}) = \r(\spn{y, z})$.
The equality $\theta i = \id_i$ now follows by the $\Sigma$-computation rule;
and it remains only to verify that $p \theta = \id_p$. Now, $p \theta$
corresponds to the judgement
\begin{equation*}
  \Gamma\c s : \Sigma(A, B) \t \lift{\vec p}(\theta(s)) : \Id_{\Gamma.A.B}\big((s.1, s.2), (s.1, s.2)\big)\text;
\end{equation*}
and we must show that in fact $\lift{\vec p}(\theta(s)) = \r(p(s))$. By
discrete identity types, it suffices to show this up to propositional equality;
and by $\Sigma$-elimination, this only when $s = \spn{y, z}$. But we calculate
that $\lift{\vec p}(\theta(\spn{y, z})) = \lift{\vec p}(\r(\spn{y, z}))
  = \r(\vec p(\spn{y, z}))$
as required.
\end{proof}

%\begin{Rk}
%It is in fact the case that Proposition~\ref{depsums2}, when combined with the
%structure detailed in Proposition~\ref{fibrmain}, subsumes
%Proposition~\ref{depsums1}: we will see how in Section~\ref{intsumtypes}.
%\end{Rk}

\begin{Prop}[Stability for dependent sums]\label{bcsums} Given $k \colon \Gamma \to \Lambda$ in $\Con_\Ss$, $A \in \Typ(\Lambda)$ and $B \in \Typ(\Lambda.A)$, we have that
$k^\ast(\Sigma_A(B)) = \Sigma_{k^\ast A}(k^\ast B)$; that $k^\ast(\eta_{A,B}) =
\eta_{k^\ast A, k^\ast B}$; and for each $f \colon B \to \Delta_A C$ in
$\Typ_\Ss(\Lambda.A)$, that $k^\ast \bar f = \overline{k^\ast f} \colon
\Sigma_{\sb k A}(\sb k B) \to \sb k C$. Moreover, reindexing along $k$ sends
the injective equivalence structure on $i_{A, B}$ to the injective equivalence
structure on $i_{\sb k A, \sb k B}$.
\end{Prop}

\begin{proof}
By the stability of dependent sum types under substitution.
\end{proof}

\subsection{Dependent product types}\label{depprodssec}
Finally, we turn to the categorical characterisation of dependent product types
in $\Ss$.

\begin{Prop}\label{depprods1}
For each context $\Gamma$ and type $\Gamma \t A \ \ty$ of $\Ss$, the weakening
$2$-functor $\Delta_A \colon \Typ_\Ss(\Gamma) \to \Typ_\Ss(\Gamma.A)$ has a
right retract biadjoint $\Pi_A$.
\end{Prop}

\begin{proof}
Once again, we suppress explicit mention of dependencies on the variables in
$\Gamma$. We must provide, for each $B \in \Typ_\Ss(\Gamma.A)$ a retract
bicoreflection $\Pi_A(B)$ of $B$ along $\Delta_A$. For this we take $\Pi_A(B)
\in \Typ_\Ss(\Gamma)$ to be given by the judgement \mbox{$\Gamma \t \Pi(A, B)\
\ty$}; and the counit map $\epsilon \colon \Delta_A \Pi_A(B) \to B$ of
$\Typ_\Ss(\Gamma.A)$ to be given by the judgement $\Gamma \c m : \Pi(A, B) \c y
: A \t m \cdot y : B(y)$. Now given a type $C \in \Typ_\Ss(\Gamma)$ and a map
$f \colon \Delta_A C \to B$ of $\Typ_\Ss(\Gamma.A)$, we are required to provide
a morphism $\bar f \colon C \to \Pi_A(B)$ of $\Typ_\Ss(\C)$ satisfying
$\epsilon \circ \Delta_A(\bar f) = f$. So if $f$ is the judgement $\Gamma\c y :
C\c z : A \t f(y, z) : B(y)$, we take $\bar f$ to be the judgement $\Gamma\c y
: C \t \lambda z.\, f(y, z) : \Pi(A, B)$. The equality $\epsilon \circ
\Delta_A(\bar f) = f$ follows by the $\beta$-rule.

It remains to show, given two morphisms $h, k \colon D \to \Pi_A(B)$ in
$\Typ_\Ss(\Gamma)$, that each $2$-cell $\alpha \colon \epsilon \circ
\Delta_A(h) \Rightarrow \epsilon \circ \Delta_A(k)$ can be written in the form
$\epsilon \circ \Delta_A(\bar \alpha)$ for a unique $\bar \alpha \colon h
\Rightarrow k$. It is here that we will make crucial use of function
extensionality. So, to give $h$, $k$ and $\alpha$ is to give judgements
\mbox{$\Gamma, C \t h : \Pi(A, B)$}; \mbox{$\Gamma, C \t k : \Pi(A, B)$}; and
\mbox{$\Gamma, C\c z : A \t \alpha(z) : \Id(h \cdot z, k \cdot z)$} (where we
now suppress explicit mention of the dependency on $C$) and so we may define
the $2$-cell $\bar \alpha \colon h \Rightarrow k$ by applying the rule
\mbox{$\Pi$\textsc{-ext}} of Table~\ref{fig3} to obtain the judgement $\Gamma,
C \t \mathsf{ext}(h, k, \alpha) : \Id(h, k)$. We must now check that $\epsilon
\circ \Delta_A(\bar \alpha) = \alpha$. Recall from~\S\ref{funextsec1} the
operation
\begin{equation*}
\inferrule{
  m, n : \Pi(A, B) \\ p : \Id(m, n) \\ a : A}
  {p \ast a : \Id(m \cdot a, n \cdot a)}
\end{equation*}
given by $p \ast a \defeq \J_{[x]\r(x \cdot a)}(m, n, p)$. It is easy to see
that $\ast$ is just the lifting of $\epsilon$ to identity types; so that
$\epsilon \circ \Delta_A(\bar \alpha)$ corresponds to the judgement
\begin{equation*}
  \Gamma\c C\c z : A \t \mathsf{ext}(h, k, \alpha) \ast z : \Id(h \cdot z\c k \cdot z)\text.
\end{equation*}
But by the rule \textsc{$\Pi$-ext-app} of Table~\ref{fig3}, we have that
$\mathsf{ext}(h, k, \alpha) \ast z = \alpha(z)$ as required. It remains to
check uniqueness of $\bar \alpha$. So suppose that we are given $\Gamma\c C \t
\beta : \Id(h, k)$ satisfying $\beta \ast z = \alpha(z)$: we must show that
$\beta = \bar \alpha$. Now, because $\beta \ast z = \alpha(z) = \bar \alpha
\ast z$, we have that
\begin{equation*}
  \Gamma\c C \c z : A \t \mathsf{ext}\big(h, k, [z]\, \beta \ast z\big)
  = \mathsf{ext}\big(h, k, [z]\, \bar \alpha \ast z\big) : \Id(h, k)\text.
\end{equation*}
Thus we will be done if we can show that
\begin{equation*}
  \Gamma\c C \c m, n : \Pi(A, B)\c k : \Id(m, n) \t
  \mathsf{ext}(m, n, [z]\,k \ast z) = k: \Id(m, n)
\end{equation*}
holds. By discrete identity types, it suffices to do this up to propositional
equality; and by $\Id$-elimination, this only in the case where $m = n$ and
\mbox{$k = \r(m)$}, so that we will be done if we can show that
\begin{equation*}
  \Gamma\c C \c m : \Pi(A, B) \t
  \mathsf{ext}(m, m, [z]\, \r(m \cdot z)) \approx \r(m) : \Id(m, m)
\end{equation*}
holds. But this follows immediately from the rule \textsc{$\Pi$-ext-comp}.
\end{proof}

\begin{Prop}[Stability for dependent products]\label{bcprods}
Given $k \colon \Gamma \to \Lambda$ in $\Con_\Ss$, $A \in \Typ(\Lambda)$ and $B
\in \Typ(\Lambda.A)$, we have that $k^\ast(\Pi_A(B)) = \Pi_{k^\ast A}(k^\ast
B)$; that $k^\ast(\epsilon_{A,B}) = \epsilon_{k^\ast A, k^\ast B}$; and for
each $f \colon \Delta_A C \to B$ in $\Typ_\Ss(\Lambda.A)$, that $k^\ast \bar f
= \overline{k^\ast f} \colon \sb k C \to \Pi_{\sb k A}(\sb k B)$.
\end{Prop}

\begin{proof}
By the stability of dependent product types under substitution.
\end{proof}

\subsection{Models of two-dimensional type theory}\label{model2dimtt}
We now abstract away from the preceding results as follows.
\begin{Defn}\label{twocompcat}
Let there be given a full split comprehension $2$-category $\mathbb C = (p
\colon \Typ \to \Con\c E \colon \Con \to \Typ^\mathbf 2)$, in the sense
of~\S\ref{compr2cat}. Then:
\begin{itemize}
\item We say that $\mathbb C$ has \emph{equality} if, for every $\Gamma \in
    \Con$ and $A \in \Typ(\Gamma)$, there is an object $\Id_A \in
    \Typ(\Gamma.A.A)$ such that $\Gamma.A.A.\Id_A$, together with its two
    projections onto $\Gamma.A$, underlies an arrow object for $\Gamma.A$
    in $\Con / \Gamma$; and these arrow objects satisfy the stability
    properties of Proposition~\ref{bcids}.
\item We say that $\mathbb C$ has \emph{units} if, for every $\Gamma \in
    \Con$, the map $\id_\Gamma \colon \Gamma \to \Gamma$ admits a retract
    bireflection $\mathbf 1_\Gamma$ along $E(\Gamma) \colon \Typ(\Gamma)
    \to \Con_\Ss / \Gamma$; and these bireflections satisfy the stability
    properties of Proposition~\ref{bcunits}.
\item We say that $\mathbb C$ has \emph{sums} if, for every $\Gamma \in
    \Con$ and $A \in \Typ(\Gamma)$, the \mbox{$2$-functor} $\Delta_A :=
    \Typ(\pi_A) \colon \Typ(\Gamma) \to \Typ(\Gamma.A)$ admits a retract
    left biadjoint $\Sigma_A$; and these biadjoints satisfy the conditions
    of Proposition~\ref{depsums2} and the stability properties of
    Proposition~\ref{bcsums}.
\item We say that $\mathbb C$ has \emph{products} if, for every $\Gamma \in
    \Con$ and $A \in \Typ(\Gamma)$, the \mbox{$2$-functor} $\Delta_A \colon
    \Typ(\Gamma) \to \Typ(\Gamma.A)$ admits a retract right biadjoint
    $\Pi_A$; and these biadjoints satisfy the stability properties of
    Proposition~\ref{bcprods}.
\item We say that $\mathbb C$ is a \emph{model of two-dimensional type
    theory} if it has equality, units, sums and products.
\end{itemize}
\end{Defn}
Thus, the results of this section can be summarised by saying that, for any
dependent type theory $\Ss$ satisfying the rules of Tables~\ref{fig1},
\ref{fig2} and~\ref{fig3}, the classifying comprehension $2$-category $\mathbb
C(\Ss)$ is a model of two-dimensional type theory. With an eye on future
applications, we end this Section by gathering together in one place a list of
the structure required for a two-dimensional model of type theory.

\begin{Defn}
A two-dimensional model of type theory $\mathbb C$ is given by:
\begin{itemize}
\item A locally groupoidal $2$-category $\Con$ of \emph{contexts}, with a
    specified $2$-terminal object;
\item A locally groupoidal $2$-category $\Typ$ of \emph{types-in-context};
\item A globally split $2$-fibration $p \colon \Typ \to \Con$ in the sense
    of Definition~\ref{2fibdef}. Spelling this out, this means that $p$ is
    a $2$-functor such that:
    \begin{enumerate}[(i)]
\item The underlying ordinary functor of $p$ is a split fibration of
    categories;
\item For every cartesian 1-cell $f \colon y \to z$ of $\Con$ and every
    $2$-cell $\alpha \colon g \Rightarrow h \colon x \to z$ of $\Typ$,
    any factorisation of $p(\alpha)$ through $p(f)$ may be lifted
    uniquely to a factorisation of $\alpha$ through $f$.
\item For each $x, y \in \Typ$, the induced functor $p_{x, y} \colon
    \Typ(x, y) \to \Con(px, py)$ is a fibration of groupoids.
\end{enumerate}
(Note that condition (iv) of Definition~\ref{2fibdef} is automatically
satisfied since every fibre category is a groupoid).
\item A \emph{comprehension} $2$-functor $E \colon \Typ \to \Con^\mathbf
    2$, equipped with:
\begin{enumerate}[(i)]
\item For each object $A \in \Typ$, a normal isofibration structure on
    $E(A)$ in the sense of Definition~\ref{isofibdefn}.
\end{enumerate}
\noindent and satisfying the following properties:
\begin{enumerate}[(i)]
\item $\cod \circ E = p$;
\item $E$ is $2$-fully faithful (i.e., an isomorphism on
    hom-groupoids);
\item $E$ sends cartesian morphisms of $\Typ$ to $2$-pullback squares
    in $\Con$;
\item The normal isofibration structures picked out by $E$ have the
    stability properties of Proposition~\ref{stabilityiso}.
\end{enumerate}
\end{itemize}
In describing the remaining, logical, structure, we use freely the conventions
of Notation~\ref{notconv}.
\begin{itemize}
\item For every $\Gamma \in \Con$ and $A \in \Typ(\Gamma)$, there is given
    an object $\Id_A \in \Typ(\Gamma.A.A)$ and a $2$-cell $\kappa \colon
    \pi_1 \Rightarrow \pi_2 \colon \Gamma.A.A.\Id_A \to \Gamma.A$ over
    $\Gamma$ which together provide an arrow object (in the sense
    of~\S\ref{idtypessec}) for $\Gamma.A$ in $\Con / \Gamma$.
\item For every $\Gamma \in \Con$, there is given a retract bireflection
    (in the sense of Proposition~\ref{characterisation}) $\mathbf 1_A$ of
    the object $\id_\Gamma \colon \Gamma \to \Gamma$ along $E(\Gamma)
    \colon \Typ(\Gamma) \to \Con_\Ss / \Gamma$.
\item For every $\Gamma \in \Con$ and $A \in \Typ(\Gamma)$, there are given
    both left and right retract biadjoints (in the sense
    of~\S\ref{digression1}) $\Sigma_A$ and $\Pi_A$ for $\Typ(\pi_A) \colon
    \Typ(\Gamma) \to \Typ(\Gamma.A)$.
\item For every $\Gamma \in \Con$, $A \in \Typ(\Gamma)$ and $B \in
    \Typ(\Gamma.A)$, there is given a choice of injective equivalence
    structure on the canonical morphism $i \colon \Gamma.A.B \to
    \Gamma.\Sigma_A B$ defined as in~\eqref{canonn}.
\item The above structures satisfy the stability properties listed in
    Propositions \ref{bcids}, \ref{bcunits}, \ref{bcsums} and
    \ref{bcprods}.
\end{itemize}
\end{Defn}

\section{The internal language of a two-dimensional model}\label{sec4}
\subsection{$2$-categorical lifting properties}\label{2catliftingprops}
In this Section, we prove a converse to the results of the previous two
Sections. Given a model $\mathbb C$ of two-dimensional type theory, we will
construct from it a dependent type theory $\Ss(\mathbb C)$ admitting the rules
of Tables~\ref{fig1},~\ref{fig2}~and~\ref{fig3}. We call this type theory the
\emph{internal language} of $\mathbb C$. The key to doing this will be to give
semantic analogues in $\mathbb C$ of each of the logical rules of \mlt. In
giving analogues of the elimination rules, we will make use of the
$2$-categorical lifting property described in Proposition~\ref{elimrules1}
below. This is again very much in the spirit of~\cite{Awodey2008Homotopy},
since this is fundamentally a result about the weak factorisation system
(injective equivalences, normal isofibrations) described in
Remark~\ref{linkage}: or rather, about an algebraic presentation of this weak
factorisation system in the style of~\cite{Grandis2006Natural}.
\begin{Prop}\label{elimrules1}
Suppose given a $2$-category $\tcat K$ and a square
\begin{equation}\label{2catsquare}
    \cd{
      A \ar[r]^f \ar[d]_i &
      C \ar[d]^p \\
      B \ar[r]_g &
      D
    }
\end{equation}
where $i$ carries the structure of an injective equivalence (cf.~\S
\ref{digression1}) and $p$ that of a normal isofibration
(cf.~Definition~\ref{isofibdefn}). From this data we can determine a canonical
diagonal filler $j \colon B \to C$ satisfying $pj = g$ and $ji = f$.
\end{Prop}
\begin{proof}
The injective equivalence structure on $i$ is given by a morphism \mbox{$k
\colon B \to A$} satisfying $ki = \id_A$ and an invertible $2$-cell $\theta
\colon \id_B \Rightarrow ik$ satisfying $\theta i = \id_i$ and $k \theta =
\id_k$. Thus we have an invertible $2$-cell
\begin{equation*}
    \cd[@!C]{
      B \ar[rr]^{fk} \ar[dr]_{g} & \rtwocell{d}{g \theta} &
      C \ar[dl]^p \\ &
      D\text,
    }
\end{equation*}
and so from the isofibration structure on $p$ we obtain a map \mbox{$j
\defeq s_{g\theta} \colon B \to C$} satisfying $pj
= g$. It remains to show that $ji = f$. By the definition of isofibration, we
have $ji = s_{g \theta} \circ i = s_{g \theta i}$; and since $s_{g \theta i} =
s_{g (\id_i)} = s_{\id_{gi}} = s_{\id_{pf}}$, we deduce by normality that $ji =
s_{\id_{pf}} = f$ as required.
\end{proof}
We now show that the liftings of the previous Proposition are stable under
pullback in a suitable sense. Note that in order for this to make sense, it is
crucial that Proposition~\ref{elimrules1} gives us a \emph{choice} of filler
for each diagram like~\eqref{2catsquare}.

\begin{Prop}\label{elimrules2}
Suppose given a morphism $h \colon X \to Y$ in a $2$-category $\tcat K$,
together with a diagram like~\eqref{2catsquare} in the slice $\tcat K / Y$.
Suppose that we are able to form the $2$-pullback of this diagram along $h$,
yielding a diagram
\begin{equation}\label{2catsquare2}
    \cd{
      h^\ast A \ar[r]^{h^\ast f} \ar[d]_{h^\ast i} &
      h^\ast C \ar[d]^{h^\ast p} \\
      h^\ast B \ar[r]_{h^\ast g} &
      h^\ast D
    }
\end{equation}
in $\tcat K / X$. Then the pullback of the canonical filler
for~\eqref{2catsquare} along $h$ is equal to the canonical filler
for~\eqref{2catsquare2}, where the injective equivalence structure on $h^\ast
i$ and the isofibration structure on $h^\ast p$ are those induced by pullback.
\end{Prop}
\begin{proof}
Let us first make clear what the induced structures on $h^\ast i$ and $h^\ast
p$ look like. The injective equivalence data for $h^\ast i$ is simply given by
applying $h^\ast$ to the corresponding data for $i$. The normal isofibration
structure on $h^\ast p$ is given as follows. Let us write $h_! \colon \tcat K /
X \to \tcat K / Y$ for the $2$-functor given by postcomposition with $h$. For
any $V \in \tcat K / Y$ whose  $2$-pullback $h^\ast V$ along $h$ exists, we
have $2$-natural bijections of categories
\begin{equation}\label{bijectionofcats}
    \tcat K / Y(h_! U, V) \cong \tcat K / X(U, h^\ast V)\text.
\end{equation}
In particular, we have bijections between diagrams of the following two forms:
\begin{equation}\label{isodiag}
    \cd[@!C@C-1.5em]{
      W \ar[rr]^{g} \ar[dr]_{f} & \rtwocell{d}{\alpha} &
      h^\ast C \ar[dl]^{h^\ast p} \\ &
      h^\ast D
    } \qquad \leftrightarrow \qquad
    \cd[@!C@C-1.5em]{
      h_! W \ar[rr]^{\bar g} \ar[dr]_{\bar f} & \rtwocell{d}{\bar \alpha} &
      C \ar[dl]^{p} \\ &
      D\text.
    }
\end{equation}
So given an $\alpha$ as on the left of~\eqref{isodiag}, we obtain a lifting for
it by first transposing to obtain a $2$-cell $\bar \alpha$ as on the right
of~\eqref{isodiag}. We then apply the isofibration structure of $p$ to obtain
$s_{\bar \alpha} \colon h_! W \to C$ and $\sigma_{\bar \alpha} \colon s_{\bar
\alpha} \Rightarrow \bar g$; and finally, we transpose back
using~\eqref{bijectionofcats} to obtain $s_\alpha \colon W \to h^\ast C$ and
$\sigma_\alpha \colon s_\alpha \Rightarrow g$. Now, consider the case where
$\alpha$ in~\eqref{isodiag} is itself of the form $h^\ast \beta$ for some
$\beta \colon u \Rightarrow pv \colon W \to D$ in $\tcat K / Y$. When this is
so, the corresponding $\bar \alpha$ is, by naturality, equal to $\beta \circ
\epsilon_W$, where $\epsilon_W \colon h_! h^\ast W \to W$ is the transpose of
$\id_{h^\ast W}$ under the bijection~\eqref{bijectionofcats}. It follows from
the definition of isofibration that $s_{\bar \alpha} = s_{\beta \circ
\epsilon_W} = s_\beta \circ \epsilon_W$ and likewise $\sigma_{\bar \alpha} =
\sigma_\beta \circ \epsilon_W$; whereupon transposing
under~\eqref{bijectionofcats} and using naturality, we have $s_{h^\ast \beta} =
h^\ast(s_\beta)$ and $\sigma_{h^\ast \beta} = h^\ast(\sigma_\beta)$. Now,
according to Proposition~\ref{elimrules1}, the canonical filler
for~\eqref{2catsquare2} is given by $s_{(h^\ast g)(h^\ast \theta)} =
s_{h^\ast(g \theta)}$; and by the above argument this is equal to $h^\ast(s_{g
\theta})$, which is precisely the pullback along $h$ of the canonical filler
for~\eqref{2catsquare}, as required.
\end{proof}

\subsection{Identity types}\label{intidtypes}
For the rest of the section, we fix a model of two-dimensional type theory
$\mathbb C$. We are going to give semantic analogues of each of the logical
constructors of \mlt\ in $\mathbb C$. We start with the identity types.

\subsubsection{Formation rule}
Given $\Gamma \in \Con$ and $A \in \Typ(\Gamma)$, we define the \emph{semantic
identity type} on $A$ to be the object $\Id_A \in \Typ(\Gamma.A.A)$ whose
existence is assured by Definition~\ref{twocompcat}.

\subsubsection{Introduction rule}
We recall that the object \mbox{$\Gamma.A.A.\Id_A \in \Con$}, together with the
maps $s \defeq \pi_1 \pi_{\Id_A}$ and \mbox{$t
\defeq \pi_2 \pi_{\Id_A} \colon \Gamma.A.A.\Id_A \to \Gamma.A$},
is an arrow object for $\Gamma.A$ in $\Con / \Gamma$. As in
Proposition~\ref{arrowobj}, we write $\kappa \colon s \Rightarrow t$ for the
corresponding universal $2$-cell. Applying universality of $\kappa$ to the
$2$-cell
\begin{equation*}
  \cd[@!C@C-1em]{
    \Gamma.A \ar[dr]_{\pi_A} \ar@/^1em/[rr]^{\id} \ar@/_1em/[rr]_{\id} \dtwocell{rr}{\id} & &
    \Gamma.A \ar[dl]^{\pi_A} \\ &
    \Gamma
  }
\end{equation*}
in $\Con/\Gamma$, we obtain a morphism $r_A \colon \Gamma.A \to
\Gamma.A.A.\Id_A$ which factorises the diagonal: we have $\pi_{\Id_A} r_A =
\delta_A \colon \Gamma.A \to \Gamma.A.A$. We call this $r_A$ the \emph{semantic
introduction rule} for $\Id_A$.

\subsubsection{Elimination and computation rules}
With reference to Table~\ref{fig1}, we require semantic analogues of the
premisses $C$ and $d$ of the rule \textsc{$\Id$-elim}. These are given by an
object $C \in \Typ(\Gamma.A.A.\Id_A)$ and a map $d \colon \Gamma.A \to
\Gamma.A.A.\Id_A.C$ of $\Con$ making the following diagram commute:
\begin{equation}\label{idsquare}
    \cd{
        \Gamma.A \ar[r]^-{d} \ar[d]_{r_A} &
        \Gamma.A.A.\Id_A.C \ar[d]^{\pi_C} \\
        \Gamma.A.A.\Id_A \ar[r]_-{\id}  &
        \Gamma.A.A.\Id_A\text.
    }
\end{equation}
To give a semantic analogue of the conclusion $\J_d$, satisfying the analogue
of the computation rule, amounts to giving a filler $J_d \colon
\Gamma.A.A.\Id_A \to \Gamma.A.A.\Id_A.C$ making both sides of~\eqref{idsquare}
commute. Now, by Proposition~\ref{normalisofibration}, $\pi_C$ is a normal
isofibration in $\Con / \Gamma$; so that if we can show that $r_A$ is an
injective equivalence in $\Con / \Gamma$, then we may obtain the required
filler $J_d$ by an application of Proposition~\ref{elimrules1}. To show that
$r_A$ is an injective equivalence in $\Con / \Gamma$, we must first give a
retraction of $r_A$ over $\Gamma$. We take this to be $t \colon
\Gamma.A.A.\Id_A \to \Gamma.A$ (though we could equally well have chosen $s$);
and we have that $t r_A = \id_{\Gamma.A}$ as required. Next we need a $2$-cell
$\theta \colon \id \Rightarrow r_A t$ over $\Gamma$ satisfying $\theta r_A =
\id_{r_A}$ and $t \theta = \id_t$. For this, we consider the following diagram
of $1$- and $2$-cells $\Gamma.A.A.\Id_A \to \Gamma.A$:
\begin{equation*}
    \cd{
        s \ar@2[r]^-{\kappa} \ar@2[d]_\kappa & s r_A t \ar@2[d]^{\id_t} \\
        t \ar@2[r]_-{\id_t} & t r_A t\text.
    }
\end{equation*}
Because $t r_A = s r_A = \id_{\Gamma.A}$, this diagram is commutative: and so
by the two-dimensional aspect of the universal property of $\Gamma.A.A.\Id_A$,
is induced by a $2$-cell $\theta \colon \id \Rightarrow r_A t$ over $\Gamma$
satisfying $s \theta = \kappa$ and $t \theta = \id_t$. It remains to verify
that $\theta r_A = \id_{r_A}$. By the uniqueness part of the universal property
of $\Gamma.A.A.\Id_A$, it suffices to show that $\kappa \circ \theta r_A =
\kappa \circ \id_{r_A}$. But here we have $\kappa \theta = \kappa (r_A t) \circ
s \theta = \id_t \circ \kappa = \kappa$ and so $\kappa \circ \theta r_A =
\kappa r_A = \kappa \circ \id_{r_A}$ as required.

\subsubsection{Stability rules}\label{idSTAB}
We now verify that the semantic identity rules given above are stable under
semantic substitution. So suppose given $f \colon \Delta \to \Gamma$ in~$\Con$
together with $A \in \Typ(\Gamma)$. We must verify three things. First we must
show that reindexing $\Gamma.A.A.\Id_{A}$ along $f$ yields $\Delta.\sb f A.\sb
f A.\Id_{\sb f A}$. This follows immediately from the stability requirements of
Proposition~\ref{bcids}. Next, we must show that the semantic introduction rule
$r_{\sb f A}$ is the reindexing along $f$ of $r_A$. This follows from the fact
that arrow object structure on $\Id_{\sb f A}$ is the reindexing of that on
$\Id_A$ along $f$. Finally, we must show that applications of the semantic
elimination rule are stable under substitution. So suppose given a diagram
like~\eqref{idsquare}. If we view this as a diagram in $\Con / \Gamma$, then we
can reindex it along $f$ to yield a diagram
\begin{equation}\label{isdq}
    \cd{
        \Delta.f^\ast A \ar[r]^-{f^\ast d} \ar[d]_{r_{\sb f A}} &
        \Delta.\sb f A.\sb f A.\Id_{\sb f A}.Cf \ar[d]^{\pi_{Cf}} \\
        \Delta.\sb f A.\sb f A.\Id_{\sb f A} \ar[r]_-{\id}  &
        \Delta.\sb f A.\sb f A.\Id_{\sb f A}
    }
\end{equation}
in $\Con / \Delta$. We must show that pulling back the assigned filler
for~\eqref{idsquare} along $f$ yields the assigned filler for~\eqref{isdq}.
Now, by the stability properties of Proposition~\ref{stabilityiso}, we know
that the isofibration structure on $\pi_{Cf}$ qua map of $\Con / \Delta$ is the
one induced by pulling back along $f$ the isofibration structure of $\pi_C$ qua
map of $\Con / \Gamma$. Moreover, by the stability of the arrow object
structure of $\Id_A$, the injective equivalence structure on $r_{\sb f A}$ is
the one induced by pulling back that of $r_A$ along $f$. The result now follows
by applying Proposition~\ref{elimrules2}.

\subsubsection{Remark} Because $\Gamma.A.A.\Id_A$ is an arrow object in $\Con /
\Gamma.A$, we will in what follows pass back and forward without comment
between morphisms $h \colon \Lambda \to \Gamma.A.A.\Id_A$ and $2$-cells $\gamma
\colon sh \Rightarrow th \colon \Lambda \to \Gamma.A$ over $\Gamma$.

\subsubsection{Discrete identity rules}
We now show that the semantic identity rules given above satisfy the semantic
equivalents of the rules in Table~\ref{fig2}. So suppose given $\Gamma \in
\Con$ and $A \in \Typ(\Gamma)$ as before. The semantic analogues of the
premisses of the rules in Table~\ref{fig2} are a pair of morphisms $a, b \colon
\Gamma \to \Gamma.A$ of $\Con$ over $\Gamma$, together with a $2$-cell
\begin{equation*}
  \cd{
    \Gamma \ar@/^1em/[r]^{p} \ar@/_1em/[r]_{q} \dtwocell{r}{\alpha} &
    *+[r]{\Gamma.A.A.\Id_A}
  }
\end{equation*}
satisfying $sp = sq = a$, $s \alpha = \id_a$, $tp = tq = b$ and $t \alpha =
\id_b$. We must show that under these circumstances we have $p = q$ and $\alpha
= \id_p$. So consider the following diagram of $1$- and $2$-cells $\Gamma \to
\Gamma.A$:
\begin{equation*}
    \cd{
        sp \ar@2[r]^-{s\alpha} \ar@2[d]_{\kappa p} & sq \ar@2[d]^{\kappa q} \\
        tp \ar@2[r]_-{t \alpha} & tq\text.
    }
\end{equation*}
It is commutative, with both sides equal to $\kappa \alpha \colon sp
\Rightarrow tq$; but since $s\alpha = \id_a$ and $t \alpha = \id_b$, we deduce
that $\kappa \alpha = \kappa p = \kappa q \colon a \Rightarrow b$. By the
uniqueness part of the universal property of $\kappa$, this entails that $p = q
\colon \Gamma \to \Gamma.A.A.\Id_A$. Moreover, we have $\kappa \alpha = \kappa
p = \kappa \id_p$, and so again by the uniqueness part of the universal
property of $\kappa$, we deduce that $\alpha = \id_p$ as required.

\subsection{Unit types}\label{intunittypes}

\subsubsection{Formation rule}
Given $\Gamma \in \Con$, we define the \emph{semantic unit type} at $\Gamma$ to
be the object $\mathbf 1_\Gamma \in \Typ(\Gamma)$ whose existence is assured by
Definition~\ref{twocompcat}.

\subsubsection{Introduction rule}\label{unitintro}
Recall that $\mathbf 1_\Gamma$ is a retract bireflection of $\id_\Gamma \colon
\Gamma \to \Gamma$ along the $2$-functor $E(\Gamma) \colon \Typ(\Gamma) \to
\Con/\Gamma$; so in particular, we have a unit map $u_\Gamma \colon \Gamma \to
\Gamma.\mathbf 1_\Gamma$ over $\Gamma$, and we call this the \emph{semantic
introduction rule} for $\mathbf 1_\Gamma$.

\subsubsection{Elimination and computation rules}
Suppose given $C \in \Typ(\Gamma.\mathbf 1_\Gamma)$ and a map $d \colon \Gamma
\to \Gamma.\mathbf 1_\Gamma.C$ of $\Con$ fitting into a commutative diagram
\begin{equation*}
    \cd{
        \Gamma \ar[r]^-{d} \ar[d]_{u} &
        \Gamma.\mathbf 1_\Gamma.C \ar[d]^{\pi_C} \\
        \Gamma.\mathbf 1_\Gamma \ar[r]_-{\id}  &
        \Gamma.\mathbf 1_\Gamma\text.
    }
\end{equation*}
The semantic elimination rule will assign to this data a filler $U \colon
\Gamma.\mathbf 1_\Gamma \to \Gamma.\mathbf 1_\Gamma.C$ making both triangles
commute. Because $\pi_C$ is an isofibration in $\Con / \Gamma$, it suffices to
show that $u_\Gamma$ is a injective equivalence in $\Con / \Gamma$, since then
we obtain the desired filler by Proposition~\ref{elimrules1}. First we must
give a retraction for $u_\Gamma$ over $\Gamma$. We take this to be $k \defeq
\pi_{\mathbf 1_\Gamma} \colon \Gamma.\mathbf 1_\Gamma \to \Gamma$, which
satisfies $k u_\Gamma = \id_{\Gamma}$ as required. We now give a $2$-cell
$\theta \colon \id_{\Gamma.\mathbf 1_\Gamma} \Rightarrow u_\Gamma k$ satisfying
$\theta u_\Gamma = \id_{u_\Gamma}$ and $k \theta = \id_k$. By the
two-dimensional aspect of the universal property of $\mathbf 1_\Gamma$, every
$2$-cell $\alpha \colon \id_{\Gamma.\mathbf 1_\Gamma} \circ u_\Gamma
\Rightarrow u_\Gamma k \circ u_\Gamma$ is of the form $\bar \alpha \circ
u_\Gamma$ for a unique $2$-cell $\bar \alpha \colon \id_{\Gamma.\mathbf
1_\Gamma} \Rightarrow u_\Gamma k$. But because $\id_{\Gamma.\mathbf 1_\Gamma}
\circ u_\Gamma = \id_{u_\Gamma} = u_\Gamma \circ \id_\Gamma = u_\Gamma k
u_\Gamma$, we have in particular the $2$-cell $\theta
\defeq \overline{\id_{u_\Gamma}} \colon \id_{\Gamma.\mathbf 1_\Gamma} \Rightarrow
u_\Gamma k$; which by definition satisfies $\theta u_\Gamma = \id_{u_\Gamma}$.
That it also satisfies $k \theta = \id_k$ follows from the fact that $\theta$
is a $2$-cell of $\Con / \Gamma$.

\subsubsection{Stability rules}
We must show that the semantic unit rules are stable under semantic
substitution. This follows by an argument entirely analogous to that
of~\S\ref{idSTAB}, but using the stability properties of
Proposition~\ref{bcunits} rather than Proposition~\ref{bcids}.

\subsection{Sum types}\label{intsumtypes}
\subsubsection{Formation rule}
Given $\Gamma \in \Con$, $A \in \Typ(\Gamma)$ and $B \in \Typ(\Gamma.A)$, we
define the \emph{semantic sum type} of $A$ and $B$ to be the object
$\Sigma_A(B) \in \Typ(\Gamma)$ whose existence is assured by
Definition~\ref{twocompcat}.

\subsubsection{Introduction rule}
$\Sigma_A(B)$ is a retract bireflection of $B \in \Typ(\Gamma.A)$ along the
$2$-functor $\Typ(\pi_A) \colon \Typ(\Gamma) \to \Typ(\Gamma.A)$; and so, as in
Proposition~\ref{depsums2}, we obtain from the unit of this bireflection a map
$i \colon \Gamma.A.B \to \Gamma.\Sigma_A(B)$ of $\Con / \Gamma$. We declare
this map to be the \emph{semantic introduction rule} for $\Sigma_A(B)$.

\subsubsection{Elimination and computation rules} We consider $C \in \Typ(\Gamma.\Sigma_A(B))$ and a
map $d \colon \Gamma.A.B \to \Gamma.\Sigma_A(B).C$ of $\Con$ fitting into a
commutative diagram
\begin{equation*}
    \cd{
        \Gamma.A.B \ar[r]^-{d} \ar[d]_{i} &
        \Gamma.\Sigma_A(B).C \ar[d]^{\pi_C} \\
        \Gamma.\Sigma_A(B) \ar[r]_-{\id}  &
        \Gamma.\Sigma_A(B)\text.
    }
\end{equation*}
To give the semantic elimination rule satisfying the semantic computation rule
is now to give a filler $E \colon \Gamma.\Sigma_A(B) \to \Gamma.\Sigma_A(B).C$
making both triangles commute. We know that $\pi_C$ is an isofibration in $\Con
/ \Gamma$, whilst Definition~\ref{twocompcat} assures us that $i$ is an
injective equivalence in $\Con / \Gamma$: thus we obtain the desired filler by
applying Proposition~\ref{elimrules1}.

\subsubsection{Stability rules}
We must show that the semantic rules for dependent sums are stable under
semantic substitution. Again, this follows by an argument analogous to that
of~\S\ref{idSTAB}, this time using the stability properties of
Proposition~\ref{bcsums}.

\subsection{Product types}\label{intprodtypes}
Finally, we give semantic analogues in $\mathbb C$ of the rules for the product
types. As in the one-dimensional case, there is a slight mismatch here between
the syntax and the semantics. This means that, in addition to the right
biadjoints to weakening, we will also need to make use of the semantic unit
types of~\S\ref{intunittypes}. See~\cite[\S5.1--5.3]{Jacobs1993Comprehension}
for a fuller discussion of this point.

\subsubsection{Formation rule}

For $\Gamma \in \Con$, $A \in \Typ(\Gamma)$ and $B \in \Typ(\Gamma.A)$, we
define the \emph{semantic product type} of $A$ and $B$ to be the object
$\Pi_A(B) \in \Typ(\Gamma)$ whose existence is assured by
Definition~\ref{twocompcat}.

\subsubsection{Application rule} $\Pi_A(B)$ is a retract bicoreflection of $B \in
\Typ(\Gamma.A)$ along $\Delta_A \defeq \Typ(\pi_A) \colon \Typ(\Gamma) \to
\Typ(\Gamma.A)$. The counit of this bicoreflection is a morphism $\epsilon
\colon \Delta_A \Pi_A(B) \to B$ of $\Typ(\Gamma.A)$. We define the
\emph{semantic application rule} for $\Pi_A(B)$ to be the corresponding
morphism $\epsilon \colon \Gamma.A.\Pi_A(B) \to \Gamma.A.B$ of $\Con /
\Gamma.A$.

\subsubsection{Abstraction and $\beta$-rules}
For this, we suppose given, as in the premiss of the abstraction rule, a
morphism $f \colon \Gamma.A \to \Gamma.A.B$ over $\Gamma.A$. We are required to
produce from this a map $\lambda(f) \colon \Gamma \to \Gamma.\Pi_A(B)$ over
$\Gamma$; which, in order for the $\beta$-rule to hold, should satisfy
$\epsilon \circ \Delta_A(\lambda(f)) = f$. So consider the unit type $\mathbf
1_{\Gamma.A} \in \Typ(\Gamma.A)$. Applying its universal property
%Recall that $\mathbf 1_{\Gamma.A}$ is a retract bireflection of
%$\id_{\Gamma.A}$ along the $2$-functor $E(\Gamma.A) \colon \Typ(\Gamma.A) \to
%\Con/\Gamma.A$. As before, we write $\eta_{\Gamma.A} \colon \Gamma.A \to
%\Gamma.A.\mathbf 1_{\Gamma.A}$ for the unit of this bireflection. Applying the
%universal property of $\mathbf 1_{\Gamma.A}$
to $f \colon \Gamma.A \to \Gamma.A.B$ yields a morphism $\overline{f} \colon
\Gamma.A.\mathbf 1_{\Gamma.A} \to \Gamma.A.B$ over $\Gamma.A$ satisfying
$\overline{f} \circ u_{\Gamma.A} = f$. We can view $\overline{f}$ as a morphism
$\mathbf 1_{\Gamma.A} \to B$ of $\Typ(\Gamma.A)$; which by the stability of
unit types under substitution is equally well a morphism $\overline{f} \colon
\Delta_A \mathbf 1_{\Gamma} \to B$ of $\Typ(\Gamma.A)$. Applying the universal
property of $\Pi_A(B)$ to this, we obtain a morphism $\overline{\overline{f}}
\colon \mathbf 1_{\Gamma} \to \Pi_A(B)$ of $\Typ(\Gamma)$ satisfying
\mbox{$\epsilon \circ \Delta_A(\overline{\overline{f}}) = \overline{f}$}. This
is equally well a morphism $\Gamma.\mathbf 1_{\Gamma} \to \Gamma.\Pi_A(B)$ over
$\Gamma$, so we can now define the map $\lambda(f) \colon \Gamma \to
\Gamma.\Pi_A(B)$ over $\Gamma$ to be $\lambda(f)
\defeq \overline{\overline{f}} \circ u_{\Gamma}$. It remains to show that we have $\epsilon \circ \Delta_A(\lambda(f))
= f$; for which we calculate that
\begin{equation*}
\epsilon \circ \Delta_A(\lambda(f)) = \epsilon \circ
\big(\Delta_A(\overline{\overline{f}}) \circ \Delta_A(u_\Gamma)\big) = \overline{f} \circ u_{\Gamma.A}
= f
\end{equation*}
as required. Here we have used the fact that, by stability of unit types under
substitution, we have $\Delta_A(u_\Gamma) = u_{\Gamma.A}$.

\subsubsection{Function extensionality rules}
We now give semantic analogues of the rules of Table~\ref{fig3}. For the first
rule $\Pi$-\textsc{ext}, we suppose given morphisms $m, n \colon \Gamma \to
\Gamma.\Pi_A(B)$ over $\Gamma$, together with a $2$-cell
\begin{equation*}
    p \colon \epsilon \circ \Delta_A(m) \Rightarrow \epsilon \circ \Delta_A(n) \colon \Gamma.A \to \Gamma.A.B
\end{equation*}
over $\Gamma.A$. We must produce from this a $2$-cell $\textsf{ext}(p) \colon m
\Rightarrow n$. First we apply the universal property of the unit type $\mathbf
1_{\Gamma}$ to $m$ and $n$ to obtain morphisms $\overline m, \overline n \colon
\Gamma.\mathbf 1_\Gamma \to \Gamma.\Pi_A(B)$ over $\Gamma$. These satisfy $m =
\overline{m} \circ u_\Gamma$ and $n = \overline{n} \circ u_\Gamma$, and so we
can view $p$ as a $2$-cell
\begin{equation*}
    p \colon \epsilon \circ \Delta_A(\overline{m}) \circ u_{\Gamma.A} \Rightarrow \epsilon \circ \Delta_A(\overline{n}) \circ u_{\Gamma.A} \colon \Gamma.A \to \Gamma.A.B\text,
\end{equation*}
where again we use stability of unit types under pullback to derive that
$\Delta_A(u_\Gamma) = u_{\Gamma.A}$. By the two-dimensional aspect of the
universal property of $\mathbf 1_{\Gamma.A}$, we have $p = \overline p \circ
u_{\Gamma.A}$ for a unique $2$-cell
\begin{equation*}
\overline p \colon \epsilon \circ \Delta_A(\overline m)
\Rightarrow \epsilon \circ \Delta_A(\overline n) \colon \Gamma.A.\mathbf 1_{\Gamma.A}
\to \Gamma.A.B\text.
\end{equation*}
Now, by the two-dimensional aspect of the universal property of $\Pi_A(B)$, we
have that $\overline p = \epsilon \circ \Delta_A(\overline{\overline p})$ for a
unique $\overline{\overline p} \colon \overline m \Rightarrow \overline n$. We
now define the $2$-cell $\mathsf{ext}(p)$ to be given by $\overline{\overline
p} \circ u_\Gamma \colon m \Rightarrow n$.

In order for $\mathsf{ext}$ to satisfy the analogue of the rule
\textsc{$\Pi$-ext-comp}, we must show that when $m = n$ and $p = \id_{\epsilon
\circ \Delta_A(m)}$, we have $\mathsf{ext}(p) = \id_m$. It suffices for this to
show that (with the above notation) $\overline{\overline p } = \id_{\overline
m} \colon \overline m \Rightarrow \overline m$; which, by applying successively
the universal properties of $\Pi_A(B)$ and $\mathbf 1_{\Gamma.A}$, follows from
the fact that $\epsilon \circ \Delta_A(\overline{\overline p}) \circ
u_{\Gamma.A} = p$ is an identity $2$-cell. Finally, we must verify that
$\mathsf{ext}$ satisfies the analogue of the rule \textsc{$\Pi$-ext-app}.
Recall from~\S\ref{depprodssec} that the operation $\ast$ appearing in
\textsc{$\Pi$-ext-app} is simply the lifting of $\epsilon \colon
\Gamma.A.\Pi_A(B) \to \Gamma.A.B$ to identity types. From this it follows that
we must verify that $\epsilon \circ \Delta_A(\mathsf{ext}(p)) = p \colon
\epsilon \circ \Delta_A(m) \Rightarrow \epsilon \circ \Delta_A(n)$. We
calculate that $\epsilon \circ \Delta_A(\mathsf{ext}(p)) = \epsilon \circ
\big(\Delta_A(\overline{\overline p}) \circ \Delta_A(u_\Gamma)\big)
     = \overline p \circ u_{\Gamma.A} = p$
as required.

\subsubsection{Stability rules}
We must now show that the semantic rules for dependent products are stable
under semantic substitution. This follows by an argument analogous to that
of~\S\ref{idSTAB}; though this time we do not need the stability properties of
isofibrations (Proposition~\ref{stabilityiso}) at all; instead, we need those
for products (Proposition~\ref{bcprods}) and also those for units
(Proposition~\ref{bcunits}).

\subsection{The internal language}\label{internalangg} We now define the type theory $\Ss(\mathbb C)$ associated
to our two-dimensional model $\mathbb C$. It is obtained by recursively
extending \mlt\ with additional inference rules. These inference rules are
``axiom'' rules with no premisses, and so may be specified by giving only their
conclusion. First we have rules introducing new types:
\begin{itemize}
\item For each $A \in \Typ(1)$ we add a judgement $\t \overline A \ \ty$.
\item For each $A \in \Typ(1)$, $B \in \Typ(1.A)$ we add a judgement
    \mbox{$x : \overline A \t \overline B (x) \ \ty$}.
\item And so on.
\end{itemize}
Then we have rules introducing new terms:
\begin{itemize}
\item For each $A \in \Typ(1)$, $a \in_1 A$, we add a judgement $\t
    \overline a : \overline A$.
\item For each $A \in \Typ(1)$, $B \in \Typ(1.A)$, $b \in_{1.A} B$, we add
    a judgement $x : \overline A \t \overline b(x) : \overline B (x)$.
\item And so on.
\end{itemize}
Here, we use the convention for global sections developed in
Notation~\ref{notconv}. Next we have rules identifying the syntactic notions of
substitution, weakening, contraction and exchange with their semantic
counterparts in $\mathbb C$. We give the case of substitution as a
representative sample. First we deal with substitution in types:
\begin{itemize}
\item For each $A \in \Typ(1)$, $B \in \Typ(1.A)$, $a \in_1 1.A$, we add a
    judgement \mbox{$\t \overline B (\overline a) = \overline{\sb a B} \
    \ty$}.
\item For each $A \in \Typ(1)$, $B \in \Typ(1.A)$, $C \in \Typ(1.A.B)$, $b
    \in_{1.A} B$, we add a judgement \mbox{$x : \overline A \t \overline C
    (x, \overline b(x)) = \overline{\sb b C} \ \ty$}.
\item For each $A \in \Typ(1)$, $B \in \Typ(1.A)$, $C \in \Typ(1.A.B)$ and
    $a \in_1 A$, we add a judgement \mbox{$y : \overline B (\overline a) \t
    \overline C (\overline a, y) = \overline{\sb {(a.\ro a)} C} \ \ty$}.
%\item For each $A \in \Typ(1)$, $B \in \Typ(1.A)$ and $C \in \Typ(1.A.B)$,
%    we adjoin a judgement $x : A \c y : B \t \overline C (x, y) \ \ty$.
\item And so on.
\end{itemize}
And now substitution in terms:
\begin{itemize}
\item For each $A \in \Typ(1)$, $B \in \Typ(1.A)$, $a \in_1 A$, $b
    \in_{1.A} B$, we add  \mbox{$\t \overline b(\overline a) =
    \overline{\sb a b} : \overline B (\overline a)$}.
\item For each $A \in \Typ(1)$, $B \in \Typ(1.A)$, $C \in \Typ(1.A.B)$, $b
    \in_{1.A} B$, $c \in_{1.A.B} C$, we add \mbox{$x : \overline A \t
    \overline c (x, \overline b(x)) = \overline{\sb b c} : \overline C(x,
    \overline c(x))$}.
\item For each $A \in \Typ(1)$, $B \in \Typ(1.A)$, $C \in \Typ(1.A.B)$, $a
    \in_{1} A$, $c \in_{1.A.B} C$, we add \mbox{$y : \overline B (\overline
    a) \t \overline c (\overline a, y) = \overline{\sb {(a.\ro a)} c} :
    \overline C (\overline a, y)$}.
\item And so on.
\end{itemize}
Finally, we have rules identifying each of the logical rules of $\mlt$ with its
semantic counterpart in $\mathbb C$. We give only the case of the identity
types; the remainder follow the same pattern. First we have the formation
rules.
\begin{itemize}
\item For each $A \in \Typ(1)$, we add \mbox{$x, y : \overline A \t
    \Id_{\overline A}(x, y) = \overline{\Id_A}(x, y) \ \ty$}.
\item For each $A \in \Typ(1)$, $B \in \Typ(1.A)$, we add
\begin{equation*}
x :
    \overline A\c y, z : \overline B(x) \t \Id_{\overline B(x)}(y, z) =
    \overline{\Id_B}(x, y, z) \ \ty\text.
\end{equation*}
\item And so on.
\end{itemize}
Next we have the introduction rule. We observe that for $A \in \Typ(\Gamma)$,
the semantic introduction rule $r_A \colon \Gamma \to \Gamma.A.A.\Id_A$ over
$\Gamma$ can be viewed as a global section $r_A \in_\Gamma
\delta_A^\ast(\Id_A)$, where $\delta_A \colon \Gamma.A \to \Gamma.A.A$ is the
diagonal morphism. Thus we may add the following rules:
\begin{itemize}
\item For each $A \in \Typ(1)$, we add \mbox{$x : \overline A \t \r(x) =
    \overline {r_A} (x) : \Id_{\overline A}(x, x)$}.
\item For each $A \in \Typ(1)$, $B \in \Typ(1.A)$, we add
\begin{equation*}
x :
    \overline A \c y : \overline B(x) \t \r(y) = \overline {r_B} (x, y) :
    \Id_{\overline B(x)}(y, y)\text.\end{equation*}
\item And so on.
\end{itemize}
Finally we come to the identity elimination rule.
\begin{itemize}
\item For each $A \in \Typ(1)$, $C \in \Typ(1.A.A.\Id_A)$ and $d \colon 1.A
    \to 1.A.A.\Id_A.C$ as in~\eqref{idsquare} (which is equally well a
    global section $d \in_{1.A} r_A^\ast C$), we add
\begin{equation*}
x, y : \overline A \c p : \Id_{\overline A}(x, y) \t \J_{\overline d}(x, y, p) = \overline {J_d} (x, y, p) : \overline C(x, y, p)\text.
\end{equation*}

\item And so on.
\end{itemize}

Now, in order for the internal language we have set up to be of any use, we
require its types and terms to denote unique elements of the model $\mathbb C$.
The next Proposition tells us that this is the case.

\begin{Prop}[Soundness]\label{soundness}
For any $B, C \in \Typ(1.A_1.A_2\dots A_n)$, if the judgement
\begin{equation*}
    x : \overline A_1.\overline A_2 \dots \overline A_n \t \overline B(x) = \overline C(x)\ \ty
\end{equation*}
is derivable, then $B = C$. Likewise, for global sections $b, c \in_{1.A_1\dots
A_n} B$, if the judgement
\begin{equation*}
    x : \overline A_1.\overline A_2 \dots \overline A_n \t \overline b(x) = \overline c(x) : \overline B(x)
\end{equation*}
is derivable, then $b = c$.
\end{Prop}
\begin{proof}
By induction on the derivation of the judgement in question, it suffices to
show that the semantic counterpart of each syntactic equality rules is
satisfied. For the non-logical equality rules, this is standard (though
delicate), and we refer the reader to~\cite{Hofmann1995Extensional}
or~\cite{Pitts2000Categorical} for the details (note that we make essential
use of the fact that the underlying $1$-fibration of $\Typ \to \Con$ is split).
The other cases we must consider are the computation rules of
Tables~\ref{fig1},~\ref{fig2}~or~\ref{fig3}, and the rules expressing stability
of the logical operations under substitution; and each of these has been dealt
with in the preceding sections.
\end{proof}
\begin{Rk}\label{aproblem}
Observe that the internal language $\Ss(\mathbb C)$ does not give us access to
all of the model $\mathbb C$: it only allows us to talk about objects of the
base $2$-category $\Con$ which have the form $1.A_0\dots A_n$ (where $1$ is the
given $2$-terminal object). There are two ways around this. We can modify the
syntax of our type theory so that contexts and context morphisms are primitive,
rather than derived, notions. Then each object or morphism of $\Con$
corresponds directly to a context or context morphism of $\Ss(\mathbb C)$.
Alternatively, we can keep our type theory the same, and instead work with
relative internal languages. Given $\Gamma \in \Con$, the \emph{relative
internal language} $\Ss_\Gamma(\mathbb C)$ is the type theory whose closed
types are objects of $\Typ(\Gamma)$, with dependent types being objects of
$\Typ(\Gamma.A)$, $\Typ(\Gamma.A.B)$ and so on. Moreover, because each morphism
$\Gamma \to \Delta$ of $\Con$ induces an interpretation (in the sense
of~\S\ref{functorialaspsects} below) $\Ss_\Delta(\mathbb C) \to
\Ss_\Gamma(\mathbb C)$, we obtain what is in an obvious sense a
``$\Con$-indexed type theory''\footnote{A finer analysis shows that this is
really a two-dimensional indexing. That is, we have a trihomomorphism
$\Con^{\coop} \to \cat{Th}$, where $\cat{Th}$ is a suitably-defined tricategory
of two-dimensional theories.}.
\end{Rk}

\subsection{Functorial aspects}\label{functorialaspsects} In Sections~\ref{sec2} and~\ref{sec3}, we
constructed from each type theory $\Ss$ incorporating $\mlt$ a two-dimensional
model $\mathbb C(\Ss)$; whilst in the preceding parts of the present Section,
we have constructed from each two-dimensional model $\mathbb C$ a type theory
$\Ss(\mathbb C)$ incorporating $\mlt$. It is natural to ask whether these
assignations give rise to a \emph{functorial semantics} in the spirit
of~\cite{Lawvere1968Some}. That is, can we define a syntactic category of type
theories and a semantic category of models for which the above assignations
underlie an equivalence of categories? We finish the paper by sketching a
positive answer to this question.

We first define a syntactic category $\cat{Th}$. Its objects are the
\emph{generalised algebraic theories}~\cite{Cartmell1986Generalised} over
$\mlt$. These are defined inductively by the following three clauses. Each
object of $\cat{Th}$ is a sequent calculus; \mbox{$\mlt \in \cat{Th}$}; and if
$\Ss \in \cat{Th}$, then so is any extension of $\Ss$. Here, an
\emph{extension} of $\Ss$ is given by adjoining a set of inference rules each
of which has no premisses and a conclusion $\mathcal J$ which obeys the
following requirements. If $\mathcal J$ is of the form $\Gamma \vdash A \ \ty$
then $A$ must be fresh for $\Ss$ and $\Gamma$ must be a well-formed context of
$\Ss$; if $\mathcal J$ is of the form $\Gamma \vdash a : A$ then $a$ must be
fresh for $\Ss$ and $\Gamma \vdash A \ \ty$ must be derivable in $\Ss$; if
$\mathcal J$ is of the form $\Gamma \vdash A = B \ \ty$ then $\Gamma \vdash A \
\ty$ and $\Gamma \vdash B \ \ty$ must be derivable in $\Ss$; and finally if
$\mathcal J$ is of the form $\Gamma \vdash a = b : A$ then $\Gamma \vdash a :
A$ and $\Gamma \vdash b : A$ must be derivable in $\Ss$. Note that the
assignation $\mathbb C \mapsto \Ss(\mathbb C)$ sends each two-dimensional model
to a GAT over $\mlt$.

The morphisms of $\cat{Th}$ are equivalence classes of interpretations. Given
$\Ss, \Tt \in \cat{Th}$, an \emph{interpretation} $F \colon \Ss \to \Tt$ is a
function $F$ taking derivable judgements of $\Ss$ to derivable judgements of
$\Tt$, subject to the following requirements. Each $F(A \ \ty)$ should have the
form $FA \ \ty$; each $F(a : A)$ should have the form $Fa : FA$; each $F(A = B
\ \ty)$ should have the form $FA = FB \ \ty$; and each $F(a = b : A)$ should
have the form $Fa = Fb : FA$. Moreover, if we suppose $F(\Gamma \vdash A \
\ty)$ has the form $F\Gamma \vdash FA \ \ty$, then each \mbox{$F(\Gamma \c x :
A \vdash B(x) \ \ty)$} should have the form $F\Gamma \c x : FA \vdash FB(x) \
\ty$; each $F(\Gamma \c x : A \vdash b(x) : B(x))$ should have the form
$F\Gamma \c x : FA \vdash Fb(x) : FB(x)$; and similarly for the two equality
judgement forms. Finally we require that $F$ should commute with all the
inference rules of $\mlt$. We give the case of the rule of $\Id$-formation for
illustration. Suppose given a derivable judgement $\Gamma \vdash A \ \ty$ in
$\Ss$. We write its image under $F$ as $F\Gamma \vdash FA \ \ty$, and the image
of $\Gamma \c x, y : A \vdash \Id_A(x, y)\ \ty$ as $F\Gamma \c x, y : FA \vdash
F\Id_A(x,y) \ \ty$. Now the following judgement should be derivable in $\Tt$:
\begin{equation*}
    F\Gamma \c x, y : FA \vdash \Id_{FA} (x, y) = F\Id_A(x,y) \ \ty\text.
\end{equation*}
The equivalence relation we impose on interpretations identifies $F, G \colon
\Ss \to \Tt$ if they differ only up to definitional equality in the obvious
sense. It is now straightforward to show that GATs and equivalence classes of
interpretations form a category $\cat{Th}$.
\begin{Rk}
Using the above notion of interpretation, we can now say what it means to give
an interpretation of a GAT $\Tt$ in a two-dimensional model $\mathbb C$:
namely, to give an interpretation (in the above sense) $\Tt \to \Ss(\mathbb
C)$. It is easy to check that this accords with the intuitive syntactic notion
we would give.
\end{Rk}

We now define a semantic category $\cat{Mod}$. Its objects are models of
two-dimensional type theory as in Definition~\ref{twocompcat}. A morphism $F
\colon \mathbb C \to \mathbb C'$ is given by a pair of $2$-functors $F_1 \colon
\Con \to \Con'$ and $F_2 \colon \Typ \to \Typ'$ rendering commutative the
squares
\begin{equation*}
    \cd{
      \Typ \ar[r]^{F_2} \ar[d]_{p} &
      \Typ' \ar[d]^{p'} \\
      \Con \ar[r]_{F_1} &
      \Con'
    } \qquad \text{and} \qquad
    \cd{
      \Typ \ar[r]^{F_2} \ar[d]_{E} &
      \Typ' \ar[d]^{E'} \\
      \Con^\mathbf 2 \ar[r]_{{F_1}^\mathbf 2} &
      (\Con')^\mathbf 2
    }
\end{equation*}
and preserving all the additional structure on the nose.
\begin{Prop}\label{givesfunctors}
The assignations $\Ss \mapsto \mathbb C(\Ss)$ and $\mathbb C \mapsto
\Ss(\mathbb C)$ underlie functors $\mathbb C(\thg) \colon \cat{Th} \to
\cat{Mod}$ and $\Ss(\thg) \colon \cat{Mod} \to \cat{Th}$.
\end{Prop}
\begin{proof}
Given an interpretation $F \colon \Ss \to \Tt$, we define functors $G_0 \colon
\Con(\Ss) \to \Con(\Tt)$ and $G_1 \colon \Typ(\Ss) \to \Typ(\Tt)$ by an obvious
structural induction over the objects and morphisms of the domain categories.
In order to extend these functors to $2$-functors, we first show by induction
that every object $\Gamma \in \Con(\Ss)$ has an accompanying arrow object given
by the identity context $\Id_\Gamma$. But now, since $F$ preserves the identity
type structure, the corresponding $G_0$ will preserve these arrow objects;
hence we may extend $G_0$ to a \mbox{$2$-functor} by regarding each $2$-cell of
$\Con(\Ss)$ as a $1$-cell into an arrow object, mapping this $1$-cell over and
then turning the resulting $1$-cell back into a $2$-cell of $\Con(\Tt)$.
Because the comprehension $2$-functors of $\mathbb C(\Ss)$ and $\mathbb C(\Tt)$
are $2$-fully faithful, this in turn determines the extension of $G_1$ to a
$2$-functor. Finally, the fact that $F$ strictly preserves the remaining
structure implies that the same is true of $(G_0, G_1)$, and so we obtain a
morphism of models $\mathbb C(\Ss) \to \mathbb C(\Tt)$ as required.

Conversely, given a morphism of models $F \colon \mathbb C \to \mathbb C'$, we
may define an interpretation $\Ss(\mathbb C) \to \Ss(\mathbb C')$ as follows.
By structural induction, every closed type $A$ of $\Ss(\mathbb C)$ is
definitionally equal to one of the form $\overline X$ for some $X \in \Typ(1)$;
and by Proposition~\ref{soundness}, this $X$ is unique. Thus we may define the
image of $A$ under the interpretation to be the type $\overline{G_1(X)}$ of
$\Ss(\mathbb C')$. Similarly, every closed term $a : A$ of $\Ss(\mathbb C)$ is
definitionally equal to one of the form $\overline x : \overline X$ for a
unique map $x \colon 1 \to 1.X$ of $\Con$, and so we may define $F(a)$ to be
the term $\overline{G_0(x)} : \overline{G_1(X)}$. This definition extends to
types and terms in non-empty contexts in an obvious way. Finally, the fact that
our morphism of models preserves all the remaining structure on the nose
implies the same for the interpretation just described.
\end{proof}

However, the functors defined in this Proposition do not give rise to an
equivalence of categories. There are two reasons for this. The first is the
issue raised in Remark~\ref{aproblem}. Observe that any two-dimensional model
in the image of $\mathbb C(\thg)$ has the property that each object $\Gamma \in
\Con$ is of the form $1.A_1\dots A_n$ for a unique (possibly empty) sequence of
objects $A_1 \in \Typ(1), \dots, A_n \in \Typ(1.A_1\dots A_{n-1})$. This is the
``tree condition'' of~\cite{Cartmell1986Generalised}. Clearly not every
two-dimensional model has this property, so that if we are to obtain an
equivalence, we must first cut down to the full sub-$2$-category
$\cat{Mod}_{tr} \subset \cat{Mod}$ on those which do. The second reason we do
not obtain an equivalence is more subtle. In order for $\cat{Th} \simeq
\cat{Mod}_{tr}$ to hold, we must certainly have for each $\Ss \in \cat{Th}$
that $\Ss(\mathbb C(\Ss)) \cong \Ss$. However, this turns out not to be the
case: we run into problems with the terms witnessing the elimination rules. As
an illustration, we will show that $\mlt \ncong \Ss(\mathbb C(\mlt))$. Because
the object $\mlt$ is initial in $\cat{Th}$, there is a unique morphism \mbox{$F
\colon \mlt \to \Ss(\mathbb C(\mlt))$}: and so it suffices to show that $F$ is
not surjective. First observe that by \mbox{$\mathbf 1$-elimination} we can
derive a judgement
\begin{equation}\label{derivablej}
    z : \mathbf 1 \t \mathrm U_\star(z) : \mathbf 1
\end{equation}
in $\mlt$. Next note that the judgements of $\Ss(\mathbb C(\mlt))$ are simply
equivalence classes of judgements of $\mlt$ with respect to definitional
equality; and so by passing to the quotient, we obtain from~\eqref{derivablej}
a judgement
\begin{equation}\label{derivablej2}
    z : [\mathbf 1] \t [\mathrm U_\star](z) : [\mathbf 1]
\end{equation}
of $\Ss(\mathbb C(\mlt))$. The crucial point is that~\eqref{derivablej2} does
not coincide with the value of $F$ at the judgement~\eqref{derivablej}. This
latter can be described as follows. First we derive a term $z : \mathbf 1 \t
\phi(z) : \Id_\mathbf 1(z, \star)$ in $\mlt$ by $\mathbf 1$-elimination, taking
$\phi(z)
\defeq \mathrm U_{\r(\star)}(z)$. Now by the description of the semantic
unit types given in~\S\ref{intunittypes}, we see that applying $F$
to~\eqref{derivablej} yields (up to definitional equality) the following
judgement in $\Ss(\mathbb C(\mlt))$:
\begin{equation}
\label{newrule}
    z : [\mathbf 1] \t [\phi(z)^\ast(\star)] : [\mathbf 1]\text.
\end{equation}
Now, if \eqref{newrule} were definitionally equal to~\eqref{derivablej2}, then
we would also have that $z : \mathbf 1 \t \mathrm U_\star (z) =
\phi(z)^\ast(\star) : \mathbf 1$ in $\mlt$, and this is not the case. Hence $F$
applied to~\eqref{derivablej} does not yield~\eqref{derivablej2}, from which it
follows by induction over derivable judgements of $\mlt$ that
\eqref{derivablej2} cannot lie in the image of $F \colon \mlt \to \Ss(\mathbb
C(\mlt))$.

%On the other hand, this judgement is not derivable in \mlt\ itself; and so
%there can be no interpretation $F \colon \Ss(\mathbb C(\mlt)) \to \mlt$, and in
%particular no isomorphism.

There are several ways in which we can resolve this issue. The first is for us
to change our notion of model so that it accords more closely with the type
theory. This is unsatisfactory, as we have then reverted to a categorical
paraphrasing of type theoretic syntax.  A second alternative is to change our
notion of type theory so that it accords more closely with the categorical
model. This involves removing the elimination rules altogether, and instead
taking the Leibniz rule, together with the injective equivalence structures on
the introduction terms, as primitives. This is unsatisfactory for a more subtle
reason. Whilst it may be reasonably straightforward to give this alternative
presentation for two-dimensional type theory, we would find as we moved towards
full intensional type theory that it would require a more and more intricate
set of rules expressing appropriate coherence properties of our new primitives.
The elegant simplicity of intensional type theory would be lost completely.

A third solution, and our preferred one, is to equip our categories of theories
and of models with more generous notions of morphism, ones which preserve some
of the structure only up to propositional, rather than definitional equality.
There is a great deal of scope in how far we go with this. In the present
paper, we make only the minimal modifications necessary to obtain the desired
equivalence. A fuller treatment would take account of the fact that our models
and theories are themselves two-dimensional structures, so that their
respective totalities should give rise not merely to equivalent categories, but
also to triequivalent $\cat{Gray}$-categories (=semi-strict $3$-categories) in
the sense of~\cite{Gordon1995Coherence}. Adopting this more comprehensive
approach would be necessary if, for instance, we wished to study the 2-category
of interpretations of some generalised algebraic theory inside a particular
two-dimensional model. However, for our present purposes, we do not need to go
this far; and so, in the interests of brevity, we do not.

The minimal modification that we will consider is given as follows. On the
syntactic side, we define a category $\cat{Th}_\psi$ with as objects GATs over
$\mlt$ and as maps \mbox{$F \colon \Ss \to \Tt$} \emph{pseudo-interpretations},
whose definition generalises that of an interpretation by dropping the
requirement that $F$ should preserve each of the following rules up to
definitional equality: \textsc{$\mathbf 1$-elim}, \textsc{$\mathbf \Id$-elim},
\textsc{$\Sigma$-elim}, and \mbox{\textsc{$\Pi$-abs}}. One may now think that,
in order to justify the name pseudo-interpretation, we should ask for $F$ to
preserve these rules at least up to propositional equality; but it turns out
that this is unnecessary, because this weaker form of preservation is a
consequence of the type-theoretic elimination rules.

On the semantic side we define a category $\cat{Mod}_\psi$ whose objects are
two-dimensional models, and whose maps $F \colon \mathbb C \to \mathbb C'$ are
\emph{pseudo-morphisms}. These are obtained by relaxing in the definition of
morphism of models the requirement that the following structure should be
preserved: the normal isofibration structures on dependent projections $\pi_A$;
the injective equivalence structures on the maps $i \colon \Gamma.A.B \to
\Gamma.\Sigma_A(B)$ associated to dependent sums; and the assignations $f
\mapsto \overline f$ on $1$-cells associated to the unit types, dependent sums,
and dependent products. Once again, we do not need to add conditions requiring
these pieces of structure to be preserved up to isomorphism, since this will be
an automatic consequence of the remaining structure. As before, we write
$(\cat{Mod}_\psi)_{tr}$ for the full subcategory of $\cat{Mod}_\psi$ on those
models satisfying the tree condition.
\begin{Prop}
The functors $\mathbb C(\thg)$ and $\Ss(\thg)$ extend to functors
$\cat{Th}_\psi \to (\cat{Mod}_\psi)_{tr}$ and $(\cat{Mod}_\psi)_{tr} \to
\cat{Th}_\psi$ respectively.
\end{Prop}
\begin{proof}
The argument of Proposition~\ref{givesfunctors} carries over almost entirely
unmodified. The only subtlety arises in defining the pseudo-morphism of models
$\mathbb C(\Ss) \to \mathbb C(\Tt)$ corresponding to a pseudo-interpretation $F
\colon \Ss \to \Tt$. As before, we define functors $G_0$ and $G_1$ by induction
on the objects and morphisms of $\mathbb C(\Ss)$; but when it comes to
extending these to $2$-functors, we encounter the problem that the
interpretation $F$, since it no longer preserves the rule \textsc{$\mathbf
\Id$-elim}, may not send identity contexts to identity contexts. However, using
the fact that \textsc{$\mathbf \Id$-elim} is preserved at least up to
\emph{propositional} equality, we may show by induction that $F$ will send an
identity context to something isomorphic to an identity context. From this, it
follows that $G_0$ will still preserve arrow objects, so that we may continue
the argument as before.
\end{proof}

But now we have that:
\begin{Prop}
The functors $\mathbb C(\thg)$ and $\Ss(\thg)$ induce an equivalence of
categories $(\cat{Mod}_\psi)_{tr} \simeq \cat{Th}_\psi$.
\end{Prop}
\begin{proof}
First observe that if $\mathbb C$ is a model satisfying the tree condition,
then the contexts and context morphisms of $\Ss(\mathbb C)$ are, up to
definitional equality, just the objects and morphisms of $\Con$, whilst the
types-in-context of $\Ss(\mathbb C)$ are just the objects of $\Typ$. From this
it follows that $\Ss(\thg)$ is fully faithful. Indeed, given a
pseudo-interpretation $F \colon \Ss(\mathbb C) \to \Ss(\mathbb C')$, our
observation allows us to define functors $G_0 \colon \Con \to \Con'$ and $G_1
\colon \Typ \to \Typ'$. As in the last proof, the pseudo-interpretation $F$
must send an identity context to something isomorphic to an identity context;
from which it follows that $G_0$ preserves arrow objects, allowing us to extend
$G_0$ and $G_1$ to $2$-functors as before. It is now easy to verify that the
resultant pair $(G_0, G_1)$ is a pseudo-morphism; and by examining the proof of
Proposition~\ref{givesfunctors}, we see that it is sent by $\Ss(\thg)$ to $F$
and that it is the unique pseudo-morphism with this property.
%
% By structural induction, we may show that
%the contexts and
%
% context of $\Ss(\mathbb C')$ of length $n$ will be definitionally equal to
%one of the form $1.\overline{B_1} \dots \overline{B_n}$; and by
%Proposition~\ref{soundness}, this context will be unique. Hence we can define a
%pseudo-morphism $\tilde F \colon \mathbb C \to \mathbb C'$ on objects by
%sending the context $1.A_1 \dots A_n$ to the unique context $1.B_1 \dots B_n$
%of $\mathbb C'$ such that $F(1.\overline{A_1} \dots \overline{A_n}) =
%1.\overline{B_1} \dots \overline{B_n}$ (note that here we use the fact that
%$\mathbb C$ satisfies the tree condition). The rest of the definition of
%$\tilde F$ proceeds similarly.

It remains to show that for each $\Ss \in \cat{Th}_\psi$ we have an isomorphism
$\Ss \cong \Ss(\mathbb C(\Ss))$ in $\cat{Th}_\psi$, and that these are natural
in $\Ss$. Now, up to definitional equality, the judgements of $\Ss(\mathbb
C(\Ss))$ are the same as those of $\Ss$; and so we obtain mutually inverse
assignations between the judgements of the former and those of the latter.
Moreover, by following through the constructions of $\mathbb C(\thg)$ and
$\Ss(\thg)$, we see that all of the logical structure of $\Ss(\mathbb C(\Ss))$
is given as in $\Ss$, with the possible exception of the rules \textsc{$\mathbf
1$-elim}, \textsc{$\mathbf \Id$-elim}, \textsc{$\Sigma$-elim}, and
\mbox{\textsc{$\Pi$-abs}} (as seen in the discussion following
Proposition~\ref{givesfunctors}). But this says precisely that these mutually
inverse assignations are pseudo-interpretations, and so give rise to a natural
isomorphism $\Ss \cong \Ss(\mathbb C(\Ss))$ in $\cat{Th}_\psi$ as required.
\end{proof}
\bibliography{rhgg2}
\end{document}